\documentclass[a4paper, 11pt, leqno]{amsart}

\usepackage{amsmath, amssymb, amsfonts, amsthm, mathrsfs}
\usepackage{graphicx} 
\usepackage{mathtools}
\usepackage[dvipsnames]{xcolor}
\usepackage{comment}
\usepackage{enumitem}

\usepackage[T1]{fontenc}
\usepackage[margin=1in]{geometry}
\usepackage[
  pagebackref,
  colorlinks,
  citecolor=Mahogany,
  linkcolor=Mahogany,
  urlcolor=Mahogany,
  filecolor=Mahogany
]{hyperref}
\usepackage{color}

\usepackage[final,notref,notcite,color]{showkeys}
\definecolor{labelkey}{rgb}{0,0,1}
\usepackage{aliascnt}
\usepackage[capitalize]{cleveref}

\theoremstyle{plain}
\newtheorem{thm}{Theorem}[section]
\newaliascnt{prop}{thm}
\newtheorem{prop}[prop]{Proposition}
\aliascntresetthe{prop}
\newaliascnt{lem}{thm}
\newtheorem{lem}[lem]{Lemma}
\aliascntresetthe{lem}
\newaliascnt{cor}{thm}
\newtheorem{cor}[cor]{Corollary}
\aliascntresetthe{cor}

\theoremstyle{definition}
\newaliascnt{dfn}{thm}
\newtheorem{dfn}[dfn]{Definition}
\aliascntresetthe{dfn}
\newaliascnt{assumption}{thm}

\aliascntresetthe{assumption}
\let\TeXPrimitiveDef\def
\newcommand{\BeginDefEnvironment}{\let\def\TeXPrimitiveDef\dfn}
\AddToHook{env/def/begin}{\let\def\BeginDefEnvironment}

\theoremstyle{remark}
\newaliascnt{rem}{thm}
\newtheorem{rem}[rem]{Remark}
\aliascntresetthe{rem}

\newenvironment{lemma}{\lem}{\endlem}
\newenvironment{corollary}{\cor}{\endcor}
\newenvironment{definition}{\dfn}{\enddfn}
\newenvironment{remark}{\rem}{\endrem}

\crefname{thm}{theorem}{theorems}
\Crefname{thm}{Theorem}{Theorems}
\crefname{prop}{proposition}{propositions}
\Crefname{prop}{Proposition}{Propositions}
\crefname{lem}{lemma}{lemmas}
\Crefname{lem}{Lemma}{Lemmas}
\crefname{cor}{corollary}{corollaries}
\Crefname{cor}{Corollary}{Corollaries}
\crefname{dfn}{definition}{definitions}
\Crefname{dfn}{Definition}{Definitions}
\crefname{assumption}{assumption}{assumptions}
\Crefname{assumption}{Assumption}{Assumptions}
\crefname{rem}{remark}{remarks}
\Crefname{rem}{Remark}{Remarks}

\newcommand{\sq}[1]{\widetilde{#1}}
\newcommand{\R}{\mathbb{R}}

\newcommand{\mc}[1]{\mathcal{#1}}
\newcommand{\p}{\partial}

\newcommand{\ov}[1]{\overline{#1}}
\newcommand{\ang}[1]{\langle #1 \rangle}

\DeclareMathOperator{\sign}{sign}

\newcommand{\Ltwo}[1]{\left\lVert #1\right\rVert_{L^2}}
\newcommand{\Linf}[1]{\begingroup\delimiterfactor=750\relax
  \left\lVert #1\right\rVert_{L^\infty}\endgroup}
\newcommand{\LinfR}[1]{\begingroup\delimiterfactor=750\relax
  \left\lVert #1\right\rVert_{L^\infty(I_R)}\endgroup}
\newcommand{\Honehalf}[1]{\left\lVert #1\right\rVert_{H^{1/2}}}
\newcommand{\Hone}[1]{\left\lVert #1\right\rVert_{H^1}}
\newcommand{\Htwo}[1]{\left\lVert #1\right\rVert_{H^2}}

\newcommand{\X}[1]{\left\lVert #1\right\rVert_{X}}

\title{Self-similar solutions for capillary water waves}

\date{\today}

\author[Carballeira]{Alvaro Carballeira}
\address{Department of Mathematics, Brown University, 02912 Providence, Rhode Island, USA}
\email{alvaro\_carballeira\_mora@brown.edu}

\author[Castro-López]{Gerard Castro-López}
\address{Department of Mathematics, Courant Institute, New York University, 10012 New York, New York,  USA}
\email{gerard\_castro\_lopez@cims.nyu.edu}

\author[García-Juárez]{Eduardo García-Juárez}
\address{Departamento de Análisis Matemático, Universidad de Sevilla, 41011 Sevilla, Spain}
\email{egarcia12@us.es}

\author[Gómez-Serrano]{Javier Gómez-Serrano}
\address{Department of Mathematics, Brown University, 02912 Providence, Rhode Island, USA}
\email{javier\_gomez\_serrano@brown.edu}

\author[Pasquariello]{Marcus Pasquariello}
\address{Department of Mathematics, Brown University, 02912 Providence, Rhode Island, USA}
\email{marcus\_pasquariello@brown.edu}

\begin{document}

\begin{abstract}
  In this paper we prove the existence of self-similar solutions to the capillary water wave problem. In particular, this implies finite-time formation of singularities in the form of a corner, and retraction from a corner due to the time-reversibility of the equations. Unlike most self-similar smooth solutions, the solutions with $2/3$ time-length scaling are seen in experiments like pinch-off retraction, bubble bursting or jet formation.
  The proof involves a new argument through the regularization of the Hamiltonian and does not require the existence of an approximate profile. We believe this strategy has further implications in the construction of self-similar solutions of dispersive PDE.
\end{abstract}

\maketitle

\setcounter{tocdepth}{1}
\tableofcontents

\section{Introduction}

\subsection{The equations}

We consider the two-dimensional free boundary Euler equations with surface tension: the capillary water-waves system 
\begin{equation}\label{eq:euler}
\begin{cases}
 \displaystyle \Delta\Phi = 0 & \text{in }\Omega(t),\\
 \displaystyle \eta_t =  \Phi_Y -  \eta_X  \Phi_X &\text{on } Y=\eta(X,t),\\
 \displaystyle\Phi_t+\frac{1}{2} |\nabla\Phi|^2 = \frac{\p_X^2\eta}{(1+ (\p_X\eta)^2)^{3/2}} &\text{on }Y=\eta(X,t),
\end{cases}
\end{equation}
where the domain occupied by the fluid is given by
\begin{equation*}
  \Omega(t) = \{(X,Y) \in \mathbb{R}^{2} : Y < \eta(X,t)\}.
\end{equation*}
The fluid is irrotational with velocity  given by $\nabla \Phi$, where $\Phi$ is the potential, and we will write $\phi(x,t) = \Phi(X, \eta(X,t), t)$ for the potential trace. 
The results of this paper concern the existence of self-similar solutions to this system. The scaling of \eqref{eq:euler} leads us to seek solutions of the form
\begin{equation}\label{eq:self_similar_coord}
 \eta(X,t) = (T-t)^{2/3} f\left( \frac{X}{(T-t)^{2/3}} \right),
 \quad
 \Phi(X,Y,t)=(T-t)^{1/3} \Psi\left( \frac{X}{(T-t)^{2/3}}, \frac{Y}{(T-t)^{2/3}} \right),
\end{equation}
where $T>0$ is fixed. We refer to the functions $f$ and $\Psi$ as the self-similar profiles. The self-similar equations become
\begin{equation*}
\begin{cases}
 \displaystyle \Delta\Psi = 0 & \text{in } \Omega_f,\\
 \displaystyle -\frac{2}{3} (f-x f') = \Psi_y - f'\Psi_x & \text{on } y = f(x),\\
 \displaystyle -\frac{1}{3} \Psi + \frac{2}{3}(x,y)\cdot \nabla \Psi + \frac12|\nabla\Psi|^2 =\frac{f''}{(1+f'^2)^{3/2}} & \text{on } y = f(x),
\end{cases}
\end{equation*} 
with $(x, y) = (T-t)^{-2/3}(X,Y)$ and $\Omega_f = \left\{ (x,y) : y < f(x) \right\}$.

This system of equations can be written in terms of the boundary traces alone using the Zakharov-Craig-Sulem \cite{Zakharov:stability-periodic-waves-deep-fluid,Craig-Sulem:numerical-simulation-gravity-waves} formulation. Denote $\varphi = \Psi(x, f(x))$ and let $G(f)$ be the Dirichlet to Neumann operator of the $\Omega_f$ domain. Since $\varphi'=(\Psi_x+f'\Psi_y)|_{y=f(x)}$, we can solve for the velocity components and obtain
\begin{equation}\label{eq:velocity_components}
 \Psi_x|_{y=f(x)} = \frac{\varphi'-f'G(f)\varphi}{1+f'^2},
 \qquad
 \Psi_y|_{y=f(x)} = \frac{G(f)\varphi+f'\varphi'}{1+f'^2}.
\end{equation}
Using these expressions, we find  
  \begin{equation}\label{eq:identity_flux}
  \Big[\frac{\Psi_x^2-\Psi_y^2}{2}+f'\Psi_x\Psi_y\Big]\Big|_{y=f(x)}=\frac12(\varphi')^2-\frac{(G(f)\varphi+f'\varphi')^2}{2(1+f'^2)},
  \end{equation}
and the self-similar equations become
\begin{equation}\label{eq:selfsimilar}
\begin{cases}
  \displaystyle -\frac{2}{3} (f-x f') = G(f)\varphi,\\
  \displaystyle -\frac{1}{3} \varphi + \frac{2}{3} x \varphi' = \frac{f''}{(1+f'^2)^{3/2}} - \frac{1}{2} (\varphi')^2 + \frac{(G(f)\varphi + f'\varphi')^2}{2(1+f'^2)}.
\end{cases}
\end{equation}
These are the equations that will be employed mostly throughout the paper. A particular feature we are going to use is that they have a Hamiltonian:
\begin{equation}\label{eq:def_hamiltonian}
  \mc{H}_s[f,\varphi] := \frac{1}{2}\int_{\R}\varphi\,G(f)\varphi \, d x + \lim_{R \to \infty}\int_{-R}^{R}\bigl(\sqrt{1+f'^2}-\sqrt{1 + s^2}\bigr) \, d x,
\end{equation}
where $s$ is the slope of $f$ at infinity.

\subsection{Historical background}

\subsubsection*{The Cauchy problem} The mathematical study of the water waves equations has a long history; we refer to the book~\cite{Lannes:water-waves-book} and to the survey~\cite{Ionescu-Pusateri:survey-water-waves} for a detailed account. In the absence of surface tension, local well-posedness for small data was established by Nalimov~\cite{Nalimov:cauchy-poisson}, Yosihara~\cite{Yosihara:gravity-waves} and Craig~\cite{Craig:existence-theory-water-waves}, and for general data in Sobolev spaces by Wu~\cite{Wu:well-posedness-water-waves-2d,Wu:well-posedness-water-waves-3d}. See also~\cite{Christodoulou-Lindblad:motion-free-surface,Lindblad:well-posedness-motion,Lannes:well-posedness-water-waves,Coutand-Shkoller:well-posedness-free-surface-incompressible,Shatah-Zeng:geometry-priori-estimates,Zhang-Zhang:free-boundary-3d-euler} for extensions to more general settings. When surface tension is taken into account, as in \eqref{eq:euler}, local well-posedness was proved by  Iguchi~\cite{Iguchi:well-posedness-capillary-gravity-waves}, Ambrose~\cite{Ambrose:well-posedness-vortex-sheet}, Ambrose--Masmoudi~\cite{Ambrose-Masmoudi:zero-surface-tension-2d-waterwaves,Ambrose-Masmoudi:zero-surface-tension-3d-waterwaves}, Schweizer~\cite{Schweizer:3d-euler-free-boundary-surface-tension}, Coutand--Shkoller~\cite{Coutand-Shkoller:well-posedness-free-surface-incompressible}, Shatah--Zeng~\cite{Shatah-Zeng:geometry-priori-estimates} and, at lower regularity, by Alazard--Burq--Zuily~\cite{Alazard-Burq-Zuily:water-wave-surface-tension}.

\begingroup
\setlength{\emergencystretch}{1em}
Concerning the long time behavior, global or almost global solutions for small, smooth and localized data were constructed for gravity waves by Wu~\cite{Wu:almost-global-wellposedness-2d,Wu:global-wellposedness-3d}, Germain--Masmoudi--Shatah~\cite{Germain-Masmoudi-Shatah:global-solutions-gravity-water-waves-annals}, Ionescu--Pusateri~\cite{Ionescu-Pusateri:global-solutions-water-waves-2d}, Alazard--Delort~\cite{Alazard-Delort:global-solutions-2d-water-waves} and Hunter--Ifrim--Tataru~\cite{Hunter-Ifrim-Tataru:holomorphic-coordinates-2d-water-waves,Ifrim-Tataru:global-existence-2d-water-waves}; for capillary waves by Germain--Masmoudi--Shatah~\cite{Germain-Masmoudi-Shatah:global-existence-capillary-water-waves}, Ionescu--Pusateri~\cite{Ionescu-Pusateri:global-solutions-water-waves-2d-surface-tension} and Ifrim--Tataru~\cite{Ifrim-Tataru:global-existence-2d-water-waves-surface-tension}; and for gravity-capillary waves by Deng--Ionescu--Pausader--Pusateri~\cite{Deng-Ionescu-Pausader-Pusateri:global-solutions-gravity-capillary-water-waves-3d}. All these results require the interface to be considerably smoother than Lipschitz, which is the critical regularity for \eqref{eq:euler}. Therefore, the corner considered in this paper is out of the reach of the available well-posedness theory.
\par\endgroup

\subsubsection*{Global continuation and degree theory} Our proof builds on global continuation arguments,  which use a topological degree and a priori bounds for the set of solutions to replace the perturbative analysis of the linearized operator. 
The idea goes back to Leray--Schauder~\cite{Leray-Schauder:topologie-equations-fonctionnelles} and serves as a starting point for the  global bifurcation theorem of Rabinowitz~\cite{Rabinowitz:global-results-nonlinear-eigenvalue} (see also Dancer~\cite{Dancer:bifurcation-theory-analytic-operators} and the book of Buffoni--Toland~\cite{Buffoni-Toland:analytic-theory-global-bifurcation}). We refer to the surveys~\cite{Toland:stokes-waves,Strauss:steady-water-waves,Haziot-Hur-Strauss-Toland-Wahlen-Walsh-Wheeler:water-waves-survey} for the extensive use of these global bifurcation arguments in the water problem: Keady--Norbury~\cite{Keady-Norbury:existence-irrotational-water-waves} constructed large amplitude periodic gravity waves, then Amick--Toland~\cite{Amick-Toland:solitary-water-waves-finite-amplitude} solitary waves, while Constantin--Strauss~\cite{Constantin-Strauss:exact-steady-periodic-water-waves-vorticity} included the case with vorticity. See also more recently Wheeler~\cite{Wheeler:large-amplitude-solitary-water-waves-vorticity} and Chen--Walsh--Wheeler~\cite{Chen-Walsh-Wheeler:stratified-solitary-water-waves} for solitary waves with vorticity and with stratification, where the problem is posed in an unbounded domain, and  Constantin--Strauss--V\u{a}rv\u{a}ruc\u{a}~\cite{Constantin-Strauss-Varvaruca:global-bifurcation-gravity-water-waves} for waves with critical layers. We also mention the result by Jia--\v{S}ver\'ak~\cite{Jia-Sverak:local-estimates-weak-ns} of forward self-similar solutions to the Navier--Stokes equations, whose proof relies on a combination of a priori estimates with the Leray--Schauder degree, which will also play an important role in our proof strategy.

\subsubsection*{Singularities} Much less is known about the formation of singularities. Castro--C\'ordoba--Fefferman--Gancedo--G\'omez-Serrano~\cite{Castro-Cordoba-Fefferman-Gancedo-GomezSerrano:splash-water-waves,Castro-Cordoba-Fefferman-Gancedo-GomezSerrano:finite-time-singularities-free-boundary-euler} proved that  \emph{splash} and \emph{splat} singularities in finite time occur in the water waves system. In these singularities,  the interface self-intersects at a point or along an arc, respectively, while it remains smooth. They later showed that the same scenario persists in the presence of surface tension~\cite{Castro-Cordoba-Fefferman-Gancedo-GomezSerrano:finite-time-singularities-water-waves-surface-tension}, while Coutand--Shkoller~\cite{Coutand-Shkoller:finite-time-splash} extended these results to three dimensions and to rotational flows. 
In all these scenarios the singularity is due to the collapse of the geometry, and both the interface and the velocity remain smooth up to the singular time. 
Moreover, these singularities cannot occur in the two-phase setting, as shown by Fefferman--Ionescu--Lie~\cite{Fefferman-Ionescu-Lie:absence-splash-singularities} and Coutand--Shkoller~\cite{Coutand-Shkoller:no-splash-vortex-sheet} (although stationary splash singularities were constructed by C\'ordoba--Enciso--Grubic~\cite{Cordoba-Enciso-Grubic21}).
The singularity constructed in this paper is of a different nature: the interface does not self-intersect, and it loses regularity because its curvature blows up at a point.

There are also results on the growth of solutions that remain regular. For the free boundary Euler equations with surface tension and nonzero vorticity, in a periodic domain of finite depth, Hu--Luo--Yao~\cite{Hu-Luo-Yao:small-scale-creation-free-boundary-euler-surface-tension} constructed initial data with a flat interface and arbitrarily small velocity for which the gradient of the vorticity grows at least double exponentially in time for as long as the solution exists. Berti--Grande--Maspero--Staffilani~\cite{Berti-Grande-Maspero-Staffilani:rogue-waves-large-deviations} characterized the tail probability of rogue waves for small Gaussian random initial data for gravity water waves, with dispersive focusing explaining their formation.

\subsubsection*{Interfaces with corners} A classical example of a  water wave with a corner is the Stokes wave of greatest height, whose existence was proved by Amick--Fraenkel--Toland~\cite{Amick-Fraenkel-Toland:stokes-conjecture-wave-extreme} and Plotnikov~\cite{Plotnikov:stokes-conjecture}. 
It is a a traveling wave for the gravity water waves system with a corner of angle $2\pi/3$ at the crest. 
For the Cauchy problem without surface tension, Kinsey--Wu~\cite{Kinsey-Wu:apriori-estimates-water-waves-crests} and Wu~\cite{Wu:well-posedness-water-waves-corners} developed a local well-posedness theory in weighted spaces that allows for interfaces with angled crests, and Agrawal~\cite{Agrawal:rigidity-corners-water-waves} proved that these singularities are rigid: an initial interface with angled crest remains so, and the angle does not change in time. C\'ordoba--Enciso--Grubic~\cite{Cordoba-Enciso-Grubic:local-wellposedness-water-waves-cusps-corners} constructed solutions with cusps and with corners whose angle varies in time, and proved in~\cite{Cordoba-Enciso-Grubic26} that angled-crested water waves can become singular in finite time. In the absence of both gravity and surface tension, Liu--Pego~\cite{Liu-Pego:local-singularities-ideal-potential-flows,Liu-Pego:search-local-singularities-ideal-potential-flows} gave numerical evidence that corners may form in a self-similar way. In the previous results~\cite{Kinsey-Wu:apriori-estimates-water-waves-crests,Wu:well-posedness-water-waves-corners,Cordoba-Enciso-Grubic:local-wellposedness-water-waves-cusps-corners}, the initial corner  persists for as long as the solution exists. The situation with surface tension is expected to be very different. Agrawal~\cite{Agrawal:angled-crested-surface-tension-wellposedness,Agrawal:angled-crested-surface-tension-zero-limit} extended the weighted energy of Kinsey--Wu to positive surface tension, and showed that the resulting energy does not allow for angled crests, although it does allow for interfaces with large curvature, which converge to angled crested solutions in the zero surface tension limit. To the best of our knowledge, there are no previous results on the evolution of an interface with a corner under surface tension, or on the formation of a corner from a smooth capillary water wave.

\subsubsection*{Numerical studies of singularities} There is an extensive numerical literature on the loss of regularity of fluid interfaces. For vortex sheets without surface tension, Moore~\cite{Moore:singularity-vortex-sheet} predicted by formal asymptotics that the curvature of the sheet becomes singular in finite time, and this was confirmed and refined in~\cite{Meiron-Baker-Orszag:vortex-sheet,Krasny:singularity-vortex-sheet-point-vortex,Shelley:singularity-vortex-sheet-spectrally-accurate,Cowley-Baker-Tanveer:moore-curvature-singularities}. 
Then Baker--Caflisch--Siegel~\cite{Baker-Caflisch-Siegel:singularity-rayleigh-taylor} studied the analogous phenomenon for the two-fluid Rayleigh--Taylor instability. For water waves, Tanveer~\cite{Tanveer:singularities-water-waves-rayleigh-taylor} 
analyzed the singularities of the analytic continuation of the conformal map for steady gravity waves
and Baker--Xie~\cite{Baker-Xie:singularities-complex-plane-deep-water-waves} extended the technique to the time-dependent case.  Fontelos--de la Hoz~\cite{Fontelos-delaHoz:singularities-water-waves} gave further evidence that singularities of Moore type develop in the two-phase setting, and also proposed a singularity in the form of a logarithmic spiral for the one-fluid problem. On the negative side, Wilkening~\cite{Wilkening:breakdown-self-similarity-standing-waves} computed standing water waves of large amplitude and found that their crests do not sharpen into a corner in a self-similar way. 
The effect of surface tension on these scenarios was investigated by Siegel~\cite{Siegel:singularity-kelvin-helmholtz-surface-tension}, de la Hoz--Fontelos--Vega~\cite{delaHoz-Fontelos-Vega:surface-tension-moore-singularity}, Fontelos--de la Hoz~\cite{Fontelos-delaHoz:singularities-water-waves}, and  Ceniceros--Hou~\cite{Ceniceros-Hou:dynamic-generation-capillary-waves}. In~\cite{delaHoz-Fontelos-Vega:surface-tension-moore-singularity,Fontelos-delaHoz:singularities-water-waves}, surface tension regularizes the curvature singularity.

\subsubsection*{Self-similar flows driven by surface tension} The self-similar ansatz \eqref{eq:self_similar_coord} goes back to Keller--Miksis~\cite{Keller-Miksis:surface-tension-driven-flows}, who computed the corresponding profiles numerically.
The study of changes in topology of an inviscid fluid now commonly uses this self-similar ansatz.
Lawrie~\cite{Lawrie:surface-tension-driven-flow-wedge} solved the linearized problem for a corner in contact with a wall, and Billingham~\cite{Billingham:surface-tension-fat-wedges-cones} studied interfaces with a corner which are close to a half-space, which is the regime considered in this paper. The merging and the retraction of these corners were computed by Keller--Milewski--Vanden-Broeck~\cite{Keller-Milewski-Vandenbroeck:merging-wetting-surface-tension}, Billingham--King~\cite{Billingham-King:surface-tension-slender-wedge-coalescence} and Sierou--Lister~\cite{Sierou-Lister:self-similar-recoil-inviscid-drops}. The same scaling governs the pinch-off of an inviscid drop~\cite{Day-Hinch-Lister:self-similar-capillary-pinchoff,Leppinen-Lister:capillary-pinchoff-inviscid-fluids}, the coalescence of drops~\cite{Duchemin-Eggers-Josserand:inviscid-coalescence-drops} and the pinching of vortex sheets with surface tension~\cite{Hou-Lowengrub-Shelley:long-time-motion-vortex-sheets-surface-tension}. 
It also describes the collapse of a film bridge~\cite{Chen-Steen:inviscid-capillary-breakup-film-bridge} and the collapse of the cavity left by a bursting bubble, which precedes the eruption of a jet~\cite{Duchemin-Popinet-Josserand-Zaleski:jet-formation-bubbles-bursting,Lai-Eggers-Deike:bubble-bursting-universal-profiles}. See also Burton--Taborek~\cite{Burton-Taborek:dimensionality-pinch-off-coalescence} for pinch-off and coalescence in two dimensions.

The work which is closest to the singularity constructed in this paper is the one of Zeff--Kleber--Fineberg--Lathrop~\cite{Zeff-Kleber-Fineberg-Lathrop:singularity-curvature-collapse-jet}, who studied, experimentally and numerically, the collapse of the depression of a standing wave in a regime in which the interface does not break. They considered the axisymmetric version of \eqref{eq:euler}, neglecting gravity near the singularity, and looked for solutions of the form \eqref{eq:self_similar_coord}, with the same exponents. They observed that the profile has to be asymptotically linear, so that the interface becomes a cone at the singular time, computed the profile numerically, and found good agreement with the experiments, in which the collapse is followed by the eruption of a thin jet. In their solutions the velocity potential grows like the square root of the distance to the origin, and the velocity of the fluid diverges at the tip of the cone at the singular time.  Remark \ref{rmk:nonzero_velocity} shows how our result covers these phenomena.

We refer to~\cite{Eggers:nonlinear-dynamics-breakup-free-surface-flows,Eggers-Fontelos:singularities-book} for a comprehensive overview. All these works are numerical or based on formal asymptotics, and the existence of the self-similar profiles has not been proved in any of these settings.

\subsection{Self-similar solutions and main results}

The main result of this paper is  the following theorem:

\begin{thm}\label{thm:main_theorem}
There exists $s_0>0$ such that for every $|s|\leq s_0$, there exists a locally smooth, even $(f,\varphi) \in C^\infty_{loc} \times C^\infty_{loc}$ solution to the self-similar equations \eqref{eq:selfsimilar} satisfying
\begin{equation*}
\begin{aligned}
  \Hone{f - s |x|} + \|\varphi\|_{H^{1/2}} + \Linf{f'}  &\leq C|s|,
\end{aligned}
\end{equation*}
and with  Hamiltonian $\mc{H}_s[f,\varphi] = 0$. The full domain potential $\Psi$ (harmonic extension of $\varphi$) satisfies
\begin{equation*}
  \left \|\frac{\nabla \Psi}{\langle x\rangle } \right \|_{L^\infty(\Omega_f)}^2 + \int_{\Omega_f} \left |\nabla \Psi \right |^2 \leq C s^2.
\end{equation*}
\end{thm}

The following corollaries refer to the physical space-time interpretation of these self-similar profiles. Without loss of generality, we set $T=0$.
At $t=0$ consider the  interface with a corner $Y<s|X|$, and denote by $\theta$ the angle of the corner measured from the fluid side. We say an angle $\theta$ is admissible if
\begin{equation}\label{eq:theta_admissible_range}
  \theta = \pi + 2 \arctan(s), \quad \quad 0 < |\theta- \pi|\leq 2\arctan(s_0).
\end{equation}
The forward evolution is given by the self-similar rescaling \eqref{eq:self_similar_coord} of the profiles $(f, \varphi)$ from Theorem \ref{thm:main_theorem} for $t<0$, while for $t>0$ it suffices to consider its time reversed profiles.
For all times except at $t=0$ the surface and the potential are locally smooth, their relative Hamiltonian \eqref{eq:def_hamiltonian} vanishes, and we have $(\eta-s|X|, \varphi) \in H^1\times H^{1/2}$.
The velocity $u= \nabla \Phi$ converges  to zero in $L^2$ in both one-sided traces at $t=0$.
From now on we will use $\theta$ and  $s$ and assume throughout that they are related by $s = \tan ((\theta-\pi)/2)$.

\begin{figure}[htbp]
  \centering
  \includegraphics{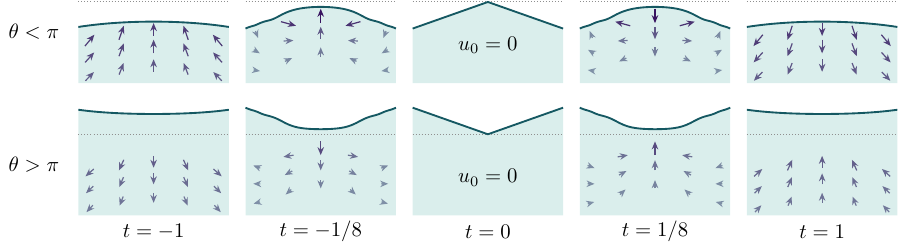}
  \caption{Finite time corner formation and immediate retraction}
  \label{fig:corner-evolution}
\end{figure}

\begin{cor}[Finite time corner formation]\label{cor:corner_formation}
  For every admissible  angle $\theta$ \eqref{eq:theta_admissible_range}, there exists a self-similar evolution with locally smooth initial data $(\eta, \phi)|_{t=-1} = (f, \varphi)$, such that as $t\uparrow 0$ the surface converges uniformly  $\eta \to s|X|$ and the velocity $u \to 0$ in $L^2$.  Moreover, at the central point of the surface $P_t = (0, \eta(0,t))$, the signed curvature and the velocity satisfy
  \begin{equation}\label{eq:kappa_u_blowup_rate}
    \kappa(0,t) = f''(0)|t|^{-2/3}, \quad u(P_t, t) = \left ( 0, -\frac{2}{3}f(0)|t|^{-1/3}\right ),
  \end{equation}
  Since $sf(0), sf''(0) >0$, the maximum of $|\kappa|$ and $|u|$ in any fix neighborhood of the origin grow at least at these respective rates. 
  See $t\leq 0$ of Figure \ref{fig:corner-evolution}. 
\end{cor}

\begin{remark}
  The formal linear far field asymptotics \eqref{eq:linear_oscillations} suggest that the singularity formation mechanism is dispersive focusing. 
  At distance $R$ and time $t$ it predicts a wave packet of leading wave number $(2R/(3|t|))^2$ and inward group velocity $R/|t|$, which suggests all the wave packets will arrive at the origin at the same time $t=0$.
\end{remark}

\begin{cor}[Immediate retraction from resting corner]\label{cor:corner_retraction}
  For every admissible $\theta$, the initial corner at rest $(\eta, \phi)|_{t=0} := (s|X|, 0)$ has a self-similar evolution which is smooth for every $t>0$. It is the time-reversal of the finite time corner formation solution, so the identities \eqref{eq:kappa_u_blowup_rate} are the same changing the sign of velocity. See $t\geq 0$ of Figure \ref{fig:corner-evolution}.
\end{cor}
\begin{remark}
  The fluid starts at rest $u(\cdot, t) \downarrow 0$ in $L^2$, even though the velocity at the tip goes to infinity pointwise $|u(P_t, t)| \to \infty$.
\end{remark}

\begin{figure}[htbp]
  \centering
  \includegraphics{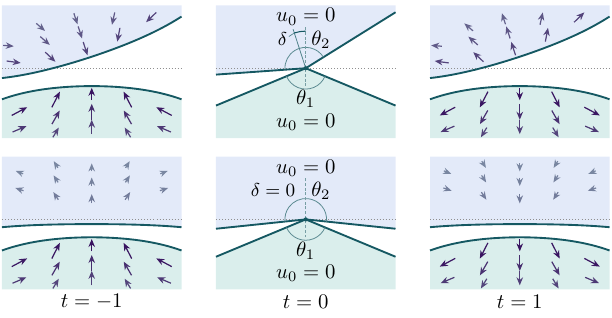}
  \caption{First contact between two fluid interfaces at a corner.}
  \label{fig:corner-contact}
\end{figure}

\begin{corollary}[Separation and first contact of two fluid corners]\label{cor:separation_two_fluids}
  Fix an admissible $\theta_1< \pi$. There exists $\Delta(\theta_1)>0$, such that for every admissible $\theta_2$ and tilt angle $\delta$ satisfying 
  \begin{equation}
    (\theta_2 -\pi)_{+} + |\delta| < \Delta(\theta_1), \quad  a_+ = \max\{0, a\},
  \end{equation}
  there are two fluid components that evolve self-similarly and whose boundaries meet for the first time at $t=0$ at the origin. At the touching point they form corners of angles $\theta_1, \theta_2$ and they have a $\delta$ tilt between their opposite bisector lines. Time reversing, we can start with any such angles and tilt, and obtain immediate separation and regularization. See Figure \ref{fig:corner-contact}.
\end{corollary}

\begin{corollary}[Corner retraction and formation in contact with a wall]\label{cor:fluid_wall}
  The Corollaries \ref{cor:corner_formation}, \ref{cor:corner_retraction} and \ref{cor:separation_two_fluids} also hold in a half plane with a solid boundary along their axis of symmetry. For Corollary \ref{cor:separation_two_fluids} take $\delta=0$.
  Each initial wedge with contact angle \(\theta/2\), where \(\theta\) is admissible, admits a self-similar evolution with $90^\circ$ fluid-wall contact angle for every $t>0$. Time reversing we obtain $\theta/2$ wedge formation from $90^\circ$ contact angle.
\end{corollary}

\begin{rem}[Gravity-capillary water waves]
  For the infinite-depth, gravity-capillary water waves problem, a change of variables to a falling moving frame removes the gravity term. For the solution $(\eta,\Phi)$ in \eqref{eq:self_similar_coord}, set
  \begin{align*}
    \eta_g(X,t) &= \eta(X,t) - \frac{1}{2}g (T-t)^2, \\
    \Phi_g(X,Y,t) &= \Phi\left(X,Y+\frac{1}{2}g (T-t)^2,t\right) + g (T-t) Y + \frac{1}{6} g^2 (T-t)^3.
  \end{align*}
  This gives a gravity-capillary solution that is self-similar in the new coordinates and such that $\eta_g(\cdot,t)$ converges to $s|X|$ as $t \to T$. With this self-similar scaling the gravity is lower order close to the singularity. In this case it only has finite energy with respect to the falling frame.
\end{rem}

The result and its corollaries presented so far have natural and interesting extensions that follow  using the same proof strategy. In the following remarks, we outline the additional treatment that they require. We plan to include some of these extensions in a future version of this manuscript, and to formalize all of them along with the main Theorem \ref{thm:main_theorem}.

\begin{rem}[Two-phase]
The same result holds in the two-phase setting, that is, for two irrotational, incompressible fluids of constant densities $\rho_\pm>0$. In the two-phase Zakharov-Craig-Sulem formulation, the self-similar interface variables are $f$ and
\begin{equation*}
    \varphi(x)=\frac{\rho_{-}\Psi_{-}(x,f(x))-\rho_{+}\Psi_{+}(x,f(x))}{\rho_{-}+\rho_{+}},
\end{equation*}
where $\Psi_{-}$ and $\Psi_{+}$ are the velocity potentials in the lower and upper fluids, respectively. The proof follows the same scheme and methods, with analogous energy and stress balances  accounting for both fluids.
\end{rem}

\begin{rem}[Non-zero terminal velocity: Jet formation] \label{rmk:nonzero_velocity} 
  In \Cref{cor:corner_formation}, we have a solution with finite kinetic energy whose velocity $u = \nabla \Phi $ converges to zero in $L^{2}$ as $t \to T$. Indeed, since $u$ has finite energy, in the self-similar coordinates one finds
  \begin{equation*}
    \int_{\Omega(t)} |u|^2 \,dX \,dY = |T-t|^{2/3} \int_{\Omega_{f}} |\nabla\Psi(x,y)|^2 \,dx \,dy \, \to \, 0 \text{ as } t\uparrow T.
  \end{equation*}
  Our result can be extended to nonzero terminal velocity $u(\cdot, T)$. Since we cannot directly work in a finite energy space for $\Psi$, we are going to work with 
  \begin{equation}
      \Psi_1(x,y) = \Psi(x,y) - \alpha Q(x,y) \quad r(x) = f(x) - s|x| - g_{s, \alpha}\langle x \rangle^{-1/2},
  \end{equation}
  where $Q(x,y) \sim \text{Re}([1+ i(x+iy)]^{1/2})$ is harmonic and smooth in the domain, and $g_{s,\alpha}$ is an explicit constant. Working with $\Psi_1$, and its trace $\varphi_1$ and the new $r$, we can obtain analogously the energy and a priori bootstrap estimate as in Sections \ref{sec:energy} and \ref{sec:uniform_bounds} when $|s|+ |\alpha|$ is small enough.

  In this case, at singular time the vertical velocity $u_2$ and the corner $f$ become 
  \begin{equation}
      u_2(X, T) \sim c_s \alpha |X|^{-1/2}, \quad  \eta(X, T) = s|X|,
  \end{equation}
  so the velocity at the tip of the corner is infinite, and the corner profile $\eta$ does not depend on $\alpha$. Thus, for the same wedge we have two self-similar solutions at the same $|\alpha|$, the negative and positive sign. This changes the corner formation and continuation picture substantially: instead of having the same dynamics before and after the singularity, this gives two solutions corresponding to opposite signs of $\alpha$.  

  Moreover with the forward sign convention, for $\alpha \gtrsim |s|$, the self-similar profiles have a higher origin: for $s<0$ we have $f_{s,\alpha}(0)>0$ and for $s>0$ we have that $f_{s,\alpha}(0) > 0$ is a local maximum instead of a minimum. In these cases, in physical variables this solution produces a jet, since the fluid velocity at the tip for $t>T$ is $~(t - T)^{-1/3}f_{s, \alpha}(0)>0$. This is the phenomenon observed in \cite{Zeff-Kleber-Fineberg-Lathrop:singularity-curvature-collapse-jet} explained in the introduction (see the next Remark \ref{rmk:axisymmetric}, since the experiment is actually observed in the three-dimensional axisymmetric setting).   
\end{rem}

\begin{rem}[Axisymmetric without swirl]\label{rmk:axisymmetric}
    All the previously mentioned results also hold for irrotational, axisymmetric flows without swirl. The proof follows the same strategy with two main differences. First, the cone does not have zero curvature in the azimuthal direction
    \begin{equation}
      z = sr, \quad \kappa = \frac{s}{r \sqrt{1+s^2}},
    \end{equation}
    away from the vertex. This makes both the curvature and, for the non-zero velocity, the quadratic Bernoulli terms not be in $L^2(r dr)$ ($1/r$ decay), so it needs a relative energy estimate.

    The other difference is the axisymmetric Laplacian: $\Delta_{r,z} = \partial_r^2 + \frac{1}{r}\partial_r + \partial_z^2$, but importantly at the small scales far from the axis it converges to the $2$d one, which is what is most essentially used in this proof .
\end{rem}

\subsection{Discussion of the proof}\label{subsec:strategy}

In other contexts \cite{GarciaJuarez-GomezSerrano-Nguyen-Pausader:self-similar-muskat,Na:self-similar-muskat-3d,Correia-Cote-Vega:asymptotics-fourier-self-similar-mkdv,Correia-Cote-Vega:self-similar-dynamics-mkdv,GomezSerrano-Pasquariello:self-similar-gsqg}, the natural strategy to construct self-similar solutions is to construct an approximate solution to the self-similar equation, and to then perform a fixed point argument in the correction to this approximate solution. In the setting of capillary water waves, where the linear dynamics are dispersive, the linearized operator takes the form
\begin{equation*}
  \mathcal L_{0,0}(f,\varphi)
  =\begin{pmatrix}
    \displaystyle -\frac23(f-xf')-|D|\varphi \\
    \displaystyle -\frac13\varphi+\frac23 x\varphi'-f''
  \end{pmatrix}.
\end{equation*}
Self-similar solutions for dispersive equations have been constructed before using the same strategy described above. See, for instance,~\cite{Correia-Cote-Vega:asymptotics-fourier-self-similar-mkdv,GomezSerrano-Pasquariello:self-similar-gsqg}. In these works, the strategy implemented consisted of building the approximate solution by taking the solution to the linearized equation, and adding sufficiently many corrections to it, so that the approximate solution would be good enough to close the fixed point argument. In the case of capillary water waves, we do expect this to be possible, but it would be a real challenge to implement. For starters, the even solution with slope $s$ to the linearized equation has asymptotics
\begin{equation}\label{eq:linear_oscillations}
\begin{aligned}
  f_{\mathrm{lin}}(x)
  &= s|x|-\frac{s}{\pi |x|^2}
     -\frac{27s}{4\sqrt{\pi}\,|x|^{7/2}}
      \cos\!\left(\frac{4|x|^3}{27}-\frac{\pi}{4}\right)
     +o\!\left(s|x|^{-7/2}\right),\\
  \varphi_{\mathrm{lin}}(x)
  &= \frac{9s}{2\sqrt{\pi}\,|x|^{5/2}}
      \sin\!\left(\frac{4|x|^3}{27}-\frac{\pi}{4}\right)
     +O\!\left(|s|\,|x|^{-4}\right)
\end{aligned}
\end{equation}
as $|x| \to \infty$. This highly-oscillatory behavior is difficult to manage in a fixed point argument, where the function spaces would need to be adapted to the oscillations, for instance in weighted $L^\infty$ spaces as in \cite{Correia-Cote-Vega:asymptotics-fourier-self-similar-mkdv,GomezSerrano-Pasquariello:self-similar-gsqg}. While it is in principle possible to implement a fixed point argument here, many technical complications such as the following arise:
\begin{enumerate}
    
    \item One needs to construct an approximate solution with sufficiently many terms as well as a correction to the phase. Such a construction is difficult to implement for several reasons, but one being the presence of quasilinear  quadratic nonlinearities (which include non-explicit terms in the case of the Dirichlet to Neumann operator.) One could potentially overcome some of this difficulty with normal forms, but we expect it to be technically challenging.
    \item One needs to bound the Dirichlet to Neumann operator in weighted $L^\infty$ spaces adapted to the oscillations.
    \item The inverse of the linearized operator $\mc{L}_{0,0}^{-1}$ maps generic Schwartz functions into $H^{\beta}\times H^{\gamma}$ only for $\beta<3/2, \gamma<1$, due to the oscillations of the homogeneous solutions, so possibly a refined argument like the ones mentioned above would have to be implemented.
\end{enumerate}
In this paper, we take a different approach. We regularize the equations by means of their Hamiltonian and prove estimates uniform in the regularization. The regularization allows us to work in a function space (see \Cref{def:X_Y_spaces}) where the local continuation from the $(0,0)$ solution is much simpler than the setting described above. 
We then use a priori estimates to obtain that the regularized solutions don't escape certain set, and hence, using a topological argument we get existence up to certain $s_0$ independent of the regularization. We fully avoid the need of an approximate profile.
Specifically, we consider the regularized equation replacing $\mathbb{J}$ by $\mathbb{J} + \epsilon\mathbb{I}$, where $\mathbb{J}$ is the antisymmetric unit matrix in \eqref{def:F_epsilon}. This makes the system elliptic and damps the oscillations at infinity, which allows us to get improved energy estimates and to get existence of solutions for small $s$ by a continuation argument. 

\textbf{\emph{Energy and stress balance.}} Since the Hamiltonian of a self-similar solution scales like $(T-t)^{2/3}$ and is conserved, it must vanish. For the regularized system, we prove in \Cref{prop:energy_balance} that
\begin{equation*}
  \mc{H}_s[f,\varphi] + \frac{3\epsilon}{2}\left( \left\|\frac{\delta\mc{H}_s}{\delta f}\right\|_{L^2}^2 + \left\|\frac{\delta\mc{H}_s}{\delta\varphi}\right\|_{L^2}^2\right) = 0.
\end{equation*}
The structure of the identity above is the main reason for our choice of the regularization: it provides a  positive quantity in the energy balance while respecting the structure at $\epsilon=0$.
This identity bounds the high frequencies of the solution by $s f(0)$, and the low frequencies are controlled by the scaling terms of the equations. Together, they yield $\|f - s|x|\|_{H^1} + \|\varphi\|_{H^{1/2}} \leq C|s|$, uniformly in $\epsilon$ (\Cref{prop:energy_bounds}). 
For a dispersive PDE, the energy estimate may not be strong enough to control the fine scales produced by the oscillations. We introduce a second quantity motivated by a balance of momentum across a vertical cut in the fluid.
This second balance is local: the  stress $\mc{R}_{s,\epsilon}$ defined in \Cref{def:relative_stress}, which is quadratic in the slope and in the velocity, satisfies $|\mc{R}_{s,\epsilon}(x)| \lesssim s^2$ pointwise away from the origin (\Cref{prop:stress_bounds}).

\textbf{\emph{Uniform slope and velocity bounds.}} 
The previous energy and stress estimates can be combined with the equations to show that stronger regularity, uniform in the regularization parameter, is controlled by the far-field slope $|s|$. 
We introduce the quantity $M(f,\varphi)$ \eqref{def:Mdef}, which  measures the slope of the interface and the weighted velocity, which informally speaking can be thought of as the slope of a ripple traveling along the interface with the speed of a point with fixed self-similar coordinate (see Remark \ref{rem:ripple}).
Proposition \ref{prop:uniform_slope_bounds}  then shows that $M(f,\varphi)$ is actually controlled by $|s|$ if $(f,\varphi)\in C^2_{\mathrm{loc}}\times C^1_{\mathrm{loc}}$ is a solution whose slope is below a fixed threshold and behaves like the wedge at infinity in $L^2$ (specifically, if it is an admissible solution \eqref{def:admissible}). Roughly speaking, the result is clear on any fixed bounded interval, thanks to the energy estimate and the assumed regularity on $(f,\varphi)$ (see Case 1 in the proof of Proposition \ref{prop:uniform_slope_bounds}). The challenge is to exclude the possibility of very large slopes appearing on very small scales as $x\to\infty$, which do not affect the energy estimate. The use of the stress \eqref{def:relative_stress} becomes fundamental at this point. First, let us explain how the dispersion relation provides the natural scales to look at.  The system \eqref{eq:euler} linearized at the flat solution gives
\begin{equation*}
    \eta_{tt}=-|\nabla|^3\eta,
\end{equation*}
so the linear evolution with the corner and zero velocity at $t=T$ gives, on the Fourier side, 
\begin{equation*}
    \widehat{\eta_{XX}}(\xi,t)=\frac{2s}{\sqrt{2\pi}}\cos\big((T-t)|\xi|^{3/2}\big),
\end{equation*}
and, taking $T-t=1$, we get the self-similar solution,
\begin{equation*}
    f_{\mathrm{lin}}''(x)=\frac{2s}{\pi}\int_0^\infty\cos(\xi x)\cos(\xi^{3/2})d\xi.
\end{equation*}
For $x\to\infty$, we obtain the leading contribution by stationary phase. The phase $\beta(\xi)=x\xi-\xi^{3/2}$ is stationary at $\xi_0=\frac49x^2$, and we have $\beta(\xi_0)=\frac{4}{27}x^3$ and $\partial_x\beta(\xi_0)=\xi_0$. 
 Therefore, we have found that for large $x$, horizontal displacements of order $1/\xi_0\approx 1/x^2$ modify the phase at order one.
This scaling motivates the inclusion of the weighted velocity in $M(f,\varphi)$:
its smallness allows us to neglect quadratic terms after rescaling.  Indeed, using the identity \eqref{eq:identity_flux}, 
 the dynamic equation in \eqref{eq:selfsimilar} can be written as
\begin{equation*}
    \frac{f''}{(1+f'^2)^{3/2}}=-\frac13\varphi+\frac23x(\Psi_x+f'\Psi_y)+\frac12\big(\Psi_x^2-\Psi_y^2\big)+f'\Psi_x\Psi_y.
\end{equation*}
For large $x$,  dividing by $x^2$, we find
\begin{equation}\label{eq:dyn_rescaled:aux}
\underbrace{\frac{f''/x^2}{(1+f'^2)^{3/2}}}_{\sim \ell'}=-\frac13\frac{\varphi}{x^2}+\frac23\frac{\Psi_x}{x}+\underbrace{\frac23f'\frac{\Psi_y}{x}+\frac12\frac{\Psi_x^2-\Psi_y^2}{x^2}+f'\frac{\Psi_x\Psi_y}{x^2}}_{\lesssim M^2},
\end{equation}
where $\ell'$ is rescaled curvature (see \eqref{eq:rescaling}).
Informally speaking, the proof of Proposition \ref{prop:uniform_slope_bounds} zooms in at a point $x_0\to\infty$ on scales $1/x_0^2$. On those scales, the stress balance relates the inclination of the interface (via the surface tension energy it causes, proportional to the square of the slope) to the difference between horizontal and vertical motion, up to $O(s^2)$. This second quantity is shown to decay with depth in the technical Lemma \ref{lem:velocity_below}. In turn, the (rescaled) dynamical equation connects variations of the slope (via the curvature) with the horizontal velocity (as can be already seen in \eqref{eq:dyn_rescaled:aux}), hence the previous decay is bootstrapped, keeping both quantities controlled by $|s|$. 
The actual proof carries this argument via a contradiction argument and exploiting the harmonicity provided by the water waves equations.

\textbf{\emph{Existence for $\epsilon > 0$.}} Using the implicit function theorem would only give solutions for $|s| \leq s_\epsilon$, since we do not know whether the linearized operator remains invertible. To pass to the limit, it is necessary to have existence of solutions with $s$ uniform in $\epsilon$. To obtain this we apply a continuation argument using the degree modulo $2$ for Fredholm maps of index zero~\cite{PejsachowiczRabier1998}. We prove that the linearization of the regularized system is Fredholm of index zero (\Cref{prop:fredholm_index_zero}) and that its set of solutions with $M \leq \delta$ is compact (\Cref{prop:compactness_solution_set}). Then, we check that the degree of the map in the set $\left\{ M \leq \delta \right\}$ is locally constant in $s$, because our a priori estimates prevent solutions from reaching the boundary of the set as $s$ varies. The trivial solution gives that the degree is one at $s = 0$, and thus the continuation argument yields a solution for all $|s| \leq s_0$ and $0 < \epsilon \leq 1$ (\Cref{prop:existence_small_unif}).

\textbf{\emph{The limit $\epsilon \to 0$.}} Finally, in \Cref{sec:main_theorem_proof} we use the uniform bounds for the energy and for $M$ to pass to the limit, locally in $C^1 \times C^0$ and weakly in the energy space. The limit is a weak solution of \eqref{eq:selfsimilar}, which is smooth on compact sets by elliptic regularity bootstrap, and it inherits the bounds of \Cref{thm:main_theorem}.

Our argument differs from the previous ones in several respects. First, there is no bifurcation. We continue the trivial solution at $s = 0$ with respect to the slope $s$. Moreover, the regularized equations are quasilinear and they are not a compact perturbation of the identity, so that we use the degree for $C^1$ Fredholm maps of index zero of Pejsachowicz--Rabier~\cite{PejsachowiczRabier1998}. Finally, the equations for traveling waves are elliptic, whereas \eqref{eq:selfsimilar} is not. For this reason the continuation is carried out for a regularized problem, using a priori bounds which are uniform in the regularization. 

\begin{remark}
  Both the contradiction argument (in Proposition \ref{prop:uniform_slope_bounds}) and the complex notation used in some parts of  the proof are just used for convenience. A constructive argument that provides an explicit $s_0$ could be done, but it would be more technical and less clear. The complex notation is used in a non-essential way to simplify expressions related to the Dirichlet-to-Neumann operator.
\end{remark}

\subsection{Outline}

In \Cref{sec:preliminaries}, we state several definitions that will be used throughout the paper, including the regularized self-similar system, the function spaces, and the modified Hamiltonian. In \Cref{sec:energy}, we  rigorously derive an energy balance (\Cref{prop:energy_balance}) and a momentum balance (\Cref{lem:stress_identity}), which in turn, provides a priori estimates for the regularized self-similar system (Propositions \ref{prop:energy_bounds} and \ref{prop:stress_bounds}). We exploit these energy and stress estimates in \Cref{sec:uniform_bounds} to obtain stronger  bounds for the regularized self-similar system, uniform in $\epsilon$ (\Cref{prop:uniform_slope_bounds}). In \Cref{sec:compactness}, we  prove that the set of solutions to the regularized self-similar system is compact in $[-s_0, s_0] \times X$ (\Cref{prop:compactness_solution_set}).  \Cref{sec:existence_regularized} proves the existence of solutions to the regularized self-similar system (\Cref{prop:existence_small_unif}), using the uniform bounds and compactness of Sections \ref{sec:uniform_bounds}-\ref{sec:compactness} to close a degree argument. Finally, in \Cref{sec:main_theorem_proof}, we  take the limit as $\epsilon \to 0$ and prove the existence of solutions to the self-similar system (\Cref{thm:main_theorem}).

\subsection{Notation}
We use $X,Y$ for the original spatial variables and $t$ for time, and $x$ and $y$ for the self-similar spatial variables.

The Japanese bracket is $\ang{y}=(1+y^2)^{1/2}$. For nonnegative quantities $a,b$, we write
$a\lesssim b$ if $a\le Cb$ for a constant $C>0$, and $a\approx b$ if both $a\lesssim b$ and
$b\lesssim a$. Additional dependencies are indicated by subscripts, as in $\lesssim_{M,N}$ or $\lesssim_K$, or
specified locally. Similarly for $\gtrsim$.

We use the Fourier transform convention
\[\widehat{h}(\xi) = (2\pi)^{-1/2}\int e^{-ix\xi}h(x) \, dx.\]

In most of the paper we ommit the $s$ dependence of the solutions $(f, \varphi)$ to make notation lighter.

\subsection{Acknowledgements}
 AC, GCL, JGS and MP were partially supported by the National Science Foundation through
Grants DMS-2247537, DMS-2434314, and DMS-2554957. 
EGJ was partially supported by the AEI through grants PID2022-140494NA-I00 and PID2025-172153NB-I00, by the RYC2021-032877 research grant, and IMUS–María de Maeztu CEX2024-001517-M, funded by MICIU/AEI/10.13039/501100011033. GCL was partially supported by a Simons Investigator award (V Vicol), and received the support of a fellowship from “la Caixa” Foundation Fellowships Area (ID 100010434), the fellowship code is LCF/BQ/EU24/12060077.

We are very grateful to Vlad Vicol for useful conversations.

\subsection{AI Statement}
The third and fourth authors began working on this problem in 2013 and 2021, respectively. Most of the effort was dedicated to closing a fixed-point construction (as in the work \cite{GomezSerrano-Pasquariello:self-similar-gsqg}). These attempts did not succeed, including with the use of GPT-5.5 in June 2026 and GPT-6.0 Astra in late September.

On September 29, trying with a toy model of the gSQG front equation, GPT-6.0 Astra came out with the new strategy: use the Hamiltonian structure to choose a viscous regularization and use degree theory to bypass the construction of the nonlinear phase. Once the existence of self-similar solutions for the toy model was obtained, it extended quickly to gSQG, first for $3/2<\alpha<2$, and after some iterations to smaller values of $\alpha$. A key step in the argument was an a priori estimate for the slope obtained by a contradiction argument done via rescaling, with the scale dictated by the dispersive relation.

The first attempts to adapt the method for the capillary water waves system failed, but emphasizing to GPT-6.0 Astra to employ these techniques and take advantage of the shared properties (Hamiltonian structure, same dispersive relation as the $\alpha=3/2$ case) a first proof closed on October 4. Subsequent work has provided further information on the solutions and extensions to related settings. 

The authors have reorganized, substantially simplified, and modified the original proof and have prepared the manuscript. 
All mathematical statements and proofs have been independently and rigorously verified by the authors, who take full responsibility for the contents of the paper.

\section{Formulation and preliminaries}\label{sec:preliminaries}
\subsection{Self-similar Zakharov-Craig-Sulem formulation}
For a wedge of slope $s$, the term in the energy corresponding to the surface tension is infinite. We define a relative Hamiltonian that subtracts the energy of the wedge. The integral over $[-R, R]$ is taken in the limit $R \to \infty$ since the integrand below may not be absolutely integrable  in our functional setting.  
 \begin{dfn}
  For a small slope $s$, define the Hamiltonian
  \[\mc{H}_s[f,\varphi] := \frac{1}{2}\int_{\R}\varphi\,G(f)\varphi \, d x + \lim_{R \to \infty}\int_{-R}^{R}\bigl(\sqrt{1+f'^2}-\sqrt{1 + s^2}\bigr) \, d x.\]
\end{dfn}
With this definition, we can write the self-similar system as a Hamiltonian system with respect to the modified Hamiltonian $\mc{H}_s$. Indeed, we can compute
\[\frac{\delta \mc{H}_s}{\delta \varphi} = G(f)\varphi, \quad \frac{\delta\mc{H}_s}{\delta f}
 =-\frac{f''}{(1+f'^2)^{3/2}}
   +\frac12\varphi'^2
   -\frac{(G(f)\varphi+f'\varphi')^2}{2(1+f'^2)}\]
so that the self-similar system can be written as
\begin{equation*}\label{eq:profile}\tag{$S_0$}
\begin{cases}
  -\displaystyle \frac{2}{3}(f-x f') = \frac{\delta \mc{H}_s}{\delta \varphi},\\
  -\displaystyle \frac{1}{3}\varphi + \frac{2}{3}x\varphi' = -\frac{\delta \mc{H}_s}{\delta f}.
\end{cases}
\end{equation*}

\subsection{Function spaces and regularization}
We fix $s$ to be the slope at infinity and we write 
\begin{equation*}
  f = s f_{\mathrm{ref}} + h, 
\end{equation*}
where $f_{\mathrm{ref}} \in C^{\infty}$ is an even, fixed function with $f_{\mathrm{ref}}(x)=|x|$ for $|x| \geq 2$. The function $h$ now represents the perturbation from the wedge, and it is the new unknown we solve for. We note here that our argument proves the existence of a solution just from perturbing the wedge $s f_{\mathrm{ref}}$. 

At the linear level, solutions to the system \eqref{eq:profile} present rapid oscillations at infinity, with behavior given by \eqref{eq:linear_oscillations}. This is one of the main difficulties of the problem. This oscillatory behavior causes a loss of control of higher derivatives of the solution since, at the linear level we already have that $f'' = O(|x|^{1/2})$. To overcome this difficulty, we introduce a regularization of the problem.
\begin{dfn}\label{dfn:regularized_system}
For $\epsilon > 0$, we say $(f,\varphi)$ is a solution of the regularized self-similar system if
\begin{equation*}\label{eq:profile_reg}\tag{$S_\epsilon$}
  \mathcal{F}_\epsilon(s, h,\varphi) = 0,
\end{equation*}
where
\begin{equation}\label{def:F_epsilon}
    \mathcal{F}_\epsilon(s, h, \varphi) :=
  \begin{pmatrix*}
    \displaystyle -\frac{2}{3} \left( f - x \p_x f \right) \\[8pt]
    \displaystyle -\frac{1}{3} \varphi + \frac{2}{3}x \p_x \varphi
  \end{pmatrix*}
    - (\mathbb{J}+\epsilon \mathbb{I}) \nabla_{(f, \varphi)} \mc{H}_s[f,\varphi],
    \quad \text{ and } \quad
  \mathbb{J} =
  \begin{pmatrix*}
    0 & 1 \\
    -1 & 0
  \end{pmatrix*}.
\end{equation}
We define $\mathcal{F}_0$ analogously, so that its zero set is equivalent to \eqref{eq:profile}.
\end{dfn}
\begin{remark}
  We are implicitly using that $f$ is a function of $s, h$ as $f= s f_\mathrm{ref} + h$ and $f_\mathrm{ref}$ is fixed.
\end{remark}

Since we are going to be working with the regularized problem, we define the following high regularity spaces. By this we mean that the oscillations described in \eqref{eq:linear_oscillations} will not belong to these spaces. We want to define $X, Y$ such that $\mathcal{F_\epsilon}: X \to Y$. 

\begin{definition}\label{def:X_Y_spaces}
  We use $H^{k}_e(\mathbb{R})$ to refer to the usual Sobolev space of even functions with $k$ derivatives. We define the spaces $X, Y$ as
  \begin{align*}
  X :&= \{(h,\varphi)\in H^3_e(\R)\times H^2_e(\R):xh',\,x\varphi'\in H^1(\R)\},\\
  \|(h,\varphi)\|_X &= \|h\|_{H^3}+\|\varphi\|_{H^2}+\|xh'\|_{H^1}+\|x\varphi'\|_{H^1},
  \end{align*}
  and
  \begin{align*}
  Y :&= H^1_e(\R) \times H^1_e(\R).
\end{align*}
\end{definition}

Throughout the rest of the argument, $\delta>0$ denotes a small threshold for the slope of the solution, independent of $\epsilon$ and $s$. We will decrease $\delta$ as needed during the argument.

\begin{lem}\label{lem:regularity_F}
  The map $\mathcal{F}_\epsilon(s, h, \varphi)$ is $C^1$ between the following spaces
  \begin{align*}
    \mathcal{F}_\epsilon(s,h,\varphi): U \subset \R \times X \to Y,
  \end{align*}
  where
  \begin{equation*}
  U := \left \{ (s, h, \varphi) \in \R \times X : \|s f_{\mathrm{ref}}' + h'\|_{L^\infty} < 2\delta \right \}.
  \end{equation*}
\end{lem}
\begin{proof}
  We go term by term in the definition of $\mathcal{F}_\epsilon$. Start with
  \begin{equation*}
    f - x f' = s(f_{\mathrm{ref}} - x f_{\mathrm{ref}}') + (h - x h').
  \end{equation*}
  Here $f_{\mathrm{ref}} - x f_{\mathrm{ref}}'$ is smooth and compactly supported, and $h$, $x h'$ belong to $H^1$ by the definition of $X$. This term is linear in $h$, so it is $C^1$ into $Y$. In the same way, the term $\varphi - x \varphi'$ is also $C^1$.

  Now we show that $(s, h, \varphi) \mapsto G(f)\varphi$ is $C^1$. The map $(s, h) \mapsto f' = sf'_{\mathrm{ref}} + h'$ is bounded and linear into $W^{1,\infty}$. Thus, the analytic dependence of $G(f): H^{2} \to H^{1}$ on $f' \in W^{1,\infty}$ of Lemma \ref{lem:DtN_analytic} implies that $(s, h) \mapsto G(f)$ is $C^1$ into $\mathcal{L}(H^2, H^1)$. We apply this lemma by taking $\delta$ sufficiently small. We conclude that $(s, h, \varphi) \mapsto G(f)\varphi$ is $C^1$.

  The last term in the equation is
  \begin{equation*}
    -\frac{f''}{(1+f'^2)^{3/2}} +\frac12\varphi'^2 - \frac{(G(f)\varphi+f'\varphi')^2}{2(1+f'^2)}.
  \end{equation*}
  We use that $f'', \varphi', G(f)\varphi \in H^{1}$ and that $H^{1}$ is an algebra. Moreover, $(1+f'^2)^{-1}$, $(1+f'^2)^{-3/2}$ are $C^{1}$ into the space $W^{1,\infty}$ and multiplication by a $W^{1,\infty}$ coefficient is bounded in $H^{1}$. In total, these terms are also $C^{1}$.
\end{proof}

The energy, stress, and bootstrap estimates in Sections \ref{sec:energy} and \ref{sec:uniform_bounds} are uniform in $0\leq\epsilon\leq1$. The compactness and degree arguments in Sections \ref{sec:compactness} and \ref{sec:existence_regularized} are done for a fixed value of $0 < \epsilon \leq 1$. 

\section{A priori energy and stress estimates}\label{sec:energy}

In this section, we prove a priori estimates for any solution $(f, \varphi)$ of the regularized self-similar system \eqref{eq:profile_reg} in the admissible class defined as follows.

\begin{dfn}\label{def:admissible}
 Given a sufficiently small $\delta > 0$, we say that a pair $(f, \varphi)$ is admissible if $f$ and $\varphi$ are even and the following conditions hold:
\begin{align}\label{sol_admisible}
\| f'\|_{L^\infty} \leq \delta, \quad f - s|x| \in H^1, \quad \varphi \in H^{1/2}.  
\end{align}
\end{dfn}

For the rest of the section, $\delta$ will be fixed, and the rest of the constants in this section are independent of $s$ and $\epsilon$. \par 
The first lemma says that, under the assumption of admissibility, the Hamiltonian $\mc{H}_s[f,\varphi]$ is well-defined and finite. Moreover, it controls both the norms $\|\varphi\|_{\dot{H}^{1/2}}$ and $\|(f - s|x|)'\|_{L^2}$.

\begin{lem}\label{lem:hamiltonian_estimates}
  Let $(f, \varphi)$ be an admissible pair. Then the Hamiltonian $\mc{H}_s[f,\varphi]$ is well-defined and finite. Moreover, we have the estimates
  \begin{equation}\label{eq:energy_hamiltonian_upper}\mc{H}_s[f,\varphi] + \frac{2sf(0)}{\sqrt{1 + s^2}} \leq \frac{1}{2}\int G(f)\varphi \, \varphi \, dx + \frac{1}{2} \int |(f - s|x|)'|^2 \, dx\end{equation}
  and 
  \begin{equation}\label{eq:energy_hamiltonian_lower}\mc{H}_s[f,\varphi] + \frac{2sf(0)}{\sqrt{1 + s^2}} \geq \frac{1}{2}\int G(f)\varphi \, \varphi \, dx + \frac{1}{2(1 + \delta^2)^{3/2}} \int |(f - s|x|)'|^2 \, dx.\end{equation}
\end{lem}

\begin{proof}
    We write
    \begin{equation}\label{eq:energy_surface_decomposition}\sqrt{1 + f'(x)^2} - \sqrt{1 + s^2} = \frac{s \, \text{sgn} \, x}{\sqrt{1 + s^2}}(f - s|x|)' + r(x).\end{equation}
    Then applying Taylor's theorem to the function $u \mapsto \sqrt{1 + u^2}$ about $s \, \text{sgn} \, x$ gives 
    \begin{equation}\label{eq:energy_surface_remainder}
    \frac{|(f(x) - s|x|)'|^2}{2(1+\delta^2)^{3/2}}
      \le r(x)\le\frac12|(f - s|x|)'|^2.
    \end{equation}
    Therefore 
    \[\mc{H}_s[f,\varphi] = \frac{1}{2} \int G(f)\varphi \, \varphi \, dx + \int r(x) \, dx + \lim_{R \to \infty}\int_{-R}^{R} \frac{s \, \text{sgn} \, x}{\sqrt{1 + s^2}}(f - s|x|)' \, dx.\]
    The integrand in the last term is even so it equals
    \[\frac{2s}{\sqrt{1 + s^2}}\lim_{R \to \infty}\int_{0}^{R}(f - s|x|)' \, dx = -\frac{2s}{\sqrt{1 + s^2}}f(0).\]
    Here we used that $f - s|x| \in H^1$ so that $f(R) - s|R| \to 0$ as $R \to \infty$. Therefore we have that
    \begin{equation}\label{eq:energy_hamiltonian_decomposition}\mc{H}_s[f,\varphi] = \frac{1}{2} \int G(f)\varphi \, \varphi \, dx + \int r(x) \, dx -\frac{2s}{\sqrt{1 + s^2}}f(0).\end{equation}
    The desired result follows.
\end{proof}

\begin{prop}[Energy balance]\label{prop:energy_balance} 
  Let $0\leq\epsilon\leq1$, and let $(f,\varphi)$ be an admissible pair with
  $(f,\varphi)\in C^2_{\mathrm{loc}}\times C^1_{\mathrm{loc}}$ and assume that \eqref{eq:profile_reg} holds pointwise.
  Then for $\epsilon>0$, both variational derivatives of $\mc{H}_s$ belong to $L^2$, and
  \begin{equation}\label{eq:energy_balance}\mc{H}_s[f,\varphi] + \frac{3\epsilon}{2}\left(
  \left\|\frac{\delta\mc{H}_s}{\delta f}\right\|_{L^2}^2
  +\left\|\frac{\delta\mc{H}_s}{\delta\varphi}\right\|_{L^2}^2\right) = 0.\end{equation}
  For $\epsilon=0$, the conclusion is $\mc{H}_s[f,\varphi]=0$. In particular, we have that 
  \begin{equation}\label{eq:energy_high_frequency}\|(f - s|x|)'\|_{L^2}^2 + \|\varphi\|_{\dot{H}^{1/2}}^2 + \epsilon \left(
  \left\|\frac{\delta\mc{H}_s}{\delta f}\right\|_{L^2}^2
  +\left\|\frac{\delta\mc{H}_s}{\delta\varphi}\right\|_{L^2}^2\right) \leq C\frac{sf(0)}{\sqrt{1 + s^2}}.\end{equation}
\end{prop}

\begin{rem}
  A key point of this proposition is that it provides a balance between height and the high-frequency behavior of the solution.
\end{rem}

\begin{rem}
  In particular, the proposition implies that $\text{sgn}(f(0)) = \text{sgn}(s)$.
\end{rem}

\begin{rem}
  To motivate the energy balance identity \eqref{eq:energy_balance}, we consider the following formal calculation. In the standard variables, the regularized system takes the form
  \[
    \begin{aligned}
    \eta_t
    &=
    \frac{\delta\mc{H}_s}{\delta\phi}
    +
    \epsilon (T - t)^{1/3}
    \frac{\delta\mc{H}_s}{\delta\eta},\\
    \phi_t
    &=
    -\frac{\delta\mc{H}_s}{\delta\eta}
    +
    \epsilon (T - t)^{-1/3}
    \frac{\delta\mc{H}_s}{\delta\phi}.
    \end{aligned}
  \]
  Taking the time derivative of the Hamiltonian then gives 
  \[\frac{d}{dt} \mc{H}_s[\eta(t),\phi(t)] = \epsilon (T - t)^{1/3}\left\|
  \frac{\delta\mc{H}_s}{\delta\eta}
  \right\|_{L^2}^{2} + \epsilon (T - t)^{-1/3} \left\|\frac{\delta\mc{H}_s}{\delta\phi}\right\|_{L^2}^{2}.\]
  (which in particular implies conservation of the Hamiltonian for $\epsilon=0$). On the other hand, with our self-similar ansatz, we also have that
  \begin{equation}\label{eq:energy_scaling}\mc{H}_s[\eta(t),\phi(t)] = (T - t)^{2/3}\mc{H}_s[f,\varphi].\end{equation}
  By taking the time derivative of this expression, and using the identities
  \[\left\|
  \frac{\delta\mc{H}_s}{\delta\eta} \right\|_{L^2}^{2} = (T - t)^{-2/3}\left\| \frac{\delta\mc{H}_s}{\delta f} \right\|_{L^2}^{2}, \quad \left\|\frac{\delta\mc{H}_s}{\delta\phi}\right\|_{L^2}^{2} = \left\| \frac{\delta\mc{H}_s}{\delta\varphi} \right\|_{L^2}^{2}\]
  (which follow from scaling), we arrive at the desired energy identity 
  \[\mc{H}_s[f,\varphi] + \frac{3\epsilon}{2}\left(
  \left\|\frac{\delta\mc{H}_s}{\delta f}\right\|_{L^2}^2
  +\left\|\frac{\delta\mc{H}_s}{\delta\varphi}\right\|_{L^2}^2\right) = 0.\]
  The point of \Cref{prop:energy_balance} is to rigorously justify this formal calculation for admissible solutions $(f, \varphi)$.
\end{rem}

Before proving \Cref{prop:energy_balance}, we first prove a lemma that will be needed to justify an integration by parts in the proof. It says that the velocity $\nabla\Psi$ of the finite-energy harmonic extension of $\varphi$ extends continuously to the graph uniformly on compact sets.

\begin{lem}\label{lem:velocity_continuity}
  Let $0 \leq \epsilon \leq 1$ and let $(f, \varphi) \in C^2_{\mathrm{loc}}\times C^1_{\mathrm{loc}}$ be an admissible pair satisfying \eqref{eq:profile_reg}. Then the finite-energy harmonic extension of $\varphi$ satisfies
  \[\Psi \in C^1_{\textup{loc}}(\ov{\Omega_f}).\]
  In particular, the velocity $\nabla\Psi$ extends continuously to the graph $\{(x, f(x)) : x \in \R\}$ uniformly on compact sets.
\end{lem}

\begin{proof}
  The idea is to show that $\Psi$ satisfies an elliptic equation with $C^1$ Neumann data. In the case that $\epsilon = 0$, this can be seen immediately from $\Delta \Psi = 0$ and the kinematic equation:
  \[G(f)\varphi = -\frac{2}{3}(f - xf').\]
  It follows that $\Psi$ satisfies
  \[\begin{cases}
    \Delta \Psi = 0 & \text{in } \Omega_f,\\
    \nabla \Psi \cdot \nu = -\frac{2(f - xf')}{3\sqrt{1 + f'^2}} & \text{on } y = f(x)
  \end{cases}, \quad \nu = \frac{(-f', 1)}{\sqrt{1 + f'^2}}\]
  and the Neumann data has enough regularity to give the conclusion. \par
  For $\epsilon > 0$, the argument is similar, but the Neumann data is different. Indeed, by adding the kinematic equation with $\epsilon$ times the dynamic equation, the terms $\epsilon \frac{\delta \mc{H}_s}{\delta f}$ cancel and we obtain
    \begin{equation}\label{eq:energy_neumann_data}
      G(f)\varphi
      =
      -\frac{2(f-xf')+\epsilon\varphi-2\epsilon x\varphi'}
            {3(1+\epsilon^2)}.
    \end{equation}
  We move the derivatives in \eqref{eq:energy_neumann_data} to the left to obtain
  \[G(f) \varphi - \frac{2 \epsilon x}{3(1 + \epsilon^2)}\varphi ' = -\frac{2(f-xf')+\epsilon\varphi}{3(1+\epsilon^2)}.\]
  Define the matrix
  \[a(x) = \frac{2\epsilon x}{3(1 + \epsilon^2)}, \quad A(x) = \begin{bmatrix}
    1 &  a(x) \\
    - a(x) & 1
  \end{bmatrix}.\]
  The Neumann data can be expressed in terms of the matrix $A$. Indeed, using 
  \[G(f)\varphi = \sqrt{1 + (f')^2} \, \nabla \Psi \cdot \nu, \, \quad \varphi' = \Psi_x + f' \Psi_y,\]
  we obtain
  \[A\nabla \Psi \cdot \nu = - \frac{2(f - xf') + \epsilon \varphi}{3(1 + \epsilon^2)\sqrt{1 + (f')^2}}\]
  on the boundary $y = f(x)$. Using the fact that $\Psi$ is harmonic, we obtain
  \[
    \operatorname{div}(A\nabla\Psi)
    =
    \partial_x(\Psi_x+a\Psi_y)
    +\partial_y(-a\Psi_x+\Psi_y)
    =
    a'\Psi_y.
  \]
  It follows that $\Psi$ satisfies the elliptic equation with $C^1_{\textup{loc}}$ Neumann data:
  \begin{equation}\label{eq:energy_conormal_problem}
    \begin{cases}
      \operatorname{div}(A\nabla\Psi) - a'\Psi_y = 0 & \text{in } \Omega_f,\\
      A\nabla\Psi\cdot\nu = g & \text{on } \partial \Omega_f
    \end{cases}
  \end{equation}
  where
  \[g = -\frac{2(f-xf')+\epsilon\varphi}{3(1+\epsilon^2)\sqrt{1+(f')^2}}.\]
  Since $\xi \cdot A \xi = |\xi|^2$, the matrix is uniformly elliptic, and we can appeal to elliptic theory to give the conclusion. There are multiple ways to get the conclusion, but the most convenient way is to use \cite[Theorem~1.7]{DongLeeKim2020}, which gives that $\Psi \in C^1_{\textup{loc}}(\ov{\Omega_f})$ as wanted.
  \footnote{The result of \cite{DongLeeKim2020} is stated for bounded domains, but a standard cutoff argument allows us to apply the result to get the local conclusion on $\Omega_f$. Indeed, the hypotheses of \cite[Theorem~1.7]{DongLeeKim2020} are satisfied since $A$ is uniformly elliptic and locally bounded, $\Omega_f$ is locally $C^2$, and $g \in C^1_{\textup{loc}}$.}
\end{proof}

\begin{proof}[Proof of \Cref{prop:energy_balance}]
  Recall that the self-similar system \eqref{eq:profile_reg} can be written as
  \begin{equation}\label{eq:energy_profile_system}
    \begin{aligned}
    -\frac{2}{3}(f-xf')
    &=
    \frac{\delta\mc{H}_s}{\delta\varphi}
    +
    \epsilon
    \frac{\delta\mc{H}_s}{\delta f},\\
    -\frac{1}{3}\varphi+\frac{2}{3}x\varphi'
    &=
    -\frac{\delta\mc{H}_s}{\delta f}
    +
    \epsilon
    \frac{\delta\mc{H}_s}{\delta\varphi}.
    \end{aligned}
  \end{equation}
  Where the variational derivatives are given by
  \begin{equation}\label{eq:energy_variational_derivatives}
    \begin{aligned}
    \frac{\delta\mc{H}_s}{\delta\varphi}
    =G(f)\varphi,\quad
    \frac{\delta\mc{H}_s}{\delta f}
    =-\frac{f''}{(1+f'^2)^{3/2}}
    +\frac12(\varphi')^2
    -\frac{(G(f)\varphi+f'\varphi')^2}{2(1+f'^2)}.
    \end{aligned}
  \end{equation}
  Define the vector field 
  \[Q(z) = \frac{1}{2}\Psi \nabla \Psi - Tz,\]
  with $T$ the tensor
\begin{equation}\label{eq:T_stress}
T = \nabla \Psi \otimes \nabla \Psi - \frac{1}{2}|\nabla \Psi|^2 I.  
  \end{equation}
  Then since $\Delta \Psi = 0$, we have that $\operatorname{div} T = 0$ and since the trace of $T$ is zero, we have that $\operatorname{div} (Tz) = 0$ as well. Therefore we have that
  \[\operatorname{div} Q = \frac{1}{2}|\nabla \Psi|^2.\]
  By multiplying the kinematic equation by $\frac{\delta \mc{H}_s}{\delta f}$ and the dynamic equation by $\frac{\delta \mc{H}_s}{\delta \varphi}$, and summing the two equations in \eqref{eq:energy_profile_system}, we arrive at the identity
  \begin{equation}\label{eq:energy_pointwise_balance}
    (f-xf')\frac{\delta\mc{H}_s}{\delta f}
    +\left(\frac12\varphi-x\varphi'\right)
    \frac{\delta\mc{H}_s}{\delta\varphi}
    +\frac{3\epsilon}{2}\left[
    \left(\frac{\delta\mc{H}_s}{\delta f}\right)^2
    +\left(\frac{\delta\mc{H}_s}{\delta\varphi}\right)^2
    \right]=0.
  \end{equation}
  Let $\chi_R(x) = \sq{\chi}(x, f(x))$ where $\sq{\chi} \in C^\infty_c(\R^2)$ is a smooth cutoff function with $0 \leq \sq{\chi} \leq 1$,  $\sq{\chi}(x,y) = 1$ for $|(x,y)| \leq R$ and $\sq{\chi}(x,y) = 0$ for $|(x,y)| \geq 2R$, and $\|\nabla \sq{\chi}_R\|_{L^\infty} \lesssim R^{-1}$. Thus if $\chi_R(x) \neq 0$ then $|x| + |f(x)| \lesssim R$. We multiply the above identity \eqref{eq:energy_pointwise_balance} by $\chi_R$ and integrate to obtain
  \begin{equation}\label{eq:energy_cutoff_balance}\int \chi_R(x) \left[(f-xf')\frac{\delta\mc{H}_s}{\delta f}
    +\left(\frac12\varphi-x\varphi'\right)
    \frac{\delta\mc{H}_s}{\delta\varphi}\right] \, dx + \frac{3\epsilon}{2}\int \chi_R(x) \left[
    \left(\frac{\delta\mc{H}_s}{\delta f}\right)^2
    +\left(\frac{\delta\mc{H}_s}{\delta\varphi}\right)^2
    \right] \, dx = 0\end{equation}
    the integrals make sense because of the local regularity of $(f, \varphi)$. Let us denote the first integral in \eqref{eq:energy_cutoff_balance} by $I_R$. We claim that 
    \[\lim_{R \to \infty}I_R = \mc{H}_s[f,\varphi].\]
    Indeed, observe that we have the identity
    \begin{equation}\label{eq:energy_boundary_flux}
    \begin{aligned}
      (f-xf')\frac{\delta\mc{H}_s}{\delta f}
    +\left(\frac12\varphi-x\varphi'\right)
    \frac{\delta\mc{H}_s}{\delta\varphi} & = Q(x, f(x)) \cdot (-f', 1) \\
    & + \sqrt{1 + (f')^2} - \sqrt{1 + s^2} \\
    & - \partial_x \left[\frac{x + ff'}{\sqrt{1 + (f')^2}} - x\sqrt{1 + s^2}\right].
    \end{aligned}
    \end{equation}
    It follows from \eqref{eq:energy_boundary_flux} that we can write $I_R$ as
    \begin{align*}
      I_R & = \int_{\R}\chi_R Q(x,f(x)) \cdot (-f', 1) \, dx + \int_{\mathbb R}\chi_R
    \left(\sqrt{1+f'^2}-\sqrt{1+s^2}\right)\,dx \\
      & + \int_{\mathbb R}\chi_R' \left[\frac{x + ff'}{\sqrt{1 + (f')^2}} - x\sqrt{1 + s^2}\right]\,dx.
    \end{align*}
    Write
    \[I_{R,1} = \int_{\R}\chi_R Q(x,f(x)) \cdot (-f', 1) \, dx + \int_{\mathbb R}\chi_R
    \left(\sqrt{1+f'^2}-\sqrt{1+s^2}\right)\,dx\]
    and 
    \[I_{R,2} = \int_{\mathbb R}\chi_R' \left[\frac{x + ff'}{\sqrt{1 + (f')^2}} - x\sqrt{1 + s^2}\right]\,dx.\]
    We claim that $I_{R,1} \to \mc{H}_s[f,\varphi]$  and $I_{R,2} \to 0$. \par
    For $I_{R,1}$, we use $Q = \frac{1}{2}\Psi \nabla \Psi - Tz$:
    \begin{align*}
      \int_{\R}\chi_R Q(x,f(x)) \cdot (-f', 1) \, dx & = \frac{1}{2}\int G(f)\varphi \, \varphi \chi_R \, dx - \int \chi_R (Tz)(x,f(x)) \cdot (-f', 1) \, dx.
    \end{align*}
    The first term converges to $\frac{1}{2}\int G(f)\varphi \, \varphi \, dx$ since $\varphi \chi_R \to \varphi$ in $H^{1/2}$. For the second term, we use the divergence theorem (which is justified by \Cref{lem:velocity_continuity}) to write it as 
    \begin{equation}\label{eq:energy_stress_flux}\int \chi_R (Tz)(x,f(x)) \cdot (-f', 1) \, dx = \int_{\Omega_f} \nabla \sq{\chi}_R \cdot (Tz) \, dz.\end{equation}
    Since $|\nabla \sq{\chi}_R||Tz| \lesssim |\nabla \Psi|^2$, and the integrand is supported in the annulus $\{z \in \R^2 : R \leq |z| \leq 2R\}$, this integral in \eqref{eq:energy_stress_flux} converges to zero. For the surface term in $I_{R,1}$, we have that
    \begin{align*}
      \int \chi_R (\sqrt{1 + f'^2} - \sqrt{1 + s^2}) \, dx & = \int [\chi_R - 1_{[-R,R]}] (\sqrt{1 + f'^2} - \sqrt{1 + s^2}) \, dx \\
      & + \int_{-R}^{R}(\sqrt{1 + f'^2} - \sqrt{1 + s^2}) \, dx.
    \end{align*}
    We claim that the first term on the right-hand side converges to zero. Indeed, by using the same decomposition as in \eqref{eq:energy_surface_decomposition}, we can write
    \[\sqrt{1 + f'(x)^2} - \sqrt{1 + s^2} = \frac{s \, \text{sgn} \, x}{\sqrt{1 + s^2}}(f - s|x|)' + r(x)\]
    where $r \in L^1$. Since $\chi_R$ is equal to $1$ for $|x| \leq R / 2$, we have that
    \[\left|\int [\chi_R - 1_{[-R,R]}] r(x) \, dx\right| \lesssim \int_{|x| \approx R} |r(x)| \, dx \to 0.\]
    For the other term, we integrate by parts and use that $f$ is even to obtain
    \[\int [\chi_R - 1_{[-R,R]}] \operatorname{sgn}(x) (f - s|x|)' \, dx = - \int \chi_R' \operatorname{sgn}(x)(f - s|x|) \, dx - 2(f(R) - sR).\]
    The first term converges to zero since $|\chi_R'| \lesssim R^{-1}$ and $f - s|x| \in L^2$. The second term converges to zero since $f - s|x| \in H^1$.
    It follows that $I_{R,1} \to \mc{H}_s[f,\varphi]$ as $R \to \infty$. \par 
    For $I_{R,2}$, we write
    \[\frac{x + ff'}{\sqrt{1 + (f')^2}} - x\sqrt{1 + s^2} = \frac{f'}{\sqrt{1 + f'^2}}(f - s|x|) + x \left[\sqrt{1 + f'^2} - \sqrt{1 + s^2} - \frac{f'}{\sqrt{1 + f'^2}}(f - s|x|)'\right].\]
    The second term is the Taylor approximation of the function $u \mapsto \sqrt{1 + u^2}$ around $u = f'$ evaluated at $s \operatorname{sgn} x$. It follows that 
    \[\left|\frac{x + ff'}{\sqrt{1 + (f')^2}} - x\sqrt{1 + s^2}\right| \lesssim |f'||f - s|x|| + |x||(f - s|x|)'|^2.\]
    By integrating and using that $|\chi_R'| \lesssim R^{-1}$, we have that
    \begin{equation}\label{eq:energy_surface_cutoff}|I_{R,2}| \lesssim R^{-1} \int_{|x| \lesssim R}|f - s|x|| \, dx + \int_{|x| \gtrsim R}|(f - s|x|)'|^2 \, dx.\end{equation}
    This term in \eqref{eq:energy_surface_cutoff} converges to zero under the assumptions that $f - s|x| \in H^1$. \par We have therefore established that $I_R \to \mc{H}_s[f,\varphi]$ as $R \to \infty$. Letting $R \to \infty$ in the identity \eqref{eq:energy_cutoff_balance},
    \[I_R + \frac{3\epsilon}{2}\int \chi_R(x) \left[
    \left(\frac{\delta\mc{H}_s}{\delta f}\right)^2
    +\left(\frac{\delta\mc{H}_s}{\delta\varphi}\right)^2
    \right] \, dx = 0\]
    and using that $I_R$ converges to the finite number $\mc{H}_s[f,\varphi]$ implies that when $\epsilon > 0$, the variational derivatives are in $L^2$ as wanted. \par
    The last assertion regarding the high-frequency control follows immediately from \Cref{lem:hamiltonian_estimates}.
    
\end{proof}

The next lemma exploits a balance between the low-frequency and high-frequency behavior of the solution to control the solution in the low frequency spaces $H^{-1/2}$ and $H^{-2}$. For instance, the kinematic equation for $\epsilon = 0$ says
\[-\frac{2}{3}(f - xf') = G(f)\varphi.\]
We will show that the left-hand side of this equation controls the low-frequency behavior of $f - s|x|$ in $H^{-1/2}$, and the right-hand side is controlled by $\|G(f)\varphi\|_{H^{-1/2}} \lesssim \|\varphi\|_{\dot{H}^{1/2}}$ which is in turn controlled by the energy balance identity.
 
\begin{lem}[Control of low frequencies]\label{lem:low_freq}
Under the hypotheses of Proposition \ref{prop:energy_balance}, we have that 
\begin{equation}\label{eq:energy_low_f}\|f - s|x|\|_{H^{-1/2}}^2  \leq C \frac{sf(0)}{\sqrt{1 + s^2}}\end{equation}
and
\begin{equation}\label{eq:energy_low_phi}\|\varphi\|_{H^{-2}} \lesssim |s| + \|(f - s|x|)'\|_{L^2} + \|\varphi\|_{\dot{H}^{1/2}}^2 + \epsilon \|\varphi\|_{\dot{H}^{1/2}}.\end{equation}
\end{lem}

\begin{proof}
The proof relies on the fact that, if $T_\lambda$ denotes the operator
\[T_\lambda(g) = \lambda g - xg',\]
then for $r \geq 0$ and $\lambda>0$, $T_\lambda$ satisfies the estimate
\begin{equation}\label{eq:energy_dilation_coercivity}\|g\|_{H^{-r}} \leq \frac{1}{\sqrt{\lambda(\lambda + 1)}} \|T_\lambda (g)\|_{H^{-r}}.\end{equation}
This can be seen by a straightforward Fourier transform argument: If $g$ is a Schwartz function, then we have that
\begin{align*}
  \|T_\lambda(g)\|_{H^{-r}}^2 & = \int_{\R} \ang{\xi}^{-2r}|(\lambda + 1) \hat{g}(\xi) + \xi \hat{g}'(\xi)|^2 \, d\xi.
\end{align*}
Expanding the square gives
\begin{align*}
  |(\lambda + 1) \hat{g}(\xi) + \xi \hat{g}'(\xi)|^2 & = (\lambda + 1)^2 |\hat{g}(\xi)|^2 + (\lambda + 1)\xi \partial_\xi(|\widehat{g}(\xi)|^2) + \xi^2 |\hat{g}'(\xi)|^2 \\
  & \geq (\lambda + 1)^2 |\hat{g}(\xi)|^2 + (\lambda + 1)\xi \partial_\xi(|\widehat{g}(\xi)|^2).
\end{align*}
Therefore, by integrating by parts, we have that
\begin{equation}\label{eq:Tlambda_low}
\begin{aligned}
  \|T_\lambda(g)\|_{H^{-r}}^2
  &\geq (\lambda+1)^2\|g\|_{H^{-r}}^2
  -(\lambda+1)\int_{\R}\partial_\xi\bigl(\xi\ang{\xi}^{-2r}\bigr)|\hat{g}(\xi)|^2\,d\xi \\
  &=\lambda(\lambda+1)\|g\|_{H^{-r}}^2
  +2r(\lambda+1)\int_{\R}\xi^2\ang{\xi}^{-2r-2}|\hat{g}(\xi)|^2\,d\xi \\
  &\geq \lambda (\lambda+1)\|g\|_{H^{-r}}^2.
\end{aligned}
\end{equation}
The estimate \eqref{eq:energy_dilation_coercivity} then follows from approximating a general $g$ by Schwartz functions. Now that we have established this low-frequency estimate, we can apply it to function $g = f - s|x|$. By applying the above estimate \eqref{eq:energy_dilation_coercivity} with $\lambda = 1$ and $r = 1/2$, we have that
\begin{align*}
  \|f - s|x|\|_{H^{-1/2}} & \lesssim \|f - xf'\|_{H^{-1/2}} \\
  & \lesssim \|G(f) \varphi\|_{H^{-1/2}} + \epsilon \left\|\frac{\delta \mc{H}_s}{\delta f}\right\|_{H^{-1/2}}.
\end{align*}
Using the estimates for the Dirichlet-to-Neumann operator (\Cref{lem:DtN_H12_bound}), we have that
\[\|G(f)\varphi\|_{H^{-1/2}} \leq \|G(f)\varphi\|_{\dot{H}^{-1/2}} \lesssim \|\varphi\|_{\dot{H}^{1/2}}.\]
Also, the variational derivative term can be bounded by the $L^2$ norm, which gives 
\begin{align*}
  \|f - s|x|\|_{H^{-1/2}}^2 & \lesssim \|\varphi\|_{\dot{H}^{1/2}}^2 + \epsilon \left\|\frac{\delta \mc{H}_s}{\delta f}\right\|_{L^2}^2
\end{align*}
(where the second term is only present if $\epsilon > 0$). By applying \Cref{prop:energy_balance}, we arrive at the bound
\[\|f - s|x|\|_{H^{-1/2}}^2 \leq C \frac{sf(0)}{\sqrt{1 + s^2}} .\]
For $\varphi$, observe that we can express the dynamic equation as 
\begin{equation}\label{eq:energy_dynamic_stress}-\frac{1}{3}\varphi + \frac{2}{3}x\varphi' = \left(\frac{f'}{\sqrt{1 + (f')^2}} - \int_{-\infty}^{f(x)}\Psi_x \Psi_y \, dy\right)' + \epsilon G(f)\varphi.\end{equation}
Therefore, by applying the low-frequency estimate \eqref{eq:energy_dilation_coercivity} with $\lambda = 1/2$ and $r = 2$, we have that
\[\|\varphi\|_{H^{-2}} \lesssim \left\|\left(\frac{f'}{\sqrt{1 + (f')^2}}\right)'\right\|_{H^{-2}} + \left\|\left(\int_{-\infty}^{f(x)}\Psi_x \Psi_y \, dy\right)'\right\|_{H^{-2}} + \epsilon \|G(f)\varphi\|_{H^{-2}}.\]
To bound the first term, we write
\[\frac{f'}{\sqrt{1 + (f')^2}} = \frac{f'}{\sqrt{1 + (f')^2}} - \frac{(s|x|)'}{\sqrt{1 + s^2}} + \frac{(s|x|)'}{\sqrt{1 + s^2}}.\]
Using that $|x|'' = 2\delta_0 \in H^{-2}$, we obtain that the derivative of the last term can be bounded by $|s|$ in $H^{-2}$. For the other term, since the function $u \mapsto u(1 + u^2)^{-1/2}$ has bounded derivative, and $\|g'\|_{H^{-2}} \lesssim \|g\|_{L^2}$, we can bound
\[\left\|\left(\frac{f'}{\sqrt{1 + (f')^2}}\right)'\right\|_{H^{-2}} \lesssim |s| + \left\|\frac{f'}{\sqrt{1 + (f')^2}} - \frac{(s|x|)'}{\sqrt{1 + s^2}}\right\|_{L^2} \lesssim |s| + \|(f - s|x|)'\|_{L^2}.\]
To bound the second term, we use that $\|g'\|_{H^{-2}} \lesssim \|g\|_{L^1}$ and Lemma \ref{lem:H12_trace} to get 
\[\left\|\left(\int_{-\infty}^{f(x)}\Psi_x \Psi_y \, dy\right)'\right\|_{H^{-2}} \lesssim \int_{\Omega_f}|\Psi_x \Psi_y| \, dxdy \lesssim \int_{\Omega_f}|\nabla \Psi|^2 \, dxdy \lesssim \|\varphi\|_{\dot{H}^{1/2}}^2.\]
Finally, for the last term, we have
\[\|G(f)\varphi\|_{H^{-2}} \lesssim \|G(f)\varphi\|_{\dot{H}^{-1/2}} \lesssim \|\varphi\|_{\dot{H}^{1/2}}.\]
It follows that
\[\|\varphi\|_{H^{-2}} \lesssim |s| + \|f - s|x|\|_{\dot{H}^{1}} + \|\varphi\|_{\dot{H}^{1/2}}^2 + \epsilon \|\varphi\|_{\dot{H}^{1/2}}\]
as wanted.
\end{proof}

Finally, by combining the high-frequency control from \Cref{prop:energy_balance} with the low-frequency control from \Cref{lem:low_freq}, we obtain uniform energy bounds for admissible solutions $(f, \varphi)$.

\begin{prop}[Uniform energy bounds]\label{prop:energy_bounds}
  Under the hypotheses of Proposition \ref{prop:energy_balance},
  \begin{equation}\label{eq:energy_bound_01}
    \|f - s|x|\|_{H^1} + \|\varphi\|_{H^{1/2}} +\|\nabla \Psi\|_{L^2(\Omega_f)} \leq C |s|
  \end{equation}
  and, for $\epsilon>0$,
  \begin{equation}\label{eq:energy_uniform_dissipation}
    \epsilon\left(
  \left\|\frac{\delta\mc{H}_s}{\delta f}\right\|_{L^2}^2
  +\left\|\frac{\delta\mc{H}_s}{\delta\varphi}\right\|_{L^2}^2\right) \leq C s^2.
  \end{equation}
  for some $C>0$, independent of $\epsilon$ and $s$. At $\epsilon=0$, the dissipative term is defined to be zero.
\end{prop}

\begin{proof}
  The proof follows from combining the high-frequency control from \Cref{prop:energy_balance} with the low-frequency control from \Cref{lem:low_freq}. Indeed, for $f - s|x|$, we can bound 
  \[\|f - s|x|\|_{H^1}^2 \lesssim \|f - s|x|\|_{H^{-1/2}}^2 + \|f - s|x|\|_{\dot{H}^{1}}^2\]
  which follows from taking the Fourier transform and applying the inequality 
  \[\ang{\xi}^2 \lesssim \ang{\xi}^{-1} + |\xi|^2.\]
  It then follows from \Cref{prop:energy_balance} and \Cref{lem:low_freq} that we have 
  \[\|f - s|x|\|_{H^1}^2 \lesssim |s||f(0)| = |s||f(0) - s|0|| \lesssim |s| \|f - s|x|\|_{H^1}.\]
  Dividing both sides by $\|f - s|x|\|_{H^1}$ gives the desired bound
  \[\|f - s|x|\|_{H^1} \lesssim |s|.\]
  In particular, this implies that $|f(0)| \lesssim |s|$.
  For the estimate on $\varphi$, \Cref{prop:energy_balance} gives
  \[\|\varphi\|_{\dot{H}^{1/2}}^2 \lesssim |s||f(0)| \lesssim s^2\]
 and \Cref{lem:low_freq} gives
  \begin{align*}
    \|\varphi\|_{H^{-2}} & \lesssim |s| + \|(f - s|x|)'\|_{L^2} + \|\varphi\|_{\dot{H}^{1/2}}^2 + \epsilon\|\varphi\|_{\dot{H}^{1/2}} \\
    & \lesssim |s| + s^2 + \epsilon |s| \lesssim |s|.
  \end{align*}
  Finally, the inequality $\ang{\xi} \lesssim \ang{\xi}^{-4} + |\xi|$ allows us to combine the low-frequency and high-frequency estimates to obtain
  \[\|\varphi\|_{H^{1/2}}^2 \lesssim \|\varphi\|_{H^{-2}}^2 + \|\varphi\|_{\dot{H}^{1/2}}^2 \lesssim s^2.\]
  The bound on the velocity $\nabla \Psi$ is then a consequence of Lemma \ref{lem:H12_trace}. 
  The final assertion \eqref{eq:energy_uniform_dissipation} of the proposition regarding the variational derivatives follows immediately from \Cref{prop:energy_balance} and $|f(0)| \lesssim |s|$.
  
\end{proof}

Roughly speaking, the previous energy estimate guarantees that the kinetic energy and the difference between the surface tension energy of the actual interface and that of the wedge are small for admissible solutions \eqref{def:admissible}.
In fact, 
the energy estimate in Proposition \eqref{prop:energy_bounds} states that the slope of the interface is controlled in $L^2$ by that of the wedge $s|x|$. Similarly,   the velocity in $L^2$ is controlled by  $s$. 
However, these estimates are too weak to control the $L^\infty$ norm of the slope for the whole profile in terms of $s$. That is, the smallness of the energy is not enough to prevent that large oscillations at small scales happen as $x\to\infty$. 

For this reason, we study the balance of momentum at a point on the interface, which is a pointwise identity. Crucially, the balance contains the slope $f'(x)$, through the surface tension force.
More specifically, we study the balance of momentum in the direction of the wedge, $(1,s)$, across a vertical cut. The choice of direction is not fundamental, we could have computed the balance of horizontal momentum, but it simplifies the proof of Proposition \ref{prop:stress_bounds} below.

\begin{dfn}\label{def:relative_stress}
  For an admissible pair $(f, \varphi)$ and $0 \leq \epsilon \leq 1$, define the  stress $\mathcal{R}_{s, \epsilon}(x)$ by
    \begin{equation}\label{eq:stress_definition}
    \begin{aligned}
    \mathcal R_{s,\epsilon}(x) & =
    \sqrt{1+s^2}
    -\frac{1+s f'(x)}{\sqrt{1+f'(x)^2}}
    \\
    &+\int_{-\infty}^{f(x)}
    \left[
    \frac{\Psi_x(x,y)^2-\Psi_y(x,y)^2}{2}
    +s\,\Psi_x(x,y)\Psi_y(x,y)
    \right]\,dy
    \\
    &-\frac{2(f(x)-sx)\varphi(x)}
            {3(1+\epsilon^2)} - \frac{\epsilon\left(
    2(f(x)-sx)^2+\varphi(x)^2
    \right)}
    {6(1+\epsilon^2)}.
  \end{aligned}
    \end{equation}
\end{dfn}

\begin{rem}
 The definition of $\mathcal{R}_{s,\epsilon}$ is motivated by the identity in Lemma \ref{lem:stress_identity} below. That identity can be formally obtained by computing the variation of momentum in a region enclosed by the interface and two vertical cuts, and then passing to self-similar coordinates.
 In fact, the first two terms in $\mathcal{R}_{s,\epsilon}$ have  clear physical interpretations: on the interface, the force due to the pressure is obtained by integration of the surface tension (normalized to be zero for the wedge), while the second term represents the pressure contribution along the vertical cuts together with the momentum crossing those cuts (in the chosen wedge direction). In fact, with $\epsilon=0$ for clarity of exposition, from $\partial_t \Phi+\frac12|\nabla\Phi|^2+p=0$, we can see that
 \begin{equation*}
     \Phi_X (\Phi_X,\Phi_Y)\cdot(1,s)+p=\Phi_X(\Phi_X+s\Phi_Y)+p=-\Phi_t+\frac{\Phi_X^2-\Phi_Y^2}{2}+s\Phi_X\Phi_Y.
 \end{equation*}
 The term $\Phi_t$ and moving to self-similar coordinates cause the last term in $\mathcal{R}_{s,\epsilon}$. Although $\mathcal{R}_{s,\epsilon}$ can be thought more as a force (in fact, the term \textit{flow force} has been used in the water waves literature for related quantities, see e.g. \cite{Benjamin-Brooke95, Basu-Martin20}),  we call it stress due to the integral term, which integrates the stress across a vertical cut in the wedge direction, $(1,s)^T T(1,0)$, with $T$ the stress tensor associated to the kinetic energy \eqref{eq:T_stress}.
\end{rem}

\begin{lemma}[Stress balance]\label{lem:stress_identity}
  Under the hypotheses of \Cref{prop:energy_balance}, the  stress satisfies the identity
  \begin{equation}\label{eq:stress_identity}\mathcal R_{s,\epsilon}'(x) = -\frac{(f'(x)-s)\varphi(x)}{1+\epsilon^2} - \frac{2\epsilon x}{3(1+\epsilon^2)}\left[(f'(x)-s)^2+\varphi'(x)^2\right].\end{equation}
\end{lemma}

\begin{proof}
  For the fluid term, harmonicity gives
  \[\partial_x\left(\frac{\Psi_x^2-\Psi_y^2}{2}+s\Psi_x\Psi_y\right)
    = \partial_y\left(-\Psi_x\Psi_y+\frac{s}{2}(\Psi_x^2-\Psi_y^2)\right).\]
  Therefore by differentiating, we have
  \begin{align*}
    &\partial_x\int_{-\infty}^{f(x)}
    \left(\frac{\Psi_x^2-\Psi_y^2}{2}+s\Psi_x\Psi_y\right)\,dy \\
    &\quad= f'\left.\left[\frac{\Psi_x^2-\Psi_y^2}{2}+s\Psi_x\Psi_y\right]\right|_{y=f(x)}
    +\int_{-\infty}^{f(x)}\partial_y\left(-\Psi_x\Psi_y+\frac{s}{2}(\Psi_x^2-\Psi_y^2)\right)\,dy \\
    &\quad= \left.\left[\frac{f'+s}{2}(\Psi_x^2-\Psi_y^2)+(sf'-1)\Psi_x\Psi_y\right]\right|_{y=f(x)}.
  \end{align*}
  On the boundary, we have
  \[\varphi'=(\Psi_x+f'\Psi_y)|_{y=f(x)}, \qquad
    G(f)\varphi=(\Psi_y-f'\Psi_x)|_{y=f(x)}.\]
  In particular,
  \[G(f)\varphi\,\varphi'=\left.\left[-f'(\Psi_x^2-\Psi_y^2)+(1-f'^2)\Psi_x\Psi_y\right]\right|_{y=f(x)},\]
 which combined with \eqref{eq:identity_flux} gives  \begin{equation}\label{eq:stress_fluid_derivative}\partial_x\int_{-\infty}^{f(x)}
      \left(\frac{\Psi_x^2-\Psi_y^2}{2}+s\Psi_x\Psi_y\right)\,dy = -(f'-s)\left[\frac{(\varphi')^2}{2}
      -\frac{(G(f)\varphi+f'\varphi')^2}{2(1+f'^2)}\right]
      -G(f)\varphi\,\varphi'.
  \end{equation}
  For the surface term, we have 
  \[\partial_x \left(\sqrt{1+s^2}
  -\frac{1+s f'}{\sqrt{1+f'^2}}\right) = (f' - s)\frac{f''}{(1 + (f')^2)^{3/2}}\]
  and for the correction terms, the product rule gives
  \[\partial_x\left(-\frac{2(f-sx)\varphi}{3(1+\epsilon^2)}\right)
    =-\frac{2}{3(1+\epsilon^2)}\left[(f'-s)\varphi+(f-sx)\varphi'\right]\]
  and
  \[\partial_x\left(-\frac{\epsilon\left(2(f-sx)^2+\varphi^2\right)}{6(1+\epsilon^2)}\right)
    =-\frac{\epsilon}{3(1+\epsilon^2)}\left[2(f-sx)(f'-s)+\varphi\varphi'\right]\]
  Summing these derivatives, we obtain
  \begin{align*}
    \mc{R}_{s,\epsilon}'
    &=(f'-s)\left[\frac{f''}{(1+f'^2)^{3/2}}-\frac{(\varphi')^2}{2}
    +\frac{(G(f)\varphi+f'\varphi')^2}{2(1+f'^2)}\right]-G(f)\varphi\,\varphi' \\
    &\quad-\frac{2}{3(1+\epsilon^2)}\left[(f'-s)\varphi+(f-sx)\varphi'\right] \\
    &\quad-\frac{\epsilon}{3(1+\epsilon^2)}\left[2(f-sx)(f'-s)+\varphi\varphi'\right].
  \end{align*}
  Using the formulas \eqref{eq:energy_variational_derivatives} for the variational derivatives of $\mc{H}_s$ and collecting the coefficients of $f'-s$ and $\varphi'$, this becomes
  \begin{equation}\label{eq:stress_variational_derivative}\mc{R}_{s,\epsilon}' =-(f'-s)\left[\frac{\delta\mc{H}_s}{\delta f}
      +\frac{2\varphi+2\epsilon(f-sx)}{3(1+\epsilon^2)}\right] -\varphi'\left[\frac{\delta\mc{H}_s}{\delta\varphi}
      +\frac{2(f-sx)+\epsilon\varphi}{3(1+\epsilon^2)}\right].\end{equation}
  Finally, solving \eqref{eq:profile_reg} in the form \eqref{eq:energy_profile_system} for the variational derivatives and using $f-xf'=(f-sx)-x(f'-s)$ gives
  \[\frac{\delta\mc{H}_s}{\delta f}
    +\frac{2\varphi+2\epsilon(f-sx)}{3(1+\epsilon^2)}
    =\frac{3\varphi-2x\varphi'+2\epsilon x(f'-s)}{3(1+\epsilon^2)},\]
  and
  \[\frac{\delta\mc{H}_s}{\delta\varphi}
    +\frac{2(f-sx)+\epsilon\varphi}{3(1+\epsilon^2)}=\frac{2x\left[(f'-s)+\epsilon\varphi'\right]}{3(1+\epsilon^2)}.\]
  Substituting these expressions and cancelling the mixed terms, we obtain
  \begin{align*}
    \mc{R}_{s,\epsilon}'
    &=-\frac{(f'-s)\left[3\varphi-2x\varphi'+2\epsilon x(f'-s)\right]}{3(1+\epsilon^2)}
    -\frac{2x\varphi'\left[(f'-s)+\epsilon\varphi'\right]}{3(1+\epsilon^2)} \\
    &=-\frac{(f'-s)\varphi}{1+\epsilon^2}
    -\frac{2\epsilon x}{3(1+\epsilon^2)}\left[(f'-s)^2+(\varphi')^2\right],
  \end{align*}
  as wanted.
\end{proof}

\begin{prop}[Stress estimate]\label{prop:stress_bounds}
    Under the hypotheses of \Cref{prop:energy_balance}, the  stress satisfies the estimate
    \begin{equation}\label{eq:stress_bound}
    |\mc{R}_{\pm s,\epsilon}(x)| \lesssim s^2, \quad \pm x \geq 2.
    \end{equation}
\end{prop}

\begin{proof}
  First let $x \geq 2$. We claim that
  \begin{equation}\label{eq:stress_integral_bound}\int_{0}^{\infty}|\mc{R}_{s,\epsilon}(x)| \, dx \lesssim s^2.\end{equation}
  To see this, write
  \[\mc{R}_{s,\epsilon} = \mc{R}_1 + \mc{R}_2 + \mc{R}_3,\]
  where $\mc{R}_1$ consists of the graph term, $\mc{R}_2$ consists of the fluid term, and $\mc{R}_3$ consists of the remaining correction terms. \par
  To bound $\mc{R}_1$, we write
  \[\mc{R}_1 = \sqrt{1+s^2}-\frac{1+sf'}{\sqrt{1+f'^2}} = (f'-s)^2 \int_0^1 \frac{\theta\,d\theta}{\left[1+\bigl(s+\theta(f'-s)\bigr)^2\right]^{3/2}}.\]
  It follows that
  \[0 \leq \mc{R}_1(x) \leq |(f - s|x|)'|^2.\]
  To bound $\mc{R}_2$, we use
  \[\left|\frac{\Psi_x^2-\Psi_y^2}{2}
  +s\,\Psi_x\Psi_y\right| \leq \frac{1 + |s|}{2}(\Psi_x^2 + \Psi_y^2)\]
  so that
  \[|\mc{R}_2(x)| \leq \frac{1 + |s|}{2}\int_{-\infty}^{f(x)}|\nabla \Psi|^2 \, dy .\]
  To bound $\mc{R}_3$, we use that $0 \leq \epsilon \leq 1$ and $2|ab| \leq a^2 + b^2$ to get
  \[\left|\mc{R}_3(x)\right| \lesssim (f - sx)^2 + \varphi(x)^2.\]
  By combining these estimates with the energy bounds from \Cref{prop:energy_bounds}, we have that
  \[\int_{0}^{\infty}|\mc{R}_{s,\epsilon}(x)| \, dx \lesssim s^2.\]
  To control $\mc{R}_s$ locally, we need to use the derivative of the stress. By differentiating $\mc{R}_{s,\epsilon}$, we obtain
  \[
    \mathcal R_{s,\epsilon}'(x)
    =
    -\frac{(f'(x)-s)\varphi(x)}{1+\epsilon^2}
    -\frac{2\epsilon x}{3(1+\epsilon^2)}
    \left[
    (f'(x)-s)^2+\varphi'(x)^2
    \right].
  \]
  By integrating this identity \eqref{eq:stress_identity}, we obtain 
  \begin{equation}\label{eq:stress_integrated_identity}
    \begin{aligned}
    \mathcal R_{s,\epsilon}(b)
    -\mathcal R_{s,\epsilon}(a)
    ={}&
    -\frac1{1+\epsilon^2}
    \int_a^b(f'-s)\varphi\,dx\\
    &-\frac{2\epsilon}{3(1+\epsilon^2)}
    \int_a^b x\left[(f'-s)^2+\varphi'^2\right]\,dx.
    \end{aligned}
  \end{equation}
  By applying the Cauchy-Schwarz inequality to the first term, and using the energy bounds from \Cref{prop:energy_bounds}, we have that
  \begin{equation}\label{eq:stress_increment_bound}\mc{R}_{s,\epsilon}(b) \leq \mc{R}_{s,\epsilon}(a) + Cs^2, \quad 0 < a < b.\end{equation}
  By setting $b = x$ and averaging over all $a \in [x-1,x]$ in \eqref{eq:stress_increment_bound}, we obtain  
  \[\mc{R}_{s,\epsilon}(x) \leq \int_{0}^{\infty}|\mc{R}_{s,\epsilon}(y)| \, dy + Cs^2 \lesssim s^2.\]
  By setting $a = x$ and averaging over all $b \in [x,x+1]$ in \eqref{eq:stress_increment_bound}, we similarly obtain
  \[\mc{R}_{s,\epsilon}(x) \gtrsim -s^2.\]
  Combining these two estimates gives the desired bound $|\mc{R}_{s,\epsilon}(x)| \lesssim s^2$ for $x \geq 2$, with no restriction on the sign of $s$. Since $f$ and $\varphi$ are even, their finite-energy harmonic extension satisfies $\Psi(-x,y)=\Psi(x,y)$. The definition of the  stress in \Cref{def:relative_stress} therefore gives
  \begin{equation}\label{eq:stress_reflection}\mc{R}_{-s,\epsilon}(-x)=\mc{R}_{s,\epsilon}(x), \quad x\geq2,\end{equation}
  which yields $|\mc{R}_{-s,\epsilon}(x)|\lesssim s^2$ for $x\leq-2$. \par
\end{proof}

\section{A priori bootstrap estimates: uniform slope and velocity bounds}\label{sec:uniform_bounds}

In this section we  prove that if an admissible pair $(f,\varphi)$ is a solution of the regularized system with $\epsilon\geq0$, then its slope and velocity are controlled by the far-field slope $|s|$ in stronger norms, uniformly in the regularizing parameter $\epsilon$.
 This will allow us in Section \ref{sec:existence_regularized} to prove the existence of regularized solutions using degree theory. Specifically, it will be used in Proposition \ref{prop:existence_small_unif} to guarantee that the branch of regularized solutions does not escape a bounded set.

We introduce the following quantity, whose control in terms of $|s|$ is the goal of the section:
\begin{equation}\label{def:Mdef}
  M(f,\varphi) := \max\Big\{\|f'\|_{L^\infty}, \sup_{x\in\mathbb{R}}\frac{|\nabla\Psi(x,f(x))|} {\langle x\rangle}\Big\}
\end{equation}
While the first quantity is obviously the slope of the interface, the weighted velocity also has a natural interpretation that  provides intuition on the overall argument. If we think of a point with fixed self-similar coordinate moving along the wedge, the second quantity involved in $M(f,\varphi)$ represents, at a linear level, the slope that a ripple traveling at that speed on the interface would  have.

\begin{rem}\label{rem:ripple}
    In fact, a ripple like $\eta(X,t)=A\cos(\xi X-\omega(\xi)t)$ with $\omega(\xi)=\xi^{3/2}$ has phase speed $U(\xi)=\omega(\xi)/\xi=\sqrt{\xi}$ and it has slope magnitude of $|A|\xi$. The kinematic condition at linear level $\Phi_Y=\eta_t$ at $y=0$ gives $\Phi(X,Y,t)=A\omega(\xi) \frac{e^{\xi Y}}{\xi}\sin{(\xi X-\omega(\xi)t)}$, thus the fluid velocity at the interface satisfies $|\nabla\Phi|=|A|\omega$, and we find the standard result that the ratio between the fluid velocity and the phase speed equals the slope.
 Denote $X(t)=(T-t)^{2/3}x$
a point that moves with fixed self-similar coordinate. For the linear evolution, stationary phase gives  $\sqrt{\xi_0}=\frac23X(t)/(T-t)=\frac23(T-t)^{-1/3}x$, and  we find
\begin{equation*}
\begin{aligned}
|\dot{X}(t)|=\frac23(T-t)^{-1/3}x,\qquad U(\xi_0)=\sqrt{\xi_0}=\frac23(T-t)^{-1/3}x.    
\end{aligned}
\end{equation*} 
Since the  fluid velocity is $|\nabla\Phi(x,f(x))|=(T-t)^{-1/3}|\nabla\Psi(x,f(x))|$, we obtain $|\nabla\Phi|/|\dot{X}|=3/2|\nabla\Psi|/x$.
\end{rem}

In the proofs of both Lemma \ref{lem:velocity_below} and Proposition \ref{prop:uniform_slope_bounds}, it will be convenient to work with the complex velocity.
Given $\varphi\in H^{\frac12}$, we define the complex potential $F:\Omega_f\to \mathbb{C}$  by
\begin{equation}\label{eq:F}
F(z)=\Psi(x,y)+i\Theta(x,y), \quad z=x+iy,
\end{equation}
where the velocity potential $\Psi:\Omega_f\to\mathbb{R}$ is the unique harmonic extension of $\varphi$ with $\nabla\Psi\in L^2(\Omega_f)$ and the stream function $\Theta$ is the harmonic conjugate of $\Psi$ (see Appendix \ref{app:harmonic}). We further denote its trace by $\Theta_b(x)=\Theta(x,f(x))$, so that $F_b(x)=F(x+if(x))=\varphi(x)+i\Theta_b(x)$. 
 We notice that $\varphi'$ and $\Theta_b'$ correspond to the tangential and normal velocities at the boundary,
\begin{equation}\label{eq:tan_nor_vel}
 \begin{aligned}
     \varphi'&=\Psi_x+f'\Psi_y,\qquad
     \Theta_b'=f'\Psi_x-\Psi_y=-G(f)\varphi.
 \end{aligned}
 \end{equation}
We also define the complex velocity
\begin{equation}\label{eq:complex_vel}
    W(z)=F'(z)=\Psi_x(x,y)-i\Psi_y(x,y),
\end{equation}
so that the previous relations are encoded into the derivative of $F_b$, 
\begin{equation}\label{eq:W_varphi_DN}
    F_b'(x)=W(x+if(x))\big(1+if'(x)\big)=\varphi'(x)-iG(f)\varphi(x).
\end{equation}
We notice that we use the same notation $W$ for its trace on the interface.

For the rest of this section, the constants may depend on the fixed $\delta$ in \eqref{def:admissible}, but are independent of $s$, $\epsilon$, and the solution $(f,\varphi)$.

\begin{lem}\label{lem:velocity_below}
    Let $0\leq \epsilon\leq 1$ and let $(f,\varphi)\in C^2_{\mathrm{loc}}\times C^1_{\mathrm{loc}}$ be an admissible pair \eqref{def:admissible} satisfying the equations  \eqref{eq:profile_reg}.
    There exists constants $C>0$ and $M_*>0$ such that, if
\begin{equation*}
M=M(f,\varphi)\leq M_*, 
\end{equation*}
where $M(f,\varphi)$ was defined in \eqref{def:Mdef}, then for every $x_0\geq4$ and every $d_0\geq0$, 
\begin{equation}\label{eq:far_column_vel}
\int_{-\infty}^{f(x_0)-d_0}|W(x_0+iy)|^2dy \leq C\Big(\frac{M^2}{1+x_0^2d_0}+s^2\Big),
\end{equation}
and
\begin{equation}\label{eq:far_potential_trace}
|\varphi(x_0)|\leq C\Big(|s|+\frac{M\log{x_0}}{x_0}\Big).
\end{equation}
\end{lem}
\begin{proof}
Let us first prove \eqref{eq:far_column_vel}. 
Denoting the points on the interface by $z(x)=x+if(x)$, we look for a bound on $|W(z(x_0)+iy)|$ with $y<0$.
If $d_0\geq c>0$ for some constant $c$, the estimate follows since $W$ is holomorphic and \begin{equation}\label{eq:W_L2_s}
    \|W\|_{L^2(\Omega_{f})}=\|\nabla\Psi\|_{L^2(\Omega_{f})}\leq C|s|,
\end{equation}
by the energy estimate \eqref{eq:energy_bound_01}.
In fact, $W$ holomorphic implies, by the mean value property applied to its real and imaginary parts, followed by Cauchy-Schwarz, that
\begin{equation*}
    \begin{aligned}
        |W(z(x_0)+iy)|^2&=\Big|\frac{1}{|B_{|y|/4}(z(x_0)+iy)|}\int_{B_{|y|/4}(z(x_0)+iy)}W(z)dz\Big|^2\\
        &\leq \frac{1}{\pi (y/4)^2}\int_{B_{|y|/4}(z(x_0)+iy)}|W(z)|^2dz\leq \frac{16}{\pi y^2}\|W\|_{L^2(\Omega_f)}^2\leq C\frac{s^2}{y^2}.
    \end{aligned}
\end{equation*}
 Therefore, when $d_0\geq c$, 
\begin{equation}\label{eq:far_col_aux}
    \begin{aligned}
        \int_{-\infty}^{-d_0}|W(z(x_0)+iy)|^2dy&\leq  C\frac{s^2}{\pi d_0}\leq Cs^2.
    \end{aligned}
\end{equation}
Next, we focus on the case when $d_0$ is arbitrarily small.
If $d_0x_0^2\leq 1$, we first notice that, for $x$ near $x_0$,  $|W(z(x))|\leq C M\langle x\rangle \leq CMx_0$ by the definition of $M$ \eqref{def:Mdef}. Moreover,  since $W$ is continuous up to the boundary by Lemma \ref{lem:velocity_continuity}, then Lemma \ref{lem:local_bulk} and \eqref{eq:W_L2_s} give 
\begin{equation}\label{eq:W_inf_bound}
    \|W\|_{L^\infty(\Omega_f\cap B_1)}\leq C Mx_0+C\|W\|_{L^2(\Omega_f)}\leq C Mx_0,
\end{equation}
and therefore we can split 
\begin{equation}\label{eq:W_near_easy}
    \begin{aligned}
    \int_{-\infty}^{-d_0}|W(z(x_0)+iy)|^2dy&\leq \int_{-\infty}^{-1/x_0^2}|W(z(x_0)+iy)|^2dy+ \int_{-1/x_0^2}^{-d_0}|W(z(x_0)+iy)|^2dy\\
    &\leq \int_{-\infty}^{-1/x_0^2}|W(z(x_0)+iy)|^2dy + CM^2.
    \end{aligned}
\end{equation}
Hence it remains to deal with the case in which $d_0x_0^2\geq1$, which in particular covers the first integral on the right-hand side above.
In this case, we look for a bound of $|W(z(x_0)+iy)|$ that decays with $x_0$. 
The idea will be to transfer the information from the interface to points below it. In particular, on the interface, the identity \eqref{eq:W_varphi_DN} connects $W(x+if(x))$ with $\varphi'(x)$, and the second equation in \eqref{eq:profile_reg} connects $x\varphi'$ with the curvature. 
In fact, let's multiply \eqref{eq:W_varphi_DN} by $z(x)-if(x_0)$, which equals $x_0$ at $x=x_0$, to obtain
\begin{equation}\label{eq:xW_boundary}
    \begin{aligned}
        \big(z(x)-if(x_0)\big)W\big(z(x)\big)z'(x)&=x\varphi'(x)+\big(f(x)-f(x_0)\big)G(f)\varphi(x)\\
        &\quad+i\Big(\varphi'(x)\big(f(x)-f(x_0)\big)-xG(f)\varphi(x)\Big).
    \end{aligned}
\end{equation}
We advance the main point in the proof: thanks to having multiplied by $z(x)-if(x_0)$, instead of simply $x$, we now notice that the term on the left could be written as
\begin{equation*}
    \frac{d}{dx}\mathscr{W}(z(x))=\mathscr{W}'(z(x))z'(x)=\mathcal{W}(z(x))z'(x),
\end{equation*}
with $\mathscr{W}'(z)=\mathcal{W}(z)=(z-if(x_0))W(z)$ holomorphic on $\Omega_f$. 
In this form, we can obtain bounds for $W$ on points below the interface from the boundary values.
We will first develop the right-hand side of \eqref{eq:xW_boundary} and absorb some additional terms before defining the corresponding holomorphic function.
We use the second equation in \eqref{eq:profile_reg} to write, 
\begin{equation}\label{eq:full_regularized01}
\begin{aligned}
   \frac{2}{3} x\varphi'(x)&=\frac{f''(x)}{(1+f'(x)^2)^{3/2}} + \frac{1}{3} \varphi(x) - \frac{1}{2} \varphi'(x)^2 + \frac{(G(f)\varphi(x) + f'(x)\varphi'(x))^2}{2(1+f'(x)^2)}   + \epsilon G(f)\varphi(x)\\
   &=\frac{f''(x)}{(1+f'(x)^2)^{3/2}} + \frac{1}{3} \varphi(x) - \frac12\text{Re}(z(x)'W(x)^2)   + \epsilon G(f)\varphi(x),
\end{aligned}
\end{equation}
where in the second line we use that, by \eqref{eq:tan_nor_vel}-\eqref{eq:complex_vel},
\begin{equation*}
    \text{Re}(z(x)'W(z(x))^2)=\varphi'(x)^2-\frac{\big(G(f)\varphi(x)+f'(x)\varphi'(x)\big)^2}{1+f'(x)^2}.
\end{equation*}
Then, we substitute $x\varphi'(x)$ in \eqref{eq:xW_boundary}:
\begin{equation*}
    \begin{aligned}
        \frac23\big(z(x)&-if(x_0)\big)W\big(z(x)\big)z'(x)=\frac{f''(x)}{(1+f'(x)^2)^{3/2}} + \frac{1}{3} \varphi(x) - \frac12\text{Re}(z(x)'W(x)^2)   + \epsilon G(f)\varphi(x)\\
        &\quad+\frac23\big(f(x)-f(x_0)\big)G(f)\varphi(x)+i\frac23\Big(\varphi'(x)\big(f(x)-f(x_0)\big)-xG(f)\varphi(x)\Big).
    \end{aligned}
\end{equation*}
We further notice from  \eqref{eq:W_varphi_DN} that
\begin{equation*}
    \text{Re}(z'(x)F_b(x))=\varphi(x)-f'(x)\Theta_b(x),\qquad \text{Re}(iz'(x)W(z(x)))=G(f)\varphi(x),
\end{equation*}
hence we can write
\begin{equation*}
    \begin{aligned}
        \frac23\big(z(x)&-if(x_0)\big)W\big(z(x)\big)z'(x)+\frac12\text{Re}\big(z(x)'W(x)^2\big)-\frac13\text{Re}(z'(x)F_b(x))-\epsilon\,\text{Re}(iz'(x)W(z(x)))\\
        &=\frac{f''(x)}{(1+f'(x)^2)^{3/2}} +\frac13f'(x)\Theta_b(x)+\frac23\big(f(x)-f(x_0)\big)G(f)\varphi(x)\\
        &\quad+i\frac23\Big(\varphi'(x)\big(f(x)-f(x_0)\big)-xG(f)\varphi(x)\Big),
    \end{aligned}
\end{equation*}
and, taking real part,
\begin{equation}\label{eq:Fboundary}
    \begin{aligned}
        \text{Re}[z'(x)\mathcal{F}(z(x))]=\frac{f''(x)}{(1+f'(x)^2)^{3/2}} +\frac13f'(x)\Theta_b(x)+\frac23\big(f(x)-f(x_0)\big)G(f)\varphi(x),
    \end{aligned}
\end{equation}
with\begin{equation*}
    \begin{aligned}
    \mathcal{F}(z(x))=\big(\frac23(z(x)-if(x_0))-i\epsilon\big)W(z(x))+\frac12W(z(x))^2-\frac13F_b(x).
    \end{aligned}
\end{equation*}
This is the value at the interface of a holomorphic function $\mathcal{F}(z)$ and  since $\Omega_f$ is simply connected, it has a primitive,
\begin{equation}\label{eq:F_in_domain}
    \begin{aligned}
    \mathscr{F}'(z)=\mathcal{F}(z)=\big(\frac23(z-if(x_0))-i\epsilon\big)W(z)+\frac12W(z)^2-\frac13F(z), \quad z\in\Omega_f.
    \end{aligned}
\end{equation}
We can use this equation to estimate $W(z(x_0)+iy)$ by estimating $\mathscr{F}'(z(x_0)+iy)$, which we will control through its boundary values. In fact,
\begin{equation*}
    \begin{aligned}
    \frac23x_0|W(z(x_0)+iy)|&\leq\big|\frac23(z(x_0)+iy-if(x_0))-i\epsilon\big||W(z(x_0)+iy)|\\
    &\leq  |\mathscr{F}'(z(x_0)+iy)|+\frac13|F(z(x_0)+iy)|+\frac12|W(z(x_0)+iy)|^2.
    \end{aligned}
\end{equation*}
Thanks to the estimate \eqref{eq:W_inf_bound}, the quadratic term in $W$ can be absorbed by taking $M_*$ small enough, 
\begin{equation}\label{eq:W_goal_F}
    \begin{aligned}
    x_0|W(z(x_0)+iy)|&\leq\big|(z(x_0)+iy-if(x_0))-i\epsilon\big||W(z(x_0)+iy)|\\
    &\leq  C\Big(|\mathscr{F}'(z(x_0)+iy)|+|F(z(x_0)+iy)|\Big).
    \end{aligned}
\end{equation}
It remains to bound the right-hand side terms.
Since the curvature is given by a derivative, we can write \eqref{eq:Fboundary} as
\begin{equation*}
    \frac{d}{dx}\text{Re}[\mathscr{F}(z(x))]=\frac{d}{dx}\Big(\frac{f'(x)}{(1+f'(x)^2)^{1/2}}\Big) +\frac13f'(x)\Theta_b(x)+\frac23\big(f(x)-f(x_0)\big)G(f)\varphi(x),
\end{equation*}
and integrating from $x_0$ to $x$, 
\begin{equation*}
    \text{Re}[\mathscr{F}(z(x))]=\frac{f'(x)}{(1+f'(x)^2)^{1/2}}-\frac{s}{(1+s^2)^{1/2}} +\int_{x_0}^x\bigg(\frac13f'(u)\Theta_b(u)+\frac23\big(f(u)-f(x_0)\big)G(f)\varphi(u)\bigg)du,
\end{equation*}
where the second term comes from the choice of the constant in $\mathscr{F}$. Recalling that $G(f)\varphi=-\Theta_b'$ \eqref{eq:tan_nor_vel}, the last two terms can be combined to obtain
\begin{equation}\label{eq:Fscr_boundary}
    \text{Re}[\mathscr{F}(z(x))]=\frac{f'(x)}{(1+f'(x)^2)^{1/2}}-\frac{s}{(1+s^2)^{1/2}} +\int_{x_0}^x f'(u)\Theta_b(u)du-\frac23(f(x)-f(x_0))\Theta_b(x).
\end{equation}
We proceed to estimate $|\mathscr{F}'|=|\nabla\text{Re}\mathscr{F}|$. We split
\begin{equation}\label{eq:F_scr_split}
\begin{aligned}
    \text{Re}[\mathscr{F}(z)]&=u_1(z)+u_2(z)+u_3(z)+R_{x_0}(z),
\end{aligned}
\end{equation}
where, for $j\in\{1,2,3\}$, $u_j(z)$ is the harmonic extension to $\Omega_f$ of $g_j(x)$, defined by 
\begin{equation*}
    \begin{aligned}
        g_1(x)&=\chi(x-x_0)\Big(\frac{f'(x)}{(1+f'(x)^2)^{1/2}}-\frac{s}{(1+s^2)^{1/2}}\Big),\\
        g_2(x)&=\chi(x-x_0)\int_{x_0}^x f'(u)\Theta_b(u)du,\\
        g_3(x)&=-\frac23(f(x)-f(x_0))\Theta_b(x)\chi(x-x_0),
    \end{aligned}
\end{equation*}
and $0\leq \chi\leq 1$ is a smooth cutoff supported on $|x|<2$ and equal to one on $|x|<1$. 
 The terms $g_j$ contain the  information of $\text{Re}[\mathscr{F}(z(x))]$ near $x_0$ and, by construction, $R_{x_0}(z(x))=0$ for $|x-x_0|\leq 1$.
By the mean value property for $\nabla u_j$, we have
\begin{equation*}
    \begin{aligned}
        \nabla u_j(z(x_0)+iy)&=\frac{1}{|B_{|y|/4}(z(x_0)+iy)|}\int_{B_{|y|/4}(z(x_0)+iy)}\nabla u_j(z)dz\\
        &=\frac{1}{|B_{|y|/4}(z(x_0)+iy)|}\int_{\partial B_{|y|/4}(z(x_0)+iy)} u_j(z) n(z)\, dz.
    \end{aligned}
\end{equation*}
The second expression together with the maximum principle gives
\begin{equation*}
    \begin{aligned}
        |\nabla u_j(z(x_0)+iy)|&\leq\frac{|\partial B_{|y|/4}(z(x_0)+iy)|}{|B_{|y|/4}(z(x_0)+iy)|}\|u_j\|_{L^\infty(B_{|y|/4}(z(x_0)+iy))}\leq \frac{C}{|y|}\|g_j\|_{L^\infty},
    \end{aligned}
\end{equation*}
while Cauchy-Schwarz on the first expression gives
\begin{equation*}
    \begin{aligned}
        |\nabla u_j(z(x_0)+iy)|&\leq \frac{1}{|B_{|y|/4}(z(x_0)+iy)|^{1/2}}\|\nabla u_j\|_{L^2(\Omega_f)}\leq\frac{C}{|y|}\|\nabla u_j\|_{L^2(\Omega_f)}.
    \end{aligned}
\end{equation*} 
By  $|s|\leq \|f'\|_{L^\infty}\leq M$, we have
\begin{equation*}
    \|g_1\|_{L^\infty}\leq CM,
\end{equation*}
while for $g_2$ we apply Cauchy-Schwarz and use the estimate on $\Theta_b$ in \eqref{eq:W_trace_varphi} and the bound \eqref{eq:energy_bound_01}, 
\begin{equation*}
\|g_2\|_{L^\infty}\leq CM\|\Theta_b\|_{L^2}\leq CM \|\varphi\|_{L^2}\leq CM|s|.
\end{equation*}
Therefore,
\begin{equation}\label{eq:nabla_u12_in}
    \begin{aligned}
        |\nabla u_1(z(x_0)+iy)|+|\nabla u_2(z(x_0)+iy)|&\leq \frac{C}{|y|}M,
    \end{aligned}
\end{equation}
 For $\nabla u_3$, Lemma  \ref{lem:H12_trace}  gives
\begin{equation*}
    \begin{aligned}
        |\nabla u_3(z(x_0)+iy)|&\leq \frac{C}{|y|}\|\nabla u_3\|_{L^2(\Omega_f)}\leq \frac{C}{|y|}\|g_3\|_{H^{1/2}},
    \end{aligned}
\end{equation*}
while by the reverse inequality in Lemma \ref{lem:H12_trace}, together with the Cauchy-Riemann equations and \eqref{eq:W_L2_s},
\begin{equation*}
    \begin{aligned}
        \|\Theta_b\|_{\dot{H}^{1/2}}\leq C\|\nabla \Theta\|_{L^2(\Omega_f)}=C\|W\|_{L^2(\Omega_f)}\leq C|s|,
    \end{aligned}
\end{equation*}
 thus
\begin{equation}\label{eq:nabla_u3_in}
    \begin{aligned}
        |\nabla u_3(z(x_0)+iy)|&\leq \frac{C}{|y|}M|s|,
    \end{aligned}
\end{equation}
  It remains to estimate $\nabla R_{x_0}(z)$ in \eqref{eq:F_scr_split}. Since by construction it is a harmonic function with $R_{x_0}(z(x))=0$ on $|x-x_0|<1$, we can apply Lemma \ref{lem:gradient_decay} with $-1/16<y<0$:  \begin{equation*}
     \begin{aligned}
         |\nabla R_{x_0}(z(x_0)+iy)|&\leq \frac{C}{|y|^{1/3}}\|R_{x_0}\|_{L^1(\Omega_f\cap B_{1}(z(x_0)))}\leq \frac{C}{|y|^{1/3}}\int_{|x-x_0|\leq 1}\int_{-2}^0|R_{x_0}(z(x)+i\lambda)|d\lambda dx.
     \end{aligned}
 \end{equation*}
We estimate the last integral by recovering $R_{x_0}(z)$ from \eqref{eq:F_scr_split},
\begin{equation}\label{eq:nabla_R_split}
     \begin{aligned}
         |\nabla R_{x_0}(z(x_0)+iy)|&\leq \frac{C}{|y|^{1/3}}\int_{|x-x_0|\leq 1}\int_{-2}^0
         \Big(|\text{Re}[\mathscr{F}(z(x)+i\lambda)]|+\sum_{j=1}^3|u_j(z(x)+i\lambda)|\Big)d\lambda dx\\
         &=\frac{C}{|y|^{1/3}}\big(I_1+I_2\big),
     \end{aligned}
 \end{equation}
and using \eqref{eq:F_in_domain} for the first term.
We write $\text{Re}[\mathscr{F}(z(x)+iy)]$ in terms of its derivative and use its value at $y=0$,
\begin{equation*}
     \begin{aligned}
         I_1&\leq \int_{|x-x_0|\leq 1}\int_{-2}^0\Big(
         |\text{Re}[\mathscr{F}(z(x))]|+\int_{y}^0|\partial_y \text{Re}[\mathscr{F}(z(x)+i\lambda)]|d\lambda \Big)dydx\\
         &\leq  \sum_{j=1}^3\|g_j\|_{L^1}+\int_{|x-x_0|\leq1}\int_{-2}^0|\nabla \text{Re}[\mathscr{F}(z(x)+iy)]| dydx,
     \end{aligned}
 \end{equation*}
so $\|g_j\|_{L^1}\leq C\|g_j\|_{L^2}\leq C|s|$ and  substituting $|\nabla \text{Re}[\mathscr{F}]|=|\mathscr{F}'|$ from \eqref{eq:F_in_domain}, we obtain
\begin{equation*}
     \begin{aligned}
         I_1&\leq C|s|+\int_{|x-x_0|\leq1}\int_{-2}^0\Big(x_0|W(z(x)+iy)|+|W(z(x)+iy)|^2+|F(z(x)+iy)|\Big) dydx\\
         &\leq C|s|+Cx_0\|W\|_{L^2(\Omega_f)}+C\|W\|_{L^2(\Omega_f)}^2+\int_{|x-x_0|\leq1}\int_{-2}^0|F(z(x)+iy)|dydx\\
         &\leq Cx_0|s|+\int_{|x-x_0|\leq1}\int_{-2}^0|F(z(x)+iy)|dydx,
     \end{aligned}
 \end{equation*}
where we used that 
$\|W\|_{L^2(\Omega_f)}\leq C|s|$ by \eqref{eq:W_L2_s}. Finally, we estimate $F(z(x)+iy)$ from its boundary values. We use its expression as a Cauchy integral followed by Cauchy-Schwarz to get
\begin{equation}\label{eq:F_in_bound}
     \begin{aligned}
         |F(z(x)+iy)|&\leq C\|F_b\|_{L^2}\Big(\int_{\mathbb{R}}\frac{dt}{(x-t)^2+y^2}\Big)^{1/2}\leq C\frac{|s|}{|y|^{1/2}},
     \end{aligned}
 \end{equation}
hence
\begin{equation*}
    \begin{aligned}
        I_1\leq Cx_0|s|.
    \end{aligned}
\end{equation*}
Since $u_j$ is the harmonic extension of $g_j$, we similarly find that
\begin{equation}\label{eq:u_in_bound}
     \begin{aligned}
         |u_j(z(x)+iy)|&\leq C\|g_j\|_{L^2}\Big(\int_{\mathbb{R}}\frac{dt}{(x-t)^2+y^2}\Big)^{1/2}\leq C\frac{|s|}{|y|^{1/2}},
     \end{aligned}
 \end{equation}
therefore, back in \eqref{eq:nabla_R_split}, we have
\begin{equation*}
    |\nabla R_{x_0}(z(x_0)+iy)|\leq C\frac{x_0|s|}{|y|^{1/3}}. 
\end{equation*}
Inserting this bound and those for $\nabla u_j$ \eqref{eq:nabla_u12_in}-\eqref{eq:nabla_u3_in} back in \eqref{eq:F_scr_split}, 
\begin{equation*}
    |\mathscr{F}'(z(x_0)+iy)|\leq C\Big(\frac{M}{|y|}+\frac{x_0|s|}{|y|^{1/3}}\Big),
\end{equation*}
and together with \eqref{eq:F_in_bound} back in \eqref{eq:W_goal_F}, 
\begin{equation}\label{eq:W_in_pointwise}
    |W(z(x_0)+iy)|\leq C\Big(\frac{M}{x_0|y|}+\frac{|s|}{|y|^{1/3}}+\frac{|s|}{x_0|y|^{1/2}}\Big)\leq C\Big(\frac{M}{x_0|y|}+\frac{|s|}{|y|^{1/3}}\Big),\qquad -\frac{1}{16}<y<0.
\end{equation}
We can thus conclude that 
\begin{equation*}
    \begin{aligned}
        \int_{-\infty}^{f(x_0)-d_0}|W(x_0+iy)|^2dy&=\int_{-\infty}^{-\frac{1}{16}}|W(z(x_0)+iy)|^2dy+\int_{-\frac{1}{16}}^{-d_0}|W(z(x_0)+iy)|^2dy\\
        &\leq Cs^2+\frac{CM^2}{x_0^2}\int_{-\frac{1}{16}}^{-d_0}\frac{dy}{|y|^2}+Cs^2\int_{-\frac{1}{16}}^{-d_0}\frac{dy}{|y|^{2/3}}\\
        &\leq C\Big(\frac{M^2}{x_0^2d_0}+s^2\Big),
    \end{aligned}
\end{equation*}
which combined with that of the case $x_0^2d_0\leq 1$ \eqref{eq:W_near_easy} gives desired estimate \eqref{eq:far_column_vel}.

It remains to prove the inequality \eqref{eq:far_potential_trace}.
Recalling from \eqref{eq:F}-\eqref{eq:complex_vel} that $\varphi(x)=\text{Re}[F_b(x)]=\text{Re}[F(z(x))]$ and $F'(z)=W(z)$, we have
\begin{equation*}
    |\varphi(x_0)|\leq |F(z(x_0))|\leq \Big|F\Big(z(x_0)-\frac{i}{16}\Big)\Big|+\int_{-\frac{1}{16}}^0 \big|W(z(x_0)+iy)\big|dy,
\end{equation*}
so using the previous bound on $F$ \eqref{eq:F_in_bound} and the bounds \eqref{eq:W_in_pointwise}, \eqref{eq:W_inf_bound} on $W$, we find
\begin{equation*}
    \begin{aligned}
    |\varphi(x_0)|&\leq  C|s|+C\int_{-\frac{1}{16}}^{-1/x_0^2}\Big(\frac{M}{x_0|y|}+\frac{|s|}{|y|^{1/3}}\Big)dy+CMx_0\int_{-1/x_0^2}^0dy\\
    &\leq C\Big(|s|+\frac{M\log{x_0}}{x_0}\Big).
    \end{aligned}
\end{equation*}

\end{proof}

\begin{prop}\label{prop:uniform_slope_bounds}
Let $0\leq \epsilon\leq 1$ and let $(f, \varphi)\in C^2_{loc}\times C^1_{loc}$ be an admissible pair \eqref{def:admissible} satisfying the equations \eqref{eq:profile_reg}. If $\delta>0$ is sufficiently small, there exists $C>0$ such that if $M(f,\varphi) \leq \delta$, then 
\begin{equation*}
  M(f, \varphi) \leq C |s|,
\end{equation*}
where $M(f,\varphi)$ was defined in \eqref{def:Mdef}.
\end{prop}

\begin{proof}
We proceed by contradiction. For $\delta = 1/n$, take $C=n$. By the contradiction assumption, for sufficiently large $n$, there is $s_n$, $\epsilon_n$, $(f_n,\varphi_n)$ such that
\begin{equation}\label{eq:Mn_assump}
    n|s_n|< M(f_n,\varphi_n)\leq \frac{1}{n}.
\end{equation}
By the definition of both quantities in $M_n=M(f_n,\varphi_n)$, there is always a point $x_n$ for which
\begin{equation}\label{eq:contradiction}
    \Big(1-\frac{1}{n}\Big)M_n\leq \max\Big\{|f_n'(x_n)|, \frac{|W_n(x_n+if_n(x_n))|} {\langle x_n\rangle}\Big\}\leq M_n,
\end{equation}
where we wrote $|W_n(x_n+if_n(x_n))|=|\nabla\Psi_n(x_n,f_n(x_n))|$. 
Since $(f_n,\varphi_n)$ are even, we can take $x_n\geq0$.
We will next show that for $n\to\infty$, no choice of $s_n$, $\epsilon_n$, $(f_n,\varphi_n)$ can satisfy the above inequality \eqref{eq:contradiction}.
More precisely, noticing from \eqref{eq:Mn_assump} that $M_n\leq 1/n\to0$ and
\begin{equation}\label{eq:snMn}
    \frac{|s_n|}{M_n}\leq \frac{1}{n}\to0,
\end{equation}
we will compute the ratios between the slope and velocity at $(x_n,f_n(x_n))$ with $M_n$ and show that they must converge to zero, contradicting \eqref{eq:contradiction}. 

We distinguish two cases: $x_n\in I$ with $I$ a bounded interval, and $x_n\to\infty$. In both cases, we will use that the energy bounds \eqref{eq:energy_bound_01} in Proposition \ref{prop:energy_bounds} give
\begin{equation}\label{eq:energy_normalized}
    \frac{\|f_n - s_n|x|\|_{H^1}}{M_n} + \frac{\|\varphi_n\|_{H^{1/2}}}{M_n} +\frac{\|W_n\|_{L^2(\Omega_{f_n})}}{M_n}\leq C \frac{|s_n|}{M_n}\to0,
\end{equation}
where we used that $\|W_n\|_{L^2(\Omega_{f_n})}=\|\nabla\Psi_n\|_{L^2(\Omega_{f_n})}\leq C \|\varphi_n\|_{H^{1/2}}$ by Lemma \ref{lem:H12_trace}. 

\textit{Case 1:} $x_n\in I$. 
First, we show that $\|f_n'\|_{L^\infty(I)}/M_n\to0$. Choose a fixed interval $I$ containing all $\left[ x_n - 8, x_n + 8 \right]$ for every $n$. The inequality \eqref{eq:energy_normalized} directly gives that $\|f_n'\|_{L^2(I)}/M_n\to0$. 
Moreover, it also gives $\|\varphi_n\|_{L^2(I)}/M_n\to0$, and since we are considering a bounded interval, the definition of $M_n$ and  the identity \eqref{eq:W_varphi_DN} show that
\begin{equation*}
    |W_n(x+if_n(x))|+|\varphi_n'(x)|+|G(f_n)\varphi_n(x)|\leq CM_n.
\end{equation*}
Therefore, $\|\varphi_n\|_{L^\infty(I)}/M_n\to0$. 
Using these estimates in the second equation of \eqref{eq:profile_reg},
\begin{equation}\label{eq:full_regularized02}
        \frac{f_n''}{(1+f_n'^2)^{3/2}}=- \frac{1}{3} \varphi_n + \frac{2}{3} x\varphi_n'+ \frac{1}{2} \varphi_n'^2 - \frac{(G(f_n)\varphi_n + f_n'\varphi_n')^2}{2(1+f_n'^2)}   - \epsilon_n G(f_n)\varphi_n,
\end{equation}
we find that $|f_n''|\leq CM_n$, hence concluding that $\|f_n'\|_{L^\infty(I)}/M_n\to0$.

It remains to show that $\|W_n\|_{L^\infty(I)}/M_n\to0$. Again, we have $\|W_n\|_{L^2(\Omega_{f_n})}/M_n\to0$ from \eqref{eq:energy_normalized}. We will then use Lemma \ref{lem:local_oblique} to obtain a modulus of continuity up to the boundary for $|W_n|/M_n$ uniform in $n$, which combined with the $L^2(\Omega_{f_n})$ convergence concludes the result. In preparation for it, we combine the two equations in \eqref{eq:profile_reg} to eliminate the curvature term, 
\begin{equation*}
  (1+\epsilon_n^2)G(f_n)\varphi_n - \frac{2}{3}\epsilon_n x \varphi_n' = - \frac{2}{3}(f_n-xf_n') - \frac{1}{3}\epsilon_n \varphi_n,
\end{equation*}
which, recalling \eqref{eq:W_varphi_DN}, we rewrite in complex form in terms of $W_n$:
\begin{equation}\label{eq:withoutK_complex}
  \text{Re}\big(A_n(x) W_n(x+if_n(x))\big)=D_n(x),
\end{equation}
with 
\begin{equation}\label{eq:AD_def}
A_n(x)=(1+if_n')\big(-\frac23\epsilon_nx+i(1+\epsilon_n^2)\big),\qquad D_n(x)=-\frac23(f_n-xf_n')-\frac13\epsilon_n\varphi_n.
\end{equation}
We denote $e^{i\theta_n(x)}=A_n(x)/|A_n(x)|$. For $x\in I$, we have that $|A_n(x)|\geq1$ and, by the previous estimates, 
\begin{equation*}
\begin{aligned}
\|D_n\|_{L^\infty(I)}+\|D_n'\|_{L^\infty(I)}&\leq C M_n,\\
\|\theta_n'\|_{L^\infty(I)}&\leq \big\|(\log{A_n})'\big\|_{L^\infty(I)}=\Big\|\frac{A_n'}{A_n}\Big\|_{L^\infty(I)}\leq CM_n+\frac23\epsilon_n\leq C.    
\end{aligned}    
\end{equation*}
Therefore, recalling that $M_n\leq \delta$,  both $\theta_n$ and $\frac{D_n}{|A_n|M_n}$ are Lipschitz on $I$ with uniform Lipschitz constants. Both components of $W_n$ are harmonic and, by Lemma \ref{lem:velocity_continuity}, continuous up to the boundary, thus Lemma \eqref{lem:local_bulk} centered at $(x_n,f_n(x_n))$ gives
\begin{equation*}
    \|W_n\|_{L^\infty(\Omega_{f_n}\cap B_4)}\leq C M_n+C\|W_n\|_{L^2(\Omega_{f_n}\cap B_{8})}\leq CM_n,
\end{equation*}
where we used that $|W_n(x+if_n(x))|\leq C M_n$ for $x\in I$. Since $\| f_n' \|_{L^{\infty}}\leq M_n\leq 1/n$ we can apply \ref{lem:local_oblique} for sufficiently large $n$ to obtain  $\|W_n\|_{C^{1/2}(\Omega_{f_n}\cap B_2)}\leq C M_n$, which provides a uniform modulus of continuity for $W_n/M_n$. 

In conclusion, we have obtained that, for any bounded interval $I$,
\begin{equation*}
    \sup_{x\in I}\max\Big\{\frac{|f_n'(x)|}{M_n}, \frac{1}{M_n}\frac{|W_n(x+if_n(x))|} {\langle x\rangle}\Big\}\to 0,
\end{equation*}
which contradicts \eqref{eq:contradiction}.

\textit{Case 2:} $x_n\to \infty$.
We study the two quantities involved in $M_n$ near the point $x_n$. First, as in Case 1, we have that for $|x-x_n|\leq1$, the bound for $\varphi$ in \eqref{eq:far_potential_trace} and
the definition of $M_n$, together with the identity \eqref{eq:W_varphi_DN}, show that
\begin{equation}\label{eq:WnCase2bound}
\begin{aligned}
    \frac{|\varphi_n(x)|}{M_n}\leq C\Big(\frac{|s_n|}{M_n}+\frac{\log{x_n}}{x_n}\Big)&\to 0,\\ |W_n(x+if_n(x))|+|\varphi_n'(x)|+|G(f_n)\varphi_n(x)|&\leq CM_nx_n.
\end{aligned}
\end{equation}
Using these bounds in the equation for $f_n''$ \eqref{eq:full_regularized02}, we obtain that, on $|x-x_n|\leq 1$,
\begin{equation*}
    |f_n''|\leq CM_nx_n^2.
\end{equation*}
We perform a change of variables $(x,y)\to(\tilde{x},\tilde{y})$ localizing around $x_n$ at scale $x_n^{-2}$:
\begin{equation}\label{eq:rescaling}
    x=x_n+\frac{\tilde{x}}{x_n^2},\qquad y=f_n(x_n)+\frac{\tilde{y}}{x_n^2}.
\end{equation}
We notice that for any fixed $N>0$, $\tilde{x}\in(-N,N)$ implies that $|x-x_n|\leq1$, for sufficiently large $n$.
In the new variables, the interface $y=f_n(x)$ becomes $\tilde{y}=\Gamma_n(\tilde{x})$, with
\begin{equation*}
    \Gamma_n(\tilde{x})=x_n^2\big(f_n(x(\tilde{x}))-f_n(x_n)\big),
\end{equation*}
and we observe that, since $\Gamma_n(0)=0$, the interface in the new variables becomes locally flat in the limit,
\begin{equation*}
 |\Gamma_n'(\tilde{x})|=|f_n'(x(\tilde{x}))|\leq M_n\to 0, \qquad |\Gamma_n(\tilde{x})|\leq  M_n|\tilde{x}|\to0.
\end{equation*}
We also introduce the rescaled version of the quantities involved in \eqref{eq:contradiction}:
\begin{equation*}
    \ell_n(\tilde{x})=\frac{f_n'(x(\tilde{x}))}{M_n},\qquad \widetilde{W}_n(\tilde{x}+i\tilde{y})=\frac{1}{M_n\langle x_n\rangle}W_n\Big(x_n+if_n(x_n)+\frac{\tilde{x}+i\tilde{y}}{x_n^2}\Big).
\end{equation*}
On a bounded interval, the previous estimate on $f_n''$ gives that
\begin{equation}\label{eq:ell_bounds}
    |\ell_n|\leq1,\qquad |\ell_n'|=\frac{|f_n''(x)|}{M_nx_n^2}\leq C.
\end{equation}
To conclude, we will compute the limits of $\ell_n$ and $\widetilde{W}_n$ (Step 1), then we will find some identities that these limits must verify, obtained by using the equation for the curvature and the stress estimate from the previous section in the limit $n\to\infty$ (Step 2). Finally, we will use those identities to reach a contradiction with \eqref{eq:contradiction} (Step 3). 

\textit{Step 1:}
First, we show the existence of the limits of $\ell_n$, $\widetilde{W}_n$. 
On each $[-N,N]$, the bounds on $\ell_n$ \eqref{eq:ell_bounds} give uniform boundedness and equicontinuity of $\ell_n$, thus Arzelà-Ascoli gives a subsequence along which  $\ell_n$ converges uniformly in $[-N,N]$. 
For the velocity, as in Case 1, we use Lemma \ref{lem:local_bulk} at $(x_n,f_n(x_n))$ with the bounds \eqref{eq:WnCase2bound} and \eqref{eq:energy_normalized}. Since for $|x_n - x| < 2$ we have $\|W_n \|_{L^{\infty}(\partial \Omega_{f_n} \cap B_2)} \leq C M_n x_n$, we obtain
\begin{equation*}
    \|W_n\|_{L^\infty(\Omega_{f_n}\cap B_1)}\leq C M_n x_n+C\|W_n\|_{L^2(\Omega_{f_n}\cap B_{2})}\leq CM_nx_n,
\end{equation*}
thus we have 
\begin{equation}\label{eq:tildeW_bulk_uniform}
\|\widetilde{W}_n\|_{L^\infty(\Omega_{\Gamma_n}\cap B_{x_n^2}(0))}\leq C    
\end{equation}
We  then use Lemma \ref{lem:local_oblique} in the rescaled coordinates. Keep $A_n, \theta_n$, and $D_n$ as in Case 1. Then, the equation \eqref{eq:withoutK_complex} reads
\begin{equation*}
  \text{Re}\big(e^{i\theta_n(x(\tilde{x}))} \widetilde{W}_n(\tilde{x}+i\Gamma_n(\tilde{x}))\big)=g_n(x(\tilde{x})),\qquad g_n(x)=\frac{D_n(x)}{M_n\langle x_n\rangle |A_n(x)|}.
\end{equation*}
Recalling their expressions \eqref{eq:AD_def}, $|A_n|\geq1$, and that $x'(\tilde{x})=x_n^{-2}$, we obtain on each bounded interval $\tilde{x}\in[-N,N]$ that
\begin{equation*}
    \begin{aligned}
        \Big|\frac{d}{d\tilde{x}}\theta_n(x(\tilde{x}))\Big|&\leq \Big|\frac{1}{A_n}\frac{d}{d\tilde{x}}A_n(x(\tilde{x}))\Big|\leq \frac{|f_n''(x)|+\frac23\epsilon_n}{x_n^2}\leq C\big(M_n+\frac{1}{x_n^2}\big)\leq C,\\
        |g_n(x)|&\leq \frac23\frac{|f_n(x)|}{\langle x_n\rangle M_n}+\frac23\frac{x|f_n'(x)|}{\langle x_n\rangle M_n}+\frac13\epsilon_n \frac{|\varphi_n(x)|}{\langle x_n\rangle M_n}\leq C,\\
        \Big|\frac{d}{d\tilde{x}}g_n(x(\tilde{x}))\Big|&\leq \Big|\frac23\frac{x}{\langle x_n\rangle}\frac{f_n''(x)}{M_nx_n^2}-\frac13\epsilon_n\frac{\varphi_n'(x)}{M_n\langle x_n\rangle x_n^2}\Big|+|g_n(x)|\Big|\text{Re}\Big(\frac{A_n'(x)}{x_n^2A_n(x)}\Big)\Big|\leq C.        
    \end{aligned}
\end{equation*}
Hence, Lemma \ref{lem:local_oblique} gives 
\begin{equation*}
  \|\widetilde{W}_n\|_{C^{1/2}(\overline{\Omega_{\Gamma_n}}\cap\overline{B_N(0)})}\leq C.
\end{equation*}
Flattening the interface by the change of variables $\tilde{y}\to\Gamma_n(\tilde{x})+\zeta$, we obtain the velocities in a common domain, 
\begin{equation}\label{eq:tildeW_pointwise}
    U_n(\tilde{x}, \zeta)=\widetilde{W}_n\big(\tilde{x}+i(\Gamma_n(\tilde{x})+\zeta)\big), \qquad \zeta\leq0.
\end{equation}
Therefore, on each $[-N,N]\times[-N,0]$, we have $\|U_n\|_{C^{1/2}([-N,N]\times[-N,0])}\leq C$. Arzelà-Ascoli along a subsequence gives uniform convergence of $U_n$ on $[-N,N]\times[-N,0]$. A diagonal argument then provides $\ell_n(\tilde{x})\to \ell(\tilde{x})$ and $U_n(\tilde{x},\zeta)\to U(\tilde{x},\zeta)$ uniformly on compact sets of $\mathbb{R}$ and $\mathbb{R}\times(-\infty,0]$, respectively. We have thus obtained
\begin{equation*}
    \widetilde{W}_n\big(\tilde{x}+i(\Gamma_n(\tilde{x})+\zeta)\big)=U_n(\tilde{x},\zeta)\to U(\tilde{x},\zeta)=\widetilde{W}(\tilde{x}+i\zeta),
\end{equation*}
uniformly on compact sets of $\mathbb{R}\times(-\infty,0]$. Moreover, since $\Gamma_n\to0$, this also implies that $\widetilde{W}_n(\tilde{x}+i\tilde{y})\to\widetilde{W}(\tilde{x}+i\tilde{y})$ for $\tilde{y}<0$, hence $\widetilde{W}$ is holomorphic in the lower half-plane, continuous up to the boundary, and uniformly bounded.

\textit{Step 2:}
Next, we will obtain some identities that these limits $\ell$, $\widetilde{W}$ must satisfy. First, we notice that \eqref{eq:contradiction}, written in the rescaled variables, and taking $n\to\infty$ yields
\begin{equation}\label{eq:contradiction_limit}
\max\{|\ell(0)|,|\widetilde{W}(0)|\}=1.
\end{equation}
\textit{Step 2.1:} Next, we use the stress estimate \eqref{eq:stress_bound} in the limit $n\to\infty$. Recalling the expression of the stress \eqref{eq:stress_definition}, we write
\begin{equation}\label{eq:stress_split}
    \begin{aligned}
    \frac{\mathcal{R}_{n}(x)}{M_n^2} & =\frac{1}{M_n^2}\Big(\sqrt{1+s_n^2}
    -\frac{1+s_n f_n'(x)}{\sqrt{1+f_n'(x)^2}}\Big)+\frac12\int_{-\infty}^{f_n(x)}
    \text{Re}\Big((1+is_n)\Big(\frac{W_n(x+iy)}{M_n}\Big)^2\Big)\,dy
    \\
    &\quad-\frac{2(f_n(x)-s_nx)\varphi_n(x)}
            {3M_n^2(1+\epsilon_n^2)} - \frac{\epsilon_n\left(
    2(f_n(x)-s_nx)^2+\varphi_n(x)^2
    \right)}
    {6M_n^2(1+\epsilon_n^2)}\\
    &=\widetilde{\mathcal{R}}_{n,1}(x)+\widetilde{\mathcal{R}}_{n,2}(x)+\widetilde{\mathcal{R}}_{n,3}(x).
  \end{aligned}
\end{equation}
On the first term, we Taylor expand and use that $|s_n|/M_n\to 0$ \eqref{eq:snMn} and $|f_n'/M_n|=|\ell_n|\to|\ell|\leq 1$ to obtain 
\begin{equation*}
    \begin{aligned}
    \widetilde{\mathcal{R}}_{n,1}(x)&=\frac12\Big(\frac{f_n'}{M_n}\Big)^2+O\Big(\frac{|s_n|}{M_n}+M_n^2\Big)\to \frac12\ell^2.
  \end{aligned}
\end{equation*}
The last term vanishes in the limit by the estimate on $\varphi_n$ \eqref{eq:WnCase2bound} and the energy estimate \eqref{eq:energy_normalized}, together with Sobolev embedding, 
\begin{equation*}
    \begin{aligned}
    |\widetilde{\mathcal{R}}_{n,3}(x)|&\leq C\Big(\frac{|f_n(x)-s_nx|}{M_n}+\frac{|\varphi_n(x)|}{M_n}\Big)^2\to 0.
  \end{aligned}
\end{equation*}
Finally, for the middle term, we perform the change of variables $y=f_n(x_n)+\tilde{y}/x_n^2$ so that
\begin{equation*}
    \begin{aligned}
        \widetilde{\mathcal{R}}_{n,2}(x(\tilde{x}))&=\frac{1}{2}\int_{-\infty}^{f_n(x(\tilde{x}))}
    \text{Re}\Big[(1+is_n)\Big(\frac{W_n(x(\tilde{x})+iy)}{M_n}\Big)^2\Big]\,dy\\
    &=\frac{1}{2}\int_{-\infty}^{\Gamma_n(\tilde{x})}
    \text{Re}\bigg[(1+is_n)\bigg(\frac{W_n\big(x_n+if_n(x_n)+(\tilde{x}+i\tilde{y})/x_n^{2}\big)}{M_n\langle x_n \rangle}\bigg)^2\bigg]\frac{\langle x_n \rangle^2}{x_n^2}d\tilde{y}\\
    &=\frac{\langle x_n \rangle^2}{2x_n^2}\int_{-\infty}^{\Gamma_n(\tilde{x})}
    \text{Re}\Big[(1+is_n)\big(\widetilde{W}_n(\tilde{x}+i\tilde{y})\big)^2\Big]d\tilde{y}\to \frac12\int_{-\infty}^{0}
    \text{Re}\big[\big(\widetilde{W}(\tilde{x}+i\tilde{y})\big)^2\big]d\tilde{y}.
    \end{aligned}
\end{equation*}
The limit in the last step above is justified by splitting the integral in two pieces: $(-\infty,-L)$ and $(-L,\Gamma_n(\tilde{x}))$. Since $\Gamma_n(\tilde{x})\to0$, the convergence of the integral over $[-L,\Gamma_n(\tilde{x})]$ follows by the pointwise convergence \eqref{eq:tildeW_pointwise} and the uniform bound \eqref{eq:tildeW_bulk_uniform}. For the first part, we undo the change of variables and  use the estimate \eqref{eq:far_column_vel}
of Lemma \ref{lem:velocity_below} with $x=x(\tilde{x})$ and $d_0=(L+\Gamma_n(\tilde{x}))/x_n^2$:
\begin{equation}\label{eq:useofLemKinCol}
    \begin{aligned}
        \int_{-\infty}^{-L}|\widetilde{W}_n(\tilde{x}+i\tilde{y})|^2d\tilde{y}=\frac{x_n^2}{M_n^2\langle x_n\rangle^2}\int_{-\infty}^{f_n(x)-d_0}|W_n(x+iy)|^2dy\leq C\Big(\frac{x_n^2}{x^2(L+\Gamma_n(\tilde{x}))}+\frac{s_n^2}{M_n^2}\Big).
    \end{aligned}
\end{equation}
Then, Fatou's lemma gives 
\begin{equation}\label{eq:kin_col}
\int_{-\infty}^{-L}|\widetilde{W}(\tilde{x}+i\tilde{y})|^2d\tilde{y}\leq C/L,    
\end{equation}
thus letting $L\to\infty$ gives the convergence.
Going back to \eqref{eq:stress_split}, we obtain that 
\begin{equation*}
    \frac{\mathcal{R}_n(x(\tilde{x}))}{M_n^2}\to \frac12\ell(\tilde{x})^2+\frac{1}{2}\int_{-\infty}^0 \text{Re}\big[\big(\widetilde{W}(\tilde{x}+i\tilde{y})\big)^2\big]d\tilde{y},
\end{equation*}
and, by the stress estimate \eqref{eq:stress_bound}, $|\mathcal{R}_n(x)|/M_n^2\leq C \big(s_n/M_n\big)^2\to 0$, hence
\begin{equation}\label{eq:l_squared}
\ell(\tilde{x})^2=-\int_{-\infty}^0 \text{Re}\big[\big(\widetilde{W}(\tilde{x}+i\tilde{y})\big)^2\big]d\tilde{y}.
\end{equation}
\textit{Step 2.2:} We get a second identity relating the limits $\ell$ and $\widetilde{W}$ from the  equation  \eqref{eq:full_regularized02}. Dividing the equation by $M_nx_n^2$,
\begin{equation}
        \frac{1}{(1+f_n'^2)^{3/2}}\frac{f_n''}{M_nx_n^2}=- \frac{1}{3} \frac{\varphi_n}{M_nx_n^2} + \frac{2}{3} \frac{x\varphi_n'}{M_nx_n^2}+ \frac{1}{2} \frac{\varphi_n'^2}{M_nx_n^2} - \frac{(G(f_n)\varphi_n + f_n'\varphi_n')^2}{2M_nx_n^2(1+f_n'^2)}   - \epsilon_n \frac{G(f_n)\varphi_n}{M_nx_n^2},
\end{equation}
the estimates \eqref{eq:WnCase2bound} and the identity $\ell_n'(\tilde{x})=f_n''(x)/(M_nx_n^2)$ give that
\begin{equation}
        \frac{\ell_n'(\tilde{x})}{(1+(f_n'(x))^2)^{3/2}}= \frac{2}{3} \frac{x\langle x_n \rangle}{x_n^2}\frac{\varphi_n'(x)}{M_n\langle x_n\rangle}+ O\Big(M_n+\frac{\epsilon_n}{x_n}+\frac{1}{x_n^2}\Big).
\end{equation}
Then, recalling \eqref{eq:tan_nor_vel}, we have
\begin{equation*}
    \frac{\varphi_n'(x)}{M_n\langle x_n \rangle}=\frac{\partial_x \Psi_n(x,f_n(x))}{M_n\langle x_n \rangle}+f_n'(x)\frac{\partial_y \Psi_n(x,f_n(x))}{M_n\langle x_n \rangle}\to \text{Re}[\widetilde{W}(\tilde{x})],
\end{equation*}
where we use that the second term is $O(M_n)$. Therefore, it holds that
\begin{equation}\label{eq:l_derivative}
    \ell'(\tilde{x})=\frac23 \text{Re}[\widetilde{W}(\tilde{x})].
\end{equation}
\textit{Step 3:} The identity \eqref{eq:l_derivative} relates variations of the slope with the horizontal velocity of the fluid at the interface. The identity \eqref{eq:l_squared} relates the stress due to surface tension (more precisely, the difference between the actual stress and that of the wedge reference) with the difference between the vertical and horizontal kinetic energies (since $\text{Re}(\widetilde{W}^2)=(\text{Re}\,\widetilde{W})^2-(\text{Im}\,\widetilde{W})^2$) of the full column below the interface. Now, we can extend these identities to points below the interface, where we already know that this last quantity decays by \eqref{eq:kin_col}. In fact, denoting the Poisson kernel,
\begin{equation*}
    P(\tilde{x},\zeta)=\frac{1}{\pi}\frac{-\zeta}{\tilde{x}^2+\zeta^2},\qquad \zeta<0,
\end{equation*}
any bounded function $g(\tilde{x},\zeta)$ harmonic on the lower half-plane and continuous up to the boundary is expressed in terms of its boundary values by $g(\tilde{x}+i\zeta)=(P(\cdot,\zeta)* g|_{\mathbb{R}})(\tilde{x})$.
Then, since $\widetilde{W}$ is holomorphic, we can extend \eqref{eq:l_squared},
\begin{equation*}
(P(\cdot,\zeta)*\ell^2)(\tilde{x})=-\int_{-\infty}^0 \text{Re}\big[\big(\widetilde{W}(\tilde{x}+i(\tilde{y}+\zeta))\big)^2\big]d\tilde{y},
\end{equation*}
and by \eqref{eq:kin_col},
\begin{equation}\label{eq:decay_01}
|(P(\cdot,\zeta)*\ell^2)(\tilde{x})|=\Big|\int_{-\infty}^\zeta \text{Re}\big[\big(\widetilde{W}(\tilde{x}+i\tilde{y})\big)^2\big]d\tilde{y}\Big|\leq\frac{C}{|\zeta|}.
\end{equation}
Next, we show that this decay implies further decay for the horizontal velocity. 
We extend \eqref{eq:l_derivative} and use the gain of decay given by derivatives of the Poisson kernel,
\begin{equation*}
\begin{aligned}
    \frac23 \big|\text{Re}[\widetilde{W}(\tilde{x}+i\zeta)]\big|&=\big|(P(\cdot,\zeta)*\ell')(\tilde{x})\big|=\Big|\frac{\partial}{\partial \tilde{x}}(P(\cdot,\zeta)*\ell)(\tilde{x})\Big|\leq \big|(\nabla P(\cdot,\zeta)*\ell)(\tilde{x})\big|\\
    &\leq \big|(|\nabla P(\cdot,\zeta)|*|\ell|)(\tilde{x})\big|\leq \frac{1}{|\zeta|}\big|( P(\cdot,\zeta)*|\ell|)(\tilde{x})\big|,
\end{aligned}
\end{equation*}
thus by Cauchy-Schwarz and \eqref{eq:decay_01},
\begin{equation*}
\begin{aligned}
    \frac23 \big|\text{Re}[\widetilde{W}(\tilde{x}+i\zeta)]\big|&\leq \frac{1}{|\zeta|}\big|( P(\cdot,\zeta)*\ell^2)(\tilde{x})\big|^{\frac12}\leq \frac{C}{|\zeta|^{3/2}}.
\end{aligned}
\end{equation*}
Since $\widetilde{W}$ is holomorphic, and by \eqref{eq:kin_col} it does not contain a constant term, this decay also holds for its imaginary part, 
\begin{equation*}
\begin{aligned}
    \big|\widetilde{W}(\tilde{x}+i\zeta)\big|&\leq  \frac{C}{|\zeta|^{3/2}}.
\end{aligned}
\end{equation*}
Inserting this improved decay back in \eqref{eq:decay_01}, we have
\begin{equation*}
|(P(\cdot,\zeta)*\ell^2)(\tilde{x})|\leq\int_{-\infty}^\zeta \big|\big(\widetilde{W}(\tilde{x}+i\tilde{y})\big)^2\big|d\tilde{y}\leq\frac{C}{\zeta^2}.
\end{equation*}
We then conclude that $\ell$ must be zero. In fact, for any $N>0$ and $|\zeta|\geq N$,
\begin{equation*}
    \begin{aligned}
        \int_{-N}^N\ell(\tilde{x})^2d\tilde{x} \leq 2|\zeta| \Big|\int_{-N}^N\frac{\zeta}{\tilde{x}^2+\zeta^2}\ell(\tilde{x})^2d\tilde{x}\Big|\leq 2\pi|\zeta| \big|(P(\cdot,\zeta)*\ell^2)(0)\big|\leq \frac{C}{|\zeta|}, 
    \end{aligned}
\end{equation*}
hence taking $\zeta\to-\infty$ we conclude that $\ell\equiv0 $, and therefore $\widetilde{W}\equiv 0$ as well, contradicting \eqref{eq:contradiction_limit}.   
\end{proof}

\section{Compactness of the solution set}\label{sec:compactness}
The goal of this section is to use the a priori bounds of \ref{sec:energy}, to prove that the set of solutions, that are bounded in the quantity we obtained uniform bounds in \ref{sec:uniform_bounds}, is compact in $\R \times X$, where the space $X$ was defined in \eqref{def:X_Y_spaces}. In this section, we fix $0 < \epsilon \leq 1$ and $s_0 > 0$.  

\begin{prop}[Compactness of the solution set]\label{prop:compactness_solution_set}
For every fixed $0 < \epsilon \leq 1$, the set of solutions
\begin{equation*}
 K_\epsilon := \left\{ (s,h,\varphi) \in \left[ -s_0, s_0 \right] \times X : \mathcal{F}_\epsilon(s,h,\varphi) = 0, \  M(sf_{\mathrm ref} + h, \varphi) \leq \delta \right\}
\end{equation*}
is compact in $\left[ -s_0, s_0 \right] \times X$.
\end{prop}

To obtain this endgoal we will first need to prove boundedness of solutions in $X$, which will allow us to obtain then converging subsequences in weaker topologies. Then we will use the equation together with the lemmas who allowed us to obtain boundedness in $X$, to obtain strong convergence in $X$. 

Notice how all the estimates we do in this section, unlike some in previous ones are going to be crucially $\epsilon$ dependent, since the bounds need to degenerate as $\epsilon \to 0$, because the solutions of the original problem don't belong to $X$.

\subsection{$\epsilon$ dependent estimates}

The first natural goal of this subsection is to obtain some explicit bounds on the solutions in $K_\epsilon$ in some lower regularity, in particular $H^2 \times H^1$ bounds.

Let's first write fully the regularized equations \eqref{eq:profile_reg} 
\begin{equation}\label{eq:full_regularized}
\begin{aligned}
  - \frac{2}{3}(f- xf') &= G(f)\varphi  - \epsilon (\kappa[f] - N[f, \varphi])\\
  - \frac{1}{3} \varphi + \frac{2}{3} x\varphi' &=  ( \kappa[f] - N[f, \varphi]) + \epsilon G(f)\varphi,
\end{aligned}
\end{equation}
where the curvature $\kappa$, and the velocity nonlinearity $N$ are defined as
\begin{align}\label{def:curvature_N}
  \kappa[f] := \frac{f''}{(1+f'^2)^{3/2}},  \quad  \quad N[f, \varphi] := \frac{1}{2} \varphi'^2 - \frac{(G(f)\varphi + f'\varphi')^2}{2(1+f'^2)}.
\end{align}

These equations imply the following equation (taking out the common term $\kappa-N$)
\begin{equation}\label{eq:simple}
  (1+\epsilon^2)G(f)\varphi - \frac{2}{3}\epsilon x \varphi' = - \frac{2}{3}(f-xf') - \frac{\epsilon}{3} \varphi,
\end{equation}
which we are going to use in several steps of the proof. 

Also the assumption $M(h, \varphi) \leq \delta$ implies
\begin{equation}\label{eq:rewritten_velocity_bounds}
  \| f' \|_{L^\infty} \leq \delta, \quad \text{ and } \quad \frac{|G(f)\varphi | + |\varphi'|}{\langle x \rangle } \leq \delta (1+\delta^2). 
\end{equation}

Lastly, before proceeding with the Lemmas , we notice that $(s, h, \varphi)$ are in the hypothesis of admissible regular solution of Section \ref{sec:energy} Proposition \ref{prop:energy_bounds}, hence we know
\begin{equation*}
  \|h\|_{H^1} + \| \varphi\|_{H^{1/2}} \leq C|s_0|, \quad \| G(f) \varphi \|_{L^2} + \|\kappa - N\|_{L^2} \leq C |s_0|/\sqrt{\epsilon},
\end{equation*}
where we have used that
\begin{equation*}
  \|h \|_{H^1} \leq \| f - s|x|\|_{H^1} + s\| f_{ref} - |x|\|_{H^1} \leq C|s|.
\end{equation*}

Since $s_0, \epsilon, \delta$ are fixed parameters, we are not going to track the dependence of the estimates on them, we just want uniform bounds once these parameters are fixed. We are therefore going to use $A \lesssim B$ as there exists a constant $C$ depending on $s_0, \epsilon, \delta$ such that $A \leq C B$.

We first start with the low regularity estimate. 

\begin{lemma}\label{lem:low_regularity}
  Let $(s,h,\varphi)$ be a solution in $K_\epsilon$. Then there exists $C_{1, \epsilon}>0$ such that
  \begin{equation*}
    \| h \|_{H^2} + \|\varphi\|_{H^1} \leq C_{1, \epsilon}
  \end{equation*}
\end{lemma}
\begin{proof}
First we just use the $L^2$ estimate for the Dirichlet to Neumann map $G(f)$ from Lemma \ref{lem:Gf_Lp} with $p=2$
\begin{equation*}
  \| \varphi' \|_{L^2} \leq \| G(f)\varphi \|_{L^2} \lesssim 1,
\end{equation*}
where we used Proposition \ref{prop:energy_bounds}, we then have $\|\varphi\|_{H^1} \lesssim 1$. Now we focus on obtaining an $L^2$ bound on $h''$.
Then, we can use equation \eqref{eq:simple} to isolate $\varphi'$, and for $x \geq 1$ obtain
\begin{equation*}
  \begin{aligned}
  |\varphi'(x)| 
  &\leq  \frac{1}{\epsilon} \left [\frac{|f|}{|x|} + |f'| + 2 \epsilon \frac{|\varphi|}{|x|}+ \frac{3}{2}(1+\epsilon^2) \frac{|G(f)\varphi|}{|x|}\right ] \\
  &\lesssim \|f' \|_{L^\infty} + \| \varphi \|_{H^1} + \left \| \frac{G(f)\varphi}{\langle x \rangle}\right \|_{L^\infty} \lesssim 1,
\end{aligned}
\end{equation*}
where we used Sobolev embedding $\| \varphi \|_{L^\infty} \lesssim \| \varphi \|_{H^1}$, and also the bounds implied by $M(f, \varphi) \leq \delta$ in \eqref{eq:rewritten_velocity_bounds}. Therefore, since we also have from \eqref{eq:rewritten_velocity_bounds} $\varphi'(x) \lesssim |x|$, we conclude that $\| \varphi'\|_{L^\infty} \lesssim 1$.
Together with the $H^1$ bound, this implies 
\begin{equation*}
  \|\varphi'\|_{L^4} \lesssim \| \varphi'\|_{L^2} + \| \varphi'\|_{L^\infty} \lesssim 1.  
\end{equation*}

To conclude our argument, we use $f'' = (1+f'^2)^{3/2} \kappa$ and
\begin{equation*}
  \begin{aligned}
  \| h'' \|_{L^2} \lesssim \|f''\|_{L^2}  &\lesssim \|(1+f'^2)^{3/2}N \|_{L^2} + \| (1+f'^2)^{3/2} (\kappa - N)\|_{L^2}\\
  &\leq (1+\delta^2)^{3/2} \left ( \| (G(f)\varphi)^2\|_{L^2} + \| \varphi'^2\|_{L^2} + \| \kappa - N\|_{L^2} \right ) \\
  &\lesssim \| \varphi'\|_{L^4}^2 + \| \varphi'\|_{L^4}^2 + \| \kappa - N\|_{L^2} \lesssim 1.
\end{aligned}
\end{equation*}
where we used the estimate $\| G(f)\varphi \|_{L^4} \lesssim \| \varphi'\|_{L^4}$ given by Lemma \ref{lem:Gf_Lp} with $p=4$.
We have then concluded that $\|h\|_{H^2} + \|\varphi \|_{H^1} \lesssim 1$. 
\end{proof}

To continue towards bounding $(h, \varphi)$ in $X$ we need to upgrade first either $h$ to $H^3$ or $\varphi$ to $H^2$. We will proceed with the latter, using again the simplified equation \eqref{eq:simple}, grouping the top order terms in $\varphi'$ to the LHS we have
\begin{equation}\label{eq:oblique}
  G(f)\varphi - \frac{2\epsilon}{3(1+\epsilon^2)} x \varphi' =  \frac{1}{1+\epsilon^2} \left (-\frac{2}{3}(f-xf') - \frac{\epsilon}{3} \varphi \right ).
\end{equation}
Notice that if we neglected the weights in $x$, it seems like the LHS has less regularity $(L^2)$, than the RHS ($H^1$). Indeed, if we look at the symbol of the LHS, we have
\begin{align*}
  \text{Sym} (G(f) - b x \partial_x) (\xi)= |\xi| - i b x  \xi \implies |\text{Sym} (G(f) - b x \partial_x) (\xi)| = \sqrt{(1+ b^2 x^2)} |\xi|.
\end{align*}
so it provide exactly the weight we need to divide for in the RHS to have things integrable. It is natural to define
\begin{equation*}
  \mathcal A_{b, f} \varphi := \frac{G(f)\varphi - b x \varphi'}{\sqrt{1+b^2x^2}},
\end{equation*}
and the next lemma just gives a rigorous proof of the estimate we expect from the symbol computation.

\begin{lem}\label{lem:general_oblique}
  Let $f$ be such that 
  \begin{equation*}
    \|f' \|_{L^\infty} \leq \delta \quad \|f''\|_{L^2} \leq c \quad f'' \in L^\infty
  \end{equation*}
  then 
  \begin{equation*}
    \|\varphi\|_{H^2} + \| G(f)\varphi\|_{H^1} \leq C_{b, c} \left ( \| \mathcal A_{b, f} \varphi\|_{H^1} + \| \varphi \|_{H^1} \right ).
  \end{equation*}
  The constant $C_{b,c}$ is independent of $\|f'' \|_{L^\infty}$ and higher $f$ norms.
\end{lem}
\begin{proof}
  Let's write $\alpha = 1/w_b$ and $\beta = xb/w_b$ (notice $\alpha^2 + \beta^2 = 1$), then the operator $A_{f,b} = \alpha G(f) - \beta \partial_x$. With this notation, we first do the estimate for the flat case $f=0$. Take $T := \alpha H - \beta$, where $H$ is the Hilbert transform. Our goal is then to bound $\Ltwo{u'} \lesssim \Ltwo{Tu'} + \Ltwo{u}$. We proceed using a $T^*T$ method, since 
  \begin{equation*}
    T^*u' = -H(\alpha u') - \beta u' 
  \end{equation*}
  then
  \begin{equation*}
  \begin{aligned}
    T^*Tu' 
    &= -H(\alpha(\alpha Hu' - \beta u'))-\beta(\alpha Hu'- \beta u') \\
    &= -[H, \alpha^2] Hu' + \alpha^2 I + [H, \alpha\beta]u' + \alpha \beta Hu' - \beta\alpha Hu' + \beta^2 u' \\
    &= u' - [H, \alpha^2] Hu' + [H, \alpha\beta]u'.
  \end{aligned}
  \end{equation*}
  Finally, we can then do 
  \begin{equation}\label{eq:flat}
    \begin{aligned}
      \Ltwo{u'} 
      &\leq \Ltwo{T^*T u'} + \Ltwo{[H, \alpha^2] (Hu)'} \Ltwo{[H, \alpha\beta]}u'\\
      &\lesssim \Ltwo{Tu'} + b\Ltwo{u},
    \end{aligned}
  \end{equation}
  where we used the classical Calder\'on  estimate $\Ltwo{[H, c]v'} \lesssim \Linf{c'}\Ltwo{v}$, and $\Linf{(\alpha^2)' + (\alpha \beta)'} \lesssim b$.

  Now we proceed to the nonflat case. Let $\Psi$ be the harmonic potential with boundary trace $\varphi$, then let's define the trace velocity $U(x) = (V, B)(x) = (\partial_x \Psi, \partial_y \Psi)(x, f(x))$. By Lemma \ref{lem:DtN_analytic}, we know that $U \in H^1$ since $f'' \in L^\infty$ and $\varphi \in H^2$. Since we have the identities 
  \begin{equation}\label{eq:varphi'_G(f)varphi_VB}
     \varphi' = V + f' B \quad G(f)\varphi = B-f'V,
  \end{equation}
  we are going to bound in $H^1$ the velocity trace $U$ and recover for $\varphi', G(f)\varphi$. We first notice the following identity
  \begin{equation*}
    G(f)V (x)= (\partial_y \partial_x \Psi - f' \partial_x^2 \Psi)(x,f(x)) = (\partial_x \partial_y \Psi + f' \partial_y^2 \Psi)(x,f(x)) = B'(x),
  \end{equation*}
  hence $\Ltwo{B'} \lesssim \Ltwo{V'}$ by Lemma \ref{lem:Gf_Lp} with $p=2$, so we only need to bound $V'$. Let's now write $A_{f,b} \varphi$ in terms of $V,B$ using the identity \eqref{eq:varphi'_G(f)varphi_VB}
  \begin{equation*}
    \mathcal{A}_{f,b}\varphi = \alpha (B-f'V) - \beta(V+f'B) = \alpha B - \beta V - f'(\alpha V + \beta B).
  \end{equation*}
  Differentiating we have
  \begin{equation*}
    \begin{aligned}
      (\mathcal{A}_{f,b} \varphi)' 
      &= \alpha B' - \beta V' - f'(\alpha V' - \beta B') \\
      &- f''(\alpha V + \beta B) + \alpha'(B-f'V) - \beta'(V+f'B),
    \end{aligned}
  \end{equation*}
  so we conclude
  \begin{equation*}
    \Ltwo{\alpha B' - \beta V'} \lesssim \Hone{\mathcal{A}_{f,b} \varphi} + \delta \Ltwo{U'} + \Ltwo{f''}^2/\delta \Ltwo{U}
  \end{equation*}
  where we used $\Ltwo{f''U} \leq c^2/\delta \Ltwo{U} + \delta \Ltwo{U'}$.
  Finally, using $G(f)V = B'$ and the estimate we obtained for the flat case \eqref{eq:flat}
  \begin{equation*}
    \begin{aligned}
    \Ltwo{V'} 
    &\lesssim \Ltwo{\alpha HV' - bV'} +\Ltwo{V} \\
    &\leq \Ltwo{\alpha B' - \beta V'} + \Ltwo{G(f)V - HV'} + \Ltwo{V} \\
    &\lesssim \Hone{\mathcal{A}_{f,b} \varphi} + \delta \Ltwo{V'} +  \Ltwo{U} ,
  \end{aligned}  
  \end{equation*}
  where we substituted the last bound, and we used the $L^2$ estimate in Lemma \ref{lem:DtN_analytic} with $f_1=f$, $f_2=0$ to bound $\Ltwo{G(f)V - HV'} \lesssim \delta\Ltwo{V'}$. Decrease $\delta$ if necessary, so that the term with coefficient $\delta$ can be absorbed. This smallness assumption on $\delta$ is uniform in $\epsilon$. Then, we obtain
  \begin{equation}\label{eq:U'_l2_bound}
    \Ltwo{U'} \lesssim \Hone{\mathcal{A}_{f,b} \varphi} + \Hone{\varphi},
  \end{equation}
  where we used $\Ltwo{U} \leq \Ltwo{\varphi'} + \Ltwo{G(f)\varphi} \lesssim \Ltwo{\varphi'}$ by \eqref{eq:varphi'_G(f)varphi_VB} and Lemma~\ref{lem:Gf_Lp} with $p=2$. Finally, using
  \begin{equation*}
    \varphi'' = V' + f' B' + f'' B\quad \text{ and } \quad (G(f)\varphi)' = B' - f' V' - f'' V,
  \end{equation*}
  both RHS are bounded in $L^2$ by the bound \eqref{eq:U'_l2_bound} which is what we wanted.
\end{proof}

Using this last lemma, we are going to arrive to an equation of the following type for $f''$, so we want to bound it in terms of the $L^2$ of the RHS, which we already have bounded.

\begin{lem}\label{lem:curvature_estimate}
  Let $F, u$ in $L^2$ satisfy (weakly) the following ODE 
  \begin{equation*}
    \epsilon u' + x a(x) u = F,
  \end{equation*}
  where $a$ is measurable and $0 < a_* < a <a_*$, where $a_*, a^*$ are finite positive numbers. 
  Then
  \begin{equation*}
    \|u' \|_{L^2} + \|xu\|_{L^2} \leq C_{\epsilon, a_*, a^*} ( \|F \|_{L^2} + \|u \|_{L^2}).
  \end{equation*}
\end{lem}
\begin{proof}
  From the equation itself we get that $u$ is in $H^1_{loc}$, we then define $\chi_R(x) = \chi(x/R)$ where $\chi$ is a cut off function equal to one in $|x|\leq 1$ and $0$ in $|x|\geq 2$. Then, we test the equation against $x \chi_R^2 u$ and get
  \begin{equation*}
    \begin{aligned}
    \Ltwo{ \chi_R x u}^2 
    &= \langle F, x\chi_R^2 u \rangle_{2}  +\frac{\epsilon}{2} \langle \chi_R^2 + 2\chi_R \chi_R' x/R, u^2  \rangle_{2}\\
    &\leq \Ltwo{F} \Ltwo{x \chi_R u} + C  \Ltwo{u}^2,
    \end{aligned}  
  \end{equation*}
  so we conclude $a_*\Ltwo{xu} \leq \Ltwo{F} + C\Ltwo{u}$. Then using the equation we obtain
  \begin{equation*}
    \epsilon \Ltwo{u'} \leq \Ltwo{F} + a^* \Ltwo{xu},
  \end{equation*}
  which finishes the proof of the lemma.
\end{proof}

\begin{lem}\label{lem:higher_regularity}
  Let $(s, h, \varphi)$ be a solution in $K_\epsilon$. Then there exists $C_\epsilon>0$ such that
  \begin{equation*}
    \|(h, \varphi)\|_X \leq C_\epsilon.
  \end{equation*}
\end{lem}
\begin{proof}
  We already have from Lemma \ref{lem:low_regularity} that $\Htwo{h} + \Hone{\varphi} \lesssim 1$, so then now we will proof $\Ltwo{h'''} + \Ltwo{\varphi''} \lesssim 1$, and deduce from the equations that $\Hone{xh'} + \Hone{x\varphi'} \lesssim 1$.

  We first want to get $\Ltwo{\varphi''}\lesssim 1$. We start by writing equation \eqref{eq:oblique} divided by $\sqrt{1+b^2x^2}$
  \begin{align*}
    \mathcal A_{f,d} \varphi = \frac{2}{3(1+\epsilon^2)} \left (\frac{x}{w_b}h'-\frac{1}{w_b}h + s \frac{xf_{ref}'-f_{ref}}{w_b}\right ) - \frac{\epsilon}{3(1+\epsilon^2) w_b} \varphi
  \end{align*}
  where $w_b := \sqrt{1+b^2x^2}$ and $b= \frac{2\epsilon}{3(1+\epsilon^2)}$. 

  We then want to apply Lemma \ref{lem:general_oblique}, we check that we are indeed in the hypothesis:
  \begin{equation*}
  \begin{aligned}
    \Linf{f'} \leq M(f, \varphi) &\leq \delta, \\
     \Ltwo{f''} \lesssim s \Ltwo{f_{ref}''} + \Ltwo{h''} &\lesssim 1,\\
    h'' \in H^1 \implies f''= sf_{ref}'' + h'' &\in L^\infty.
  \end{aligned}
  \end{equation*}
  We then deduce
  \begin{align}\label{eq:varphi_Gfvarphi_bounds_H2H1}
    \Htwo{\varphi} + \Hone{G(f)\varphi} \lesssim s\Ltwo{xf_{ref}'-f_{ref}} + \Linf{\frac{x}{w_b}} \Ltwo{h'} + \Ltwo{h} + \Hone{\varphi} \lesssim 1.
  \end{align}

  To obtain $h''' \in L^2$ we are going to differentiate the first equation in \eqref{eq:full_regularized}, write all the local $f$ dependent terms on the LHS, and substitute the definition of the curvature and its derivative, to obtain
  \begin{equation}\label{eq:curvature_eq}
    \epsilon \kappa ' + \frac{2}{3}(1+f'^2)^{3/2} x\kappa = (G(f)\varphi)' + \epsilon N'.
  \end{equation}
  We want to apply Lemma \ref{lem:curvature_estimate}, clearly $a = \frac23(1+f'^2)^{3/2}$ is bounded above and below, then also $\kappa \in L^2$ since $h \in H^3$, and then lastly using the definition of $N$ \eqref{def:curvature_N}
  \begin{equation}\label{eq:N_H1}
    \Hone{N}  \lesssim \left ( \Hone{\varphi'}^2 + \Hone{G(f)\varphi}^2 \right ) \left ( 1 + \Ltwo{f''}\right ) \lesssim 1,
  \end{equation}
  where in the last step we used the previous bounds for $\varphi, G(f)\varphi$ \eqref{eq:varphi_Gfvarphi_bounds_H2H1}, and the already lower regularity bound for $f''$ in $L^2$ from \ref{lem:low_regularity}. We can then indeed apply Lemma \ref{lem:curvature_estimate} to conclude
  \begin{equation}\label{eq:curvature'_L2}
    \Ltwo{\kappa'} \lesssim \Ltwo{(G(f)\varphi) ' } + \Ltwo{N'} + \Ltwo{f''} \lesssim 1,
  \end{equation} 
  where in the first inequality we used $|\kappa| \leq |f''|$. Then we expand $f'''$
  \begin{equation}\label{eq:f'''_bound}
  \begin{aligned}
    \Ltwo{f'''} 
    &=  \Ltwo{(1+f'^2)^{3/2}\kappa ' + 3f'f''(1+f'^2)^{1/2}\kappa} \\
    &\lesssim (1+\delta^2)^{3/2} \Ltwo{\kappa'} + \delta(1+\delta^2)^{1/2} \Ltwo{f''} \Hone{\kappa} \lesssim 1.
  \end{aligned}
  \end{equation}
  Clearly, as we did many times before $f'''$ is bounded in $L^2$ if and only if $h'''$ is. Now we only need to bound $xh', x\varphi'$ in $H^1$. We write the equations \eqref{eq:full_regularized}
  \begin{equation}\label{eq:xh'_xvarphi}
  \begin{aligned}
    xh' &= h + s(f_{ref}-xf_{ref}') + \frac{3}{2}\left ( G(f)\varphi  - \epsilon (\kappa[f] - N[f, \varphi]) \right )\\
     x\varphi' &= \frac{1}{2}\varphi + \frac{3}{2} \left (( \kappa[f] - N[f, \varphi]) + \epsilon G(f)\varphi \right ),
  \end{aligned}
  \end{equation}
  and we then indeed check that all the terms on the right-hand side are bounded in $H^1$ already: $N$ in \eqref{eq:N_H1}, $\kappa$ in \eqref{eq:curvature'_L2} and Lemma \ref{lem:low_regularity}, $G(f)\varphi$ in \eqref{eq:varphi_Gfvarphi_bounds_H2H1}. We can then conclude $\X{(h, \varphi)} \lesssim 1$.
\end{proof}

\subsection{Proof of Proposition \ref{prop:compactness_solution_set}}

\begin{proof}
  We are going to do the same steps that we did to prove boundedness in $X$ in the previous Lemma \ref{lem:higher_regularity}. First prove convergence in $H^2 \times H^1$, then  in $L^2$ for $\varphi''$, then also in $L^2$ for $f'''$ and finally $H^1$ for the weighted quantities $x\varphi', xh'$.

  Assume we have a sequence $(s_j, h_j, \varphi_j) \in K_\epsilon$, and $f_j := s_j f_{ref} + h_j$, then by Lemma \ref{lem:higher_regularity} we know $(h_j, \varphi_j)$ is bounded in $X$, which is reflexive, hence
  \begin{equation*}
    s_j \to s, \quad \quad  (h_j, \varphi_j) \rightharpoonup (h, \varphi) \text{ in } X,
  \end{equation*}
  where $(s_j, h_j, \varphi_j)$ is a subsequence of the original sequence, we are going to abuse notation and use the same notation for subsequences all the time.  
  Let's define $\kappa_j = \kappa[f_j], N_j = N[f_j, \varphi_j]$ and $f= sf_{ref} + h$, $\kappa = \kappa[f], N = N[f, \varphi]$.

  \textit{Step 1:} \textbf{Convergence in $H^2 \times H^1$}. Since $(h_j, \varphi_j)$ is uniformly bounded in $X$, by Rellich compactness theorem, we have strong convergence in $H^2\times H^1$ in every compact interval. Also
  \begin{equation*}
    \Ltwo {h_j' 1_{|x|\geq R}} \leq \Ltwo{xh_j'}/R \quad \text{ and }\quad \Ltwo{h_j'' 1_{|x|\geq R}} \leq \Ltwo{(xh_j')'-h_j'}/R,
  \end{equation*}
   since $xh_j', x\varphi_j'$ are uniformly bounded in $H^1$ we can take a diagonal subsequence such that $(h_j', \varphi_j') \to (h',\varphi')$ strongly in $H^1\times L^2$. The bounds at infinity are the same for $\varphi', \varphi''$ but the local strong convergence is weaker. To get to $H^2 \times H^1$ strong convergence we need to recover the low frequencies.

  Since $f_j'-f'=(s_j-s)f_{\mathrm{ref}}'+h_j'-h'$, the convergence $s_j\to s$ and the embedding $H^1(\R)\hookrightarrow L^\infty(\R)$ give $f_j'\to f'$ uniformly. Also, $f_j''-f''=(s_j-s)f_{\mathrm{ref}}''+h_j''-h''\to0$ in $L^2$ since $f_{\mathrm{ref}}''$ is compactly supported. The third derivatives are bounded in $L^2$, so
  \begin{equation}
    \Linf{f_j''-f''}^2 \leq 2\Ltwo{f_j''-f''} \Ltwo{f_j'''-f'''} \to 0,
  \end{equation}
  hence $f_j'\to f'$ in $W^{1,\infty}$. Lemma \ref{lem:DtN_analytic} now gives
  \begin{equation}\label{eq:Gf_j-Gf_conv_H1}
    \|G(f_j) - G(f)\|_{\mathcal{L}(H^2, H^1)} \to 0.
  \end{equation}
  By Lemma \ref{lem:Gf_Lp} with $p=2$, we then have
  \begin{equation*}
    \begin{aligned}
      \Ltwo{G(f_j)\varphi_j - G(f)\varphi } &\leq  \|G(f_j) - G(f)\|_{\mathcal{L}(H^2, H^1)} \Htwo{\varphi_j} + \Ltwo{G(f)(\varphi_j-\varphi)} \\
      &\lesssim  \|G(f_j) - G(f)\|_{\mathcal{L}(H^2, H^1)} + \Ltwo{\varphi_j'-\varphi'} \to 0.
    \end{aligned}
  \end{equation*}
  Also, recalling the definition of $\kappa, N$ \ref{def:curvature_N}, we have $\kappa_j \to \kappa \text{ in } L^2$, since $f_j' \to f'$ in $L^\infty$ and $f_j'' \to f''$ in $L^2$. Lastly $N_j \to N$ in $L^2$, because $\varphi_j', G(f_j)\varphi_j$ are uniformly bounded in $H^1$ hence in $L^\infty$, $f_j' \to f'$ in $L^\infty$ and $\varphi_j', G(f_j)\varphi_j$ converge in $L^2$.
  
  Finally, we write the equations \eqref{eq:full_regularized}
  \begin{align*}
    \tfrac23(h_j-xh_j') &=-G(f_j)\varphi_j+\epsilon(\kappa_j-N_j) -\frac23s_j(f_{ref}-xf_{ref}'),\\
    \tfrac13\varphi_j-\tfrac23x\varphi_j' &=-\kappa_j+N_j-\epsilon G(f_j)\varphi_j.
  \end{align*}
  Since we have proved that the RHS converges strongly in $L^2$, and the LHS converges weakly, we know that $(f, \varphi)$ is a weak solution to the equations. Then writing the LHS of the difference, using notation and the computation of \eqref{eq:Tlambda_low} we get
  \begin{align*}
    \Ltwo{h_j-h} \leq \Ltwo{T_{1}(h_j-h)} \to 0 \quad  \Ltwo{\varphi_j-\varphi} \lesssim \Ltwo{T_{1/2}(\varphi_j-\varphi)} \to 0,
  \end{align*}
  since the RHS converges strongly in $L^2$. We have concluded 
  \begin{equation}\label{eq:h_varphi_H2H1_conv}
    (h_j, \varphi_j) \to (h, \varphi) \quad \text{ in } \quad H^2 \times H^1. 
  \end{equation}

  \textit{Step 2:} \textbf{Convergence of $\varphi'_j$  and $G(f_j)\varphi_j$ in $H^1$.}
  If we subtract the reduced equation of $j$-th term and the limit \eqref{eq:simple} we get
  \begin{equation*}
    \begin{aligned}
    \mathcal{A}_{b,f_j} (\varphi_j-\varphi)&= 
    \frac{2x}{3w_b} \frac{h_j'-h'}{1+\epsilon^2} - \frac{\frac23(h_j-h) - \frac\epsilon3(\varphi_j-\varphi)}{(1+\epsilon^2)w_b} + (s_j-s)\frac{(xf_{ref}'-f_{ref})}{(1+\epsilon^2)w_b}\\
    &+ \frac{(G(f)-G(f_j))\varphi}{(1+\epsilon^2)w_b},
  \end{aligned}
  \end{equation*}
  with the first and second terms converging in $H^1$ because of Step 1 \eqref{eq:h_varphi_H2H1_conv}, the third one by convergence of $s_j$ and the fourth one because of \eqref{eq:Gf_j-Gf_conv_H1}. Then, applying Lemma \ref{lem:general_oblique}
  \begin{equation}\label{eq:varphi_j_varphi_H2_conv}
    \Htwo{\varphi_j-\varphi} + \Hone{G(f_j)(\varphi_j-\varphi)} \lesssim \Hone{A_{b, f_j}(\varphi_j-\varphi)} + \Hone{\varphi_j-\varphi} \to 0, 
  \end{equation}
  where we could apply the lemma since the bounds on $f$ are satisfied by all $f_j$.

  To get the convergence of $G(f_j)\varphi_j$ we write
  \begin{equation}\label{eq:Gj-G_H1_conv}
    \Hone{ G(f_j)\varphi_j - G(f)\varphi} \leq \Hone{G(f_j)(\varphi_j-\varphi)} + \Hone{(G(f)-G(f_j))\varphi} \to 0,
  \end{equation}
  because of the last inequality \eqref{eq:varphi_j_varphi_H2_conv} and the equation in Step 1 \eqref{eq:Gf_j-Gf_conv_H1}.

  \textit{Step 3:} \textbf{Convergence of $h_j$ in $H^3$.}
  We follow the same process as in the boundedness in $X$, we obtain the equation for the curvature difference $\kappa_j-\kappa$, defining $a_j := \frac23(1+f_j'^2)^{3/2}$ and analogously $a$, then
  \begin{equation*}
    \epsilon (\kappa_j-\kappa)' + x a_j (\kappa_j-\kappa) = (G(f_j)\varphi_j)'-(G(f)\varphi)' + \epsilon(N_j'-N') + x(a-a_j)\kappa.
  \end{equation*}
  The first term tends to zero in $L^2$ by \eqref{eq:Gj-G_H1_conv}, and the last does so since $a_j\to a$ uniformly and $x\kappa\in L^2$. By Steps 1 and 2, $f_j'\to f'$ in $W^{1,\infty}$, while $\varphi_j'\to\varphi'$ and $G(f_j)\varphi_j\to G(f)\varphi$ in $H^1$. Since $H^1(\R)$ is an algebra and multiplication by $W^{1,\infty}$ functions is bounded on $H^1$, \eqref{def:curvature_N} gives $N_j\to N$ in $H^1$. Lemma \ref{lem:curvature_estimate} therefore gives
  \begin{equation*}
    \Ltwo{(\kappa_j-\kappa)'} + \Ltwo{x(\kappa_j-\kappa)} \to 0.
  \end{equation*}
  Expanding as in \eqref{eq:f'''_bound} gives $f_j'''\to f'''$ in $L^2$. Since $s_j\to s$ and $f_{\mathrm{ref}}'''\in L^2$, we obtain $h_j'''\to h'''$ in $L^2$, completing the $H^3$ convergence.

  \textit{Step 4:} \textbf{Convergence of $xh_j', x\varphi_j'$ and $M(f,\varphi)$ bounds.}
  Writing the equations as in \eqref{eq:xh'_xvarphi} and taking its difference with the $j-th$ one, we obtain that we have already proved the RHS converges in $H^1$ so 
  \begin{equation*}
    xh_j' \to xh' \text{ in }  H^1 \quad \text{ and } \quad  x\varphi_j' \to x\varphi' \text{ in }  H^1.
  \end{equation*}
  By Steps 1 and 2, $f_j'$, $G(f_j)\varphi_j$ and $\varphi_j'$ converge uniformly on $\R$. The boundary velocity formulas \eqref{eq:velocity_components} therefore give
  \begin{equation*}
    M(f,\varphi) = \lim_{j \to \infty} M(f_j, \varphi_j) \leq \delta,
  \end{equation*}
  and clearly $|s| = |\lim_{j} s_j| \leq s_0$, so $(h, \varphi)$ converges in $X$ and belongs to $K_\epsilon$.
\end{proof}

\section{Solving the regularized equation}\label{sec:existence_regularized}

\subsection{Linear analysis}

The goal of this section is to prove, for each fixed $0 < \epsilon \leq 1$, the existence of solutions $(s, h(s), \varphi(s))$ to the regularized system \eqref{eq:profile_reg} for $|s| \leq s_0$, with $s_0$ independent of $\epsilon$. Applying the implicit function theorem would only give a local branch of solutions, with an interval of existence that may depend on $\epsilon$. We discard this approach and instead we apply a continuation argument using degree theory, which requires the linearized operator to be Fredholm of index zero at every solution for which we use the degree.

Fix a point $(s,h,\varphi)$ in $U$, the domain of $\mathcal F_\epsilon$, and put $f=sf_{\mathrm{ref}}+h$. For a variation $(\tilde{h},\tilde{\varphi})\in X$, we have
\begin{equation*}
  D_{(h,\varphi)} \mathcal F_\epsilon(s, h, \varphi)[\tilde{h}, \tilde{\varphi}]
  =
  \begin{pmatrix*}
    \displaystyle -\frac{2}{3} \left(\tilde{h} - x \p_x \tilde{h}\right) \\[8pt]
    \displaystyle -\frac{1}{3} \tilde{\varphi} + \frac{2}{3}x \p_x \tilde{\varphi}
  \end{pmatrix*}
  -(\mathbb J+\epsilon\mathbb I) D_{(f,\varphi)}^2\mc{H}_s[f,\varphi]
    \begin{pmatrix} \tilde{h} \\ \tilde{\varphi} \end{pmatrix}.
\end{equation*}
The Hessian of the Hamiltonian has the explicit form 
\begin{equation*}
  D_{(f,\varphi)}^2\mc{H}_s[f,\varphi]
  \begin{pmatrix} \tilde{h} \\ \tilde{\varphi} \end{pmatrix}
  =
  \begin{pmatrix}
    \displaystyle -\partial_x\left(\frac{(\tilde{h})'}{(1+f'^2)^{3/2}}\right)+\Psi_x(\tilde{\varphi})'-\Psi_yG(f)\tilde{\varphi} +\Psi_y[G(f),\tilde{h}]\Psi_y \\
    \displaystyle G(f)\tilde{\varphi}-[G(f),\tilde{h}]\Psi_y-\Psi_x(\tilde{h})'
  \end{pmatrix},
\end{equation*}
where $\Psi_x$ and $\Psi_y$ are the velocity components in \eqref{eq:velocity_components}, and $[A,B]=AB-BA$. Note that to express the derivative of the Dirichlet-to-Neumann operator we have used the chain rule, the fact that $\Psi$ is harmonic, and the identity $G(f)\Psi_y = -\partial_x \Psi_x$ to get  
\begin{align*}
  D_{(f, \varphi)} (G(f)\varphi)[\tilde{h}, \tilde{\varphi}] &= G(f)(\tilde{\varphi} - \tilde{h} \Psi_y) + \tilde{h}(\Psi_{yy} - f' \Psi_{xy})|_{y = f(x)} - (\tilde{h})'\Psi_x \\
    &= G(f)\tilde{\varphi} - G(f)(\tilde{h}\Psi_y) - \partial_x(\tilde{h} \Psi_x) \\
    &= G(f)\tilde{\varphi} - G(f)(\tilde{h}\Psi_y) + \tilde{h} G(f)\Psi_y - \Psi_x(\tilde{h})' \\
    &= G(f)\tilde{\varphi} - [G(f), \tilde{h}]\Psi_y - \Psi_x(\tilde{h})'.
\end{align*}

We now extract from the linearized operator the dilation terms, the flat Dirichlet-to-Neumann operator $\Lambda = |D| = H\partial_x$ and the term with the curvature coefficient at infinity $(1 + s^2)^{-3/2}$. For every $\epsilon>0$ and $s\in\R$, these terms correspond to the operator
\begin{equation}\label{eq:L}
 \mathcal L_{\epsilon,s}(h,\varphi)
 =\begin{pmatrix}
     \displaystyle -\frac23(h-x h')+\frac{\epsilon h''}{(1+s^2)^{3/2}}-\Lambda\varphi\\
     \displaystyle -\frac13\varphi+\frac23 x\varphi' - \frac{h''}{(1+s^2)^{3/2}}-\epsilon\Lambda\varphi
   \end{pmatrix}
\end{equation}
After subtracting this, the remaining terms in $D_{(h,\varphi)}\mathcal F_\epsilon(s,h,\varphi)[\tilde{h},\tilde{\varphi}]$ are the highest derivative terms with variable coefficients plus a compact perturbation. In total we have the decomposition
\begin{equation}\label{eq:derivative_decomposition}
 D_{(h,\varphi)}\mathcal F_\epsilon(s,h,\varphi) = \mathcal L_{\epsilon,s}+\begin{pmatrix}-\epsilon\\1\end{pmatrix}\mathcal R+\mathcal K,
\end{equation}
where $\mathcal R$ is the scalar operator with the remaining variable coefficients
\begin{equation}\label{eq:R}
  \mathcal{R}(\tilde{h} , \tilde{\varphi}) = \bigl((1+s^2)^{-3/2}-(1+f'^2)^{-3/2}\bigr)(\tilde{h})'' + \Psi_x(\tilde{\varphi})'-\Psi_y\Lambda\tilde{\varphi},
\end{equation}
and $\mathcal{K}$ contains the remaining compact terms
\begin{equation}\label{eq:K}
  \begin{aligned}
  \mathcal{K}(\tilde{h} , \tilde{\varphi})
  &={}\begin{pmatrix}-\epsilon\\1\end{pmatrix}\left[
    \frac{3f'f''}{(1+f'^2)^{5/2}}(\tilde{h})'
    +\Psi_y[G(f),\tilde{h}]\Psi_y
    -\Psi_y(G(f)-\Lambda)\tilde{\varphi}
  \right]\\
  &\quad-\begin{pmatrix}1\\\epsilon\end{pmatrix}\left[
    (G(f)-\Lambda)\tilde{\varphi}
    -[G(f),\tilde{h}]\Psi_y-\Psi_x(\tilde{h})'
  \right].
  \end{aligned}
\end{equation}
The next lemmas show first that the operator $\mathcal L_{\epsilon,s}$ is invertible, and then that the terms in $\mathcal{K}$ are compact. The analysis of the operator $\mathcal{R}$ will be done later and will give the final step to prove that the linearized operator is Fredholm of index zero.

\begin{lem}[Inverse of linear operator]\label{lem:inverse_linear}
  The operator $\mathcal L_{\epsilon,s}:X\to Y$ is invertible for all $\epsilon > 0$ and $s \in \mathbb{R}$. Its inverse has norm locally bounded for $(\epsilon, s) \in (0, +\infty) \times \mathbb{R}$. 
\end{lem}
\begin{proof}
  The definition of $X$ shows that $\mathcal L_{\epsilon,s}:X\to Y$ is bounded. We first prove an a priori estimate for smooth, compactly supported solutions of the linear equation $\mathcal L_{\epsilon,s}(h,\varphi)=(g_1,g_2)$. Existence and uniqueness of solutions will follow from a density argument.

  \textit{Step 1:} To obtain the a priori estimate, we first get a bound for the low frequencies of the solution. When doing this, we disregard the terms with $\epsilon$. The regularization terms are only needed to control the high frequencies. We pair $-\mathcal L_{\epsilon,s}(h,\varphi)=-(g_1,g_2)$ with $h/(1+s^2)^{3/2}$ and $\langle D\rangle^{-1}\varphi$ in $L^{2}$ to cancel the cross terms in the equations and get a coercive estimate. For the dilation terms, integration by parts gives
\begin{align*}
 \operatorname{Re}\langle\tfrac23(h-x h'),h\rangle_{L^2}
   &=\frac43\int_{\R}|\widehat h|^2\,d\xi
     +\frac13\int_{\R}\xi\partial_\xi|\widehat h|^2\,d\xi\\
   &=\left(\frac43-\frac13\right)\int_{\R}|\widehat h|^2\,d\xi
     =\|h\|_{L^2}^2\,,\\
 \operatorname{Re}\langle\tfrac13\varphi-\tfrac23x\varphi',
                         \langle D\rangle^{-1}\varphi\rangle_{L^2}
   &=\int_{\R}\langle \xi\rangle^{-1}|\widehat\varphi|^2\,d\xi
     +\frac13\int_{\R}\xi\langle \xi\rangle^{-1}
                         \partial_\xi|\widehat\varphi|^2\,d\xi\\
   &=\int_{\R}\left[\langle \xi\rangle^{-1}
          -\frac13\left(\langle \xi\rangle^{-1}
                       -\xi^2\langle \xi\rangle^{-3}\right)\right]
                       |\widehat\varphi|^2\,d\xi\\
   &=\int_{\R}\left(\frac23+\frac{\xi^2}{3\langle \xi\rangle^2}\right)
                \langle \xi\rangle^{-1}|\widehat\varphi|^2\,d\xi\,.
\end{align*}
Using Plancherel's identity for the rest of the terms in the equation, we obtain
\begin{equation}\label{eq:low_freq_pairing}
\begin{aligned}
 &\operatorname{Re}\left\langle-\mathcal L_{\epsilon,s}(h,\varphi),
       \left(\frac{h}{(1+s^2)^{3/2}},\langle D\rangle^{-1}\varphi\right)
                         \right\rangle_{L^2}\\
 &=\int_{\R}\left[\frac{|\widehat h|^2}{(1+s^2)^{3/2}}
       +\left(\frac23+\frac{\xi^2}{3\langle \xi\rangle^2}\right)
                  \langle \xi\rangle^{-1}|\widehat\varphi|^2\right]\,d\xi
 +\epsilon\int_{\R}\left(\frac{\xi^2|\widehat h|^2}{(1+s^2)^3}
       +|\xi|\langle \xi\rangle^{-1}|\widehat\varphi|^2\right)\,d\xi\\
 &\quad+\int_{\R}\frac{|\xi|-\xi^2\langle \xi\rangle^{-1}}{(1+s^2)^{3/2}}
                 \operatorname{Re}(\widehat\varphi\,\overline{\widehat h})\,d\xi\,.
\end{aligned}
\end{equation}
The last term is the cross term we wanted to cancel. If we had used $|D|^{-1}$ in the pairing, this term would have cancelled. However, this multiplier is not bounded at the origin. We estimate the mismatch between $|D|^{-1}$ and $\langle D\rangle^{-1}$ by 
\begin{equation}\label{eq:mixed_coefficient}
 0\leq\frac{(|\xi|-\xi^2\langle \xi\rangle^{-1})/(1+s^2)^{3/2}}{(1+s^2)^{-3/4}\langle \xi\rangle^{-1/2}} = \frac{|\xi|}{(1+s^2)^{3/4}\langle\xi\rangle^{1/2}(\langle \xi\rangle+|\xi|)} \le \frac12\,.
\end{equation}
The equality uses $(\langle\xi\rangle-|\xi|)(\langle \xi\rangle+|\xi|)=1$. The final inequality follows from $\langle \xi\rangle\ge\max(1,|\xi|)$ and $1+s^2\ge1$. Now, \eqref{eq:mixed_coefficient} and Young's inequality give
\begin{equation}\label{eq:mixed_estimate}
\begin{aligned}
 \left|\int_{\R}\frac{|\xi|-\xi^2\langle \xi\rangle^{-1}}{(1+s^2)^{3/2}}
             \operatorname{Re}(\widehat\varphi\,\overline{\widehat h})\,d\xi\right|
 &\le\frac12\int_{\R}
       \frac{|\widehat h|}{(1+s^2)^{3/4}}
       \langle \xi\rangle^{-1/2}|\widehat\varphi|\,d\xi \\
 &\le\frac14\int_{\R}\left(\frac{|\widehat h|^2}{(1+s^2)^{3/2}}
                      +\langle \xi\rangle^{-1}|\widehat\varphi|^2\right)\,d\xi\,.
\end{aligned}
\end{equation}
Substituting \eqref{eq:mixed_estimate} into \eqref{eq:low_freq_pairing}, 
\begin{equation}\label{eq:low_freq_coercivity}
\begin{aligned}
 &\operatorname{Re}\left\langle-\mathcal L_{\epsilon,s}(h,\varphi),
               \left(\frac{h}{(1+s^2)^{3/2}},\langle D\rangle^{-1}\varphi\right)
                    \right\rangle_{L^2}\\
 &\qquad\gtrsim\int_{\R}\left(\frac{|\widehat h|^2}{(1+s^2)^{3/2}}+\langle \xi\rangle^{-1}|\widehat\varphi|^2\right)\,d\xi 
 + \epsilon\int_{\R}
       \left(\frac{\xi^2|\widehat h|^2}{(1+s^2)^3}
                       +|\xi|\langle \xi\rangle^{-1}|\widehat\varphi|^2\right)\,d\xi\,.
\end{aligned}
\end{equation}
For the forcing terms, Cauchy-Schwarz gives
\begin{equation}\label{eq:low_freq_forcing_pairing}
\begin{aligned}
 &\left|\left\langle(g_1,g_2),
       \left(\frac{h}{(1+s^2)^{3/2}},\langle D\rangle^{-1}\varphi\right)
                         \right\rangle_{L^2}\right|\\
 &\quad\le\left(\frac{\|g_1\|_{L^2}^2}{(1+s^2)^{3/2}}
                         +\|g_2\|_{H^{-1/2}}^2\right)^{1/2}
       \left(\frac{\|h\|_{L^2}^2}{(1+s^2)^{3/2}}
                         +\|\varphi\|_{H^{-1/2}}^2\right)^{1/2}\,.
\end{aligned}
\end{equation}
Thus, combining \eqref{eq:low_freq_coercivity} and \eqref{eq:low_freq_forcing_pairing}, discarding the $\epsilon$ terms and using $\|g_2\|_{H^{-1/2}}\le\|g_2\|_{L^2}$, we obtain the low-frequency estimate
\begin{equation}\label{eq:low_freq_estimate}
 \|h\|_{L^2}^2+\|\varphi\|_{H^{-1/2}}^2\lesssim_s\|g_1\|_{L^2}^2+\|g_2\|_{L^2}^2\,.
\end{equation}

  \textit{Step 2:} We now bound $h$ in $H^3$ and $\varphi$ in $H^2$. To do this we argue in the same way and pair the equation with $\langle D\rangle^4h/(1+s^2)^{3/2}$ and $\langle D\rangle^3\varphi$. Integration by parts in the two dilation terms now gives
\begin{equation*}
\begin{aligned}
 &\operatorname{Re}\int_{\R}\left(\frac43\widehat h
       +\frac23\xi\partial_\xi\widehat h\right)
       \frac{\langle \xi\rangle^4\overline{\widehat h}}{(1+s^2)^{3/2}}\,d\xi
        \quad=\int_{\R}\left(\frac43\langle \xi\rangle^4
       -\frac13\partial_\xi(\xi\langle \xi\rangle^4)\right)
       \frac{|\widehat h|^2}{(1+s^2)^{3/2}}\,d\xi\\
 &\quad=\int_{\R}\langle \xi\rangle^4\frac{|\widehat h|^2}{(1+s^2)^{3/2}}\,d\xi
       -\frac43\int_{\R}\xi^2\langle \xi\rangle^2
       \frac{|\widehat h|^2}{(1+s^2)^{3/2}}\,d\xi\,,\\
 &\operatorname{Re}\int_{\R}\left(\widehat\varphi
       +\frac23\xi\partial_\xi\widehat\varphi\right)
       \langle \xi\rangle^3\overline{\widehat\varphi}\,d\xi
       \quad=\int_{\R}\left(\langle \xi\rangle^3
       -\frac13\partial_\xi(\xi\langle \xi\rangle^3)\right)
       |\widehat\varphi|^2\,d\xi\\
 &\quad=\int_{\R}\langle \xi\rangle^4
       \left(\frac23+\frac{\xi^2}{3\langle \xi\rangle^2}\right)
       \langle \xi\rangle^{-1}|\widehat\varphi|^2\,d\xi
       -\frac43\int_{\R}\xi^2\langle \xi\rangle
       |\widehat\varphi|^2\,d\xi\,.
\end{aligned}
\end{equation*}
The two negative terms come from integrating by parts through the weight. Working now as in the estimate \eqref{eq:low_freq_coercivity} with the additional factor $\langle \xi \rangle^{4}$ we obtain 
\begin{align}\label{eq:high_freq_coercivity}
 &\int_{\R}\langle \xi\rangle^4
       \left(\frac{|\widehat h|^2}{(1+s^2)^{3/2}}
                    +\langle \xi\rangle^{-1}|\widehat\varphi|^2\right)\,d\xi
 +\epsilon\int_{\R}\langle \xi\rangle^4
       \left(\frac{\xi^2|\widehat h|^2}{(1+s^2)^3}
              +|\xi|\langle \xi\rangle^{-1}|\widehat\varphi|^2\right)\,d\xi
       \nonumber\\
 &\qquad\lesssim
       \left|\left\langle (g_1,g_2),
       \left(\frac{\langle D\rangle^4h}{(1+s^2)^{3/2}},
                          \langle D\rangle^3\varphi\right)\right\rangle_{L^2}\right|
  +\int_{\R}\xi^2\langle \xi\rangle^2
       \left(\frac{|\widehat h|^2}{(1+s^2)^{3/2}}
                    +\langle \xi\rangle^{-1}|\widehat\varphi|^2\right)\,d\xi\,.
\end{align}
Since $\xi^2\langle \xi\rangle^2\le\langle \xi\rangle^4$, the last term in \eqref{eq:high_freq_coercivity} is bounded by $\|h\|_{H^2}^2+\|\varphi\|_{H^{3/2}}^2$. We now define the following intermediate quantities to simplify notation:
\begin{equation*}
 T=\|h\|_{H^3}^2+\|\varphi\|_{H^2}^2\,, \qquad D=\int_{\R} \langle \xi\rangle^4\left(\frac{\xi^2|\widehat h|^2}{(1+s^2)^3}+|\xi|\langle \xi\rangle^{-1}|\widehat\varphi|^2\right)\,d\xi\,.
\end{equation*}
By comparing the weights at frequencies $|\xi|\ge1$ and $|\xi|\le1$,
\begin{equation}\label{eq:T_D_comparison}
 T\lesssim_s D+\int_{\R}\left(\frac{|\widehat h|^2}{(1+s^2)^{3/2}} + \langle \xi\rangle^{-1}|\widehat\varphi|^2\right)\,d\xi\,.
\end{equation}
Splitting at a sufficiently large frequency also gives, for every $\eta>0$, the following interpolation estimate used to control the last term in \eqref{eq:high_freq_coercivity},
\begin{equation}\label{eq:extra_term_interpolation}
 \|h\|_{H^2}^2+\|\varphi\|_{H^{3/2}}^2 \leq \eta T+C_{\eta,s}\int_{\R}\left(\frac{|\widehat h|^2}{(1+s^2)^{3/2}}+\langle \xi\rangle^{-1}|\widehat\varphi|^2\right)\,d\xi\,.
\end{equation}

The forcing terms satisfy
\begin{equation*}
 \left|\left\langle g_1,\frac{\langle D\rangle^4h}{(1+s^2)^{3/2}}\right\rangle_{L^2}\right|\le\|g_1\|_{H^1}\|h\|_{H^3}\,,
 \qquad 
 |\langle g_2,\langle D\rangle^3\varphi\rangle_{L^2}|\le\|g_2\|_{H^1}\|\varphi\|_{H^2}\,.
\end{equation*}
Substituting this and \eqref{eq:extra_term_interpolation} into \eqref{eq:high_freq_coercivity} and discarding the first nonnegative term, we obtain
\begin{equation*}
 \epsilon D\lesssim_s\eta T+C_{\eta,s}\int_{\R}\left(\frac{|\widehat h|^2}{(1+s^2)^{3/2}}+\langle \xi\rangle^{-1}|\widehat\varphi|^2\right)\,d\xi+(\|g_1\|_{H^1}+\|g_2\|_{H^1})T^{1/2}\,.
\end{equation*}
Multiplying by $\epsilon$ the estimate \eqref{eq:T_D_comparison} now gives
\begin{align*}
 \epsilon T
 &\lesssim_s\epsilon D+\epsilon\int_{\R}\left(
       \frac{|\widehat h|^2}{(1+s^2)^{3/2}}
       +\langle \xi\rangle^{-1}|\widehat\varphi|^2\right)\,d\xi\\
 &\lesssim_s\eta T+(C_{\eta,s}+\epsilon)\int_{\R}\left(
       \frac{|\widehat h|^2}{(1+s^2)^{3/2}}
       +\langle \xi\rangle^{-1}|\widehat\varphi|^2\right)\,d\xi
       +(\|g_1\|_{H^1}+\|g_2\|_{H^1})T^{1/2}\,.
\end{align*}
The implicit constants are independent of $\eta$ and $\epsilon$. Choose $\eta>0$ sufficiently small in terms of $\epsilon,s$ to absorb the $\eta T$ contribution into the left-hand side. Young's inequality absorbs the $T^{1/2}$ contribution as well. Together with the low-frequency estimate \eqref{eq:low_freq_estimate}, we obtain 
\begin{align*}
 T&\lesssim_{\epsilon,s}\int_{\R}\left(
       \frac{|\widehat h|^2}{(1+s^2)^{3/2}}
       +\langle \xi\rangle^{-1}|\widehat\varphi|^2\right)\,d\xi
       +(\|g_1\|_{H^1}+\|g_2\|_{H^1})^2\\
 &\lesssim_{\epsilon,s}(\|g_1\|_{H^1}+\|g_2\|_{H^1})^2\,.
\end{align*}
Taking square roots yields the high frequency estimate
\begin{equation*}
  \|h\|_{H^3}+\|\varphi\|_{H^2}\lesssim_{\epsilon,s}\|g_1\|_{H^1}+\|g_2\|_{H^1}\,.
\end{equation*}
The weighted components of the $X$ norm follow from
\begin{align*}
 \frac23x h'&=\frac23h
       -\frac{\epsilon h''}{(1+s^2)^{3/2}}
       +\Lambda\varphi+g_1\,,\\
 \frac23x\varphi'&=\frac{h''}{(1+s^2)^{3/2}}
       +\frac13\varphi+\epsilon\Lambda\varphi+g_2\,,
\end{align*}
since every term on the right is bounded in $H^1$. In total, we have the a priori estimate
\begin{equation}\label{eq:inverse_estimate}
 \|(h,\varphi)\|_X\lesssim_{\epsilon,s}\|(g_1,g_2)\|_Y\,.
\end{equation}
This a priori estimate is valid for smooth, compactly supported solutions. The steps related to integration by parts for functions in $X$ can be justified by introducing a cutoff function $\chi_{r,R}$ with support in $[-r, R]$ and passing to the limit as $r \to 0$ and $R \to \infty$.

\textit{Step 3:} We now show that $\mathcal L_{\epsilon, s}$ defines an isomorphism. First we show that for $(g_1, g_2) \in Y$, the linear equation has a solution in $X$. In Fourier variables, the equation $\mathcal L_{\epsilon,s}(h,\varphi)=(g_1,g_2)$ reads, for $\xi>0$,
\begin{align*}
 \frac23\xi\partial_\xi\widehat h+
 \left(\frac43+\frac{\epsilon\xi^2}{(1+s^2)^{3/2}}\right)\widehat h
 +\xi\widehat\varphi&=-\widehat g_1\,,\\
 \frac23\xi\partial_\xi\widehat\varphi
 -\frac{\xi^2}{(1+s^2)^{3/2}}\widehat h
 +(1+\epsilon\xi)\widehat\varphi&=-\widehat g_2\,.
\end{align*}
For fixed $r>0$, we can solve this linear ODE on $[r,\infty)$ with initial conditions  $\widehat{h}_r(r) = \widehat{\varphi}_r(r) = 0$. We set this solution to be zero for $0 < \xi < r$ and extend it evenly. Now fix $R>r$ and apply the estimates of Steps 1 and 2 on the interval $(r,R)$ to get those same bounds uniformly in $r, R$. The boundary terms at $r$ vanish, and the boundary term at $R$ for the low-frequency estimate is
\begin{equation*}
  \frac{R}{3} \left( \frac{|\widehat h_r(R)|^2}{(1+s^2)^{3/2}} + \langle R \rangle^{-1} |\widehat\varphi_r(R)|^2 \right) \geq 0,
\end{equation*}
so it can be discarded. The boundary term coming from the high-frequency estimate can also be discarded in the same way. Taking $R \to +\infty$ gives a solution $(h_r,\varphi_r)\in X$ for the equation with right hand side $\mathbf1_{\{|\xi|>r\}}(\widehat g_1,\widehat g_2)$. The inverse Fourier transform of this term converges in $Y$ to $(g_1,g_2)$, so using \eqref{eq:inverse_estimate} gives that $(h_r,\varphi_r)$ converges in $X$ as $r \to 0$ to a solution with right-hand side $(g_1,g_2)$. This shows existence of solutions. Estimate \eqref{eq:inverse_estimate} gives uniqueness and boundedness of the inverse.
\end{proof}

\begin{lem}[Compactness of $\mathcal{K}$]\label{lem:K_compact}
Let $\epsilon>0$, $s\in\R$, and $(h,\varphi)\in X$. If $(sf_{\mathrm{ref}}+h,\varphi)$ is admissible, then the operator $\mathcal K$ in \eqref{eq:K} is compact from $X$ to $Y$.  
\end{lem}
\begin{proof}
  We use the notation $X = X_h \times X_\varphi$ to refer to the two factors in Definition \ref{def:X_Y_spaces} corresponding to $h$ and $\varphi$. The operator $\mathcal{K}$ is a sum of terms in \eqref{eq:K}. We show that each one of them is compact. First, the operator $\tilde{h} \mapsto (\tilde{h})'$ is compact from $X_h$ to $H^1$ because of local Rellich compactness and the control of the term $\|x\tilde{h}'\|_{H^1}$ in the $X_h$ norm. Therefore, the operators
  \begin{equation*}
    \tilde{h} \mapsto \frac{3f'f''}{(1+f'^2)^{5/2}}(\tilde{h})', \qquad \tilde{h} \mapsto \Psi_x(\tilde{h})'
  \end{equation*}
  are compact from $X_h$ to $H^1$ because the coefficients define bounded multiplication operators.

  Now we show that the operator $\tilde{h} \mapsto [G(f), \tilde{h}]\Psi_y$ is compact from $X_h$ to $H^1_e$. Lemma \ref{lem:DtN_commutator} gives
  \begin{equation}\label{eq:DtN_commutator_velocity}
    \|[G(f),\tilde h]\Psi_y\|_{H^1} \lesssim \|(\tilde h)'\|_{W^{1,\infty}}\|\Psi_y\|_{H^1}.
  \end{equation}
  This is the Dirichlet-to-Neumann version of the classical Coifman--Meyer commutator estimate \cite{CoifmanMeyer1978}. Since $\tilde h\mapsto(\tilde h)'$ is compact from $X_h$ to $H^1$ and bounded into $H^2$, by interpolation and the embedding $H^{7/4}\hookrightarrow W^{1,\infty}$, it is also compact into $W^{1,\infty}$. Therefore the bound \eqref{eq:DtN_commutator_velocity} gives that $\tilde h\mapsto[G(f),\tilde h]\Psi_y$ is compact from $X_h$ to $H^1_e$. Since multiplication by $\Psi_y$ is bounded on $H^1$, the operator $\tilde h\mapsto\Psi_y[G(f),\tilde h]\Psi_y$ is also compact.

The remaining terms in $\mathcal{K}$ concern the operator $G(f) - \Lambda: X_\varphi \to H^1_e$. We will show this operator is compact. First suppose $f=sf_{\mathrm{ref}}+h_0$ with $h_0\in C_c^\infty$ even. In this case, using \eqref{eq:DtN_integral_formula} and \eqref{eq:half_laplacian_integral_formula}, we can express it as the integral operator 
  \begin{equation*}
    (G(f) - \Lambda)\tilde{\varphi}(x) = \Re \int_{\mathbb{R}} \partial_x K_f(x,y) (\tilde{\varphi}(y) + i \psi(y)) \, dy,
  \end{equation*}
  where the function $\psi$ is the trace at the boundary of the harmonic conjugate in $L^{2}$ of the harmonic extension of $\tilde{\varphi}$ in $\Omega_f$. The kernel is 
  \begin{equation*}
     K_f(x,y) = \frac{1}{\pi} \left[ \frac{1 + if'(y)}{(x-y) + i(f(x)-f(y))} - \frac{1}{x-y} \right] = \frac{i\bigl(f'(y)(x-y)-f(x)+f(y)\bigr)}{\pi\bigl(x-y+i(f(x)-f(y))\bigr)(x-y)}.
  \end{equation*}
  This kernel is smooth since the numerator has a zero of order two at $x=y$. By taking derivatives and using that outside of a compact interval $f(x) = s|x|$, we obtain the bounds
  \begin{equation*}
    \left| \partial_x K_f(x,y) \right| + \left| \partial_x^2 K_f(x,y) \right| \lesssim \frac{1}{( 1 + \left| x \right| + \left| y \right|)^2}
  \end{equation*}
  Thus, the kernel defines a Hilbert-Schmidt operator from $L^2$ to $H^1$. Since we have the bound $\| \tilde{\varphi} + i \psi \|_{L^{2}} \lesssim \| \tilde{\varphi} \|_{L^{2}}$, we conclude that $G(f) - \Lambda$ is compact from $X_\varphi$ to $H^1_e$ in this case. 

  For a general admissible graph $f$, we can approximate it by a sequence of graphs $f_n = sf_{\mathrm{ref}} + h_n$ with $h_n \in C_c^\infty$ even and $\| f_n' - f' \|_{W^{1,\infty}} \to 0$. The corresponding operators $G(f_n) - \Lambda$ converge to $G(f) - \Lambda$ in operator norm, since Lemma \ref{lem:DtN_analytic} gives $G(f_n) \to G(f)$ in $\mathcal{L}(H^{2}, H^{1})$, hence also in $\mathcal{L}(X_\varphi, H^{1})$. This shows $G(f) - \Lambda : X_\varphi \to H^1_e$ is compact and thus the operator $\tilde{\varphi} \mapsto \Psi_y(G(f) - \Lambda) \tilde{\varphi}$ is also compact.
\end{proof}

  It remains to treat the term $\mathcal R$. The key point is that, when we add  $\epsilon$ times the second equation to the first, there is a cancellation in $\mathcal{R}$ that allows us to write the highest order terms as a multiplication by  $((1+s^2)/(1+f'^2))^{3/2}$. This can be inverted because it is a function bounded away from zero. Every other term in the remaining expression is compact, so it does not contribute to the Fredholm index.

\begin{prop}\label{prop:fredholm_index_zero}
  Fix $\epsilon > 0$ and $s \in \mathbb{R}$. Let $(h,\varphi) \in X$ be such that the pair $(sf_{\mathrm{ref}} + h, \varphi)$ is admissible. Then the linearized operator $D_{(h,\varphi)}\mathcal{F}_{\epsilon}(s, h, \varphi): X \to Y$ is a Fredholm operator of index zero.      
\end{prop}
\begin{proof}
  We use again the notation $X = X_h \times X_\varphi$. By Lemma \ref{lem:inverse_linear}, the inverse linear operator acting on right-hand sides of the form $(-\epsilon g, g)$,
  \begin{equation*}
    \mathcal S_{\epsilon,s}:H^1_e\to X,
    \qquad
    \mathcal S_{\epsilon,s}g:=\mathcal L_{\epsilon,s}^{-1}
    \begin{pmatrix}-\epsilon g\\g\end{pmatrix}
  \end{equation*}
  is bounded. Write $\mathcal S=\mathcal S_{\epsilon,s}$ for the remainder of the proof. The decomposition \eqref{eq:derivative_decomposition} becomes
  \begin{equation*}
    D_{(h,\varphi)}\mathcal F_\epsilon(s,h,\varphi) = \mathcal L_{\epsilon,s}(I+\mathcal S\mathcal R)+\mathcal K.
  \end{equation*}
  Since $\mathcal L_{\epsilon,s}$ is invertible and $\mathcal K$ is compact by Lemma \ref{lem:K_compact}, it suffices to prove that $I+\mathcal S\mathcal R$ is Fredholm of index zero.

  The $I+AB$, $I+BA$ principle \cite[Theorem~6]{Barnes1998} shows that
  \begin{equation*}
    I + \mathcal S\mathcal R \qquad \text{and} \qquad I + \mathcal R\mathcal S
  \end{equation*}
  are Fredholm simultaneously and have the same index. We study the second operator, which acts on the simpler space $H^1_e$.  Because of the highest order derivative terms, the operator $\mathcal{R} \mathcal{S}$ need not be compact. We decompose $(I + \mathcal{R} \mathcal{S})$ into an invertible multiplication operator and a compact remainder. Given $g \in H^1_e$, denote by $(\tilde{h}_g, \tilde{\varphi}_g) = \mathcal{S} g \in X$ the solution of $\mathcal{L}_{\epsilon,s}(\tilde{h}_g, \tilde{\varphi}_g) = (-\epsilon g, g)$. Using the first component in \eqref{eq:L}, we get
  \begin{equation*}
    (\tilde{h}_g)'' =  - (1 + s^2)^{3/2} g + \frac{(1 + s^{2})^{3/2}}{\epsilon} \left( \frac{2}{3} (\tilde{h}_g - x (\tilde{h}_g)') + \Lambda \tilde{\varphi}_g \right).
  \end{equation*}
  Substituting this into $(I + \mathcal R\mathcal S)g$ we obtain, using \eqref{eq:R}, the following expression for this operator:
  \begin{equation}\label{eq:RS_decomposition}
  \begin{aligned}
    g + \mathcal{R}(\tilde{h}_g, \tilde{\varphi}_g) =& \left( \frac{1+s^{2}}{1+f'^2} \right)^{3/2} g + \Psi_x (\tilde{\varphi}_g)' - \Psi_y \Lambda \tilde{\varphi}_g \\
    +& \frac{1}{\epsilon} \left( 1 + s^{2} \right)^{3/2} \left( \left( 1 + s^{2} \right)^{-3/2} - \left( 1+f'^2 \right)^{-3/2} \right) \left( \frac{2}{3} (\tilde{h}_g - x (\tilde{h}_g)') + \Lambda \tilde{\varphi}_g  \right).
  \end{aligned}
\end{equation}
The first term is a multiplication by $((1+s^{2})/(1+f'^2))^{3/2}$. Since $f$ is admissible, this term has a positive lower bound and belongs to $W^{1,\infty}$ so it defines a bounded, invertible operator on $H^1_e$. Now we show that every other term in \eqref{eq:RS_decomposition} is compact, and this will finish the proof.

We show that $g \mapsto (\tilde{\varphi}_g)'$ is compact. We use cancellation in the equation to get some extra regularity for $\tilde{\varphi}_g$. Multiplying the second equation in $\mathcal{L}_{\epsilon,s}(\tilde{h}_g, \tilde{\varphi}_g) = (-\epsilon g, g)$ by $\epsilon$ and adding it to the first equation gives
\begin{equation*}
  (1+\epsilon^2)\Lambda\tilde{\varphi}_g=-\frac{2}{3}\bigl(\tilde{h}_g-x(\tilde{h}_g)'\bigr)-\frac{\epsilon}{3}\tilde{\varphi}_g + \frac{2\epsilon}{3}x(\tilde{\varphi}_g)'.
\end{equation*}
Differentiating twice, taking the $L^{2}$ inner product with $(\tilde{\varphi}_g)''$ and integrating by parts the term $\langle (\tilde{\varphi}_g)'', x(\tilde{\varphi}_g)'''\rangle$ using a regularization argument gives
\begin{equation*}
  (1 + \epsilon^2) \| (\tilde{\varphi}_g)'' \|_{\dot{H}^{1/2}}^2 = \frac{2}{3}\langle (\tilde{h}_g)'', (\tilde{\varphi_g})'' \rangle + \frac{2}{3} \langle (\tilde{h}_g)''', x(\tilde{\varphi_g})'' \rangle + \frac{2\epsilon}{3} \| (\tilde{\varphi}_g)'' \|_{L^{2}}^{2}.
\end{equation*}
The right hand side is controlled using Cauchy-Schwarz and the norm of $X$. Consequently, we have $\| \tilde{\varphi}_g \|_{H^{5/2}} \lesssim_{\epsilon, s} \| g \|_{H^{1}}$. Since $(\tilde{\varphi}_g)'$ is bounded in $H^{3/2}$, it is compact in $H^1$ by Rellich compactness and the control of $\|x(\tilde{\varphi}_g)'\|_{H^1}$ in the $X_\varphi$ norm. Therefore, the operator $g \mapsto (\tilde{\varphi}_g)'$ is compact from $H^1_e$ to $H^1$. Since multiplication by $\Psi_x$ is bounded in $H^1$, the operator $g \mapsto \Psi_x(\tilde{\varphi}_g)'$ is compact. In the same way, the operator $g \mapsto \Psi_y \Lambda \tilde{\varphi}_g$ is compact because $\Lambda = H\partial_x$.

Now, the operator $\tilde{h}_g \mapsto \tilde{h}_g - x(\tilde{h}_g)'$ is bounded from $X_h$ to $H^1$ and from $X_h$ to $H^2([-R, R])$ for $R > 0$. Let $\chi_R = \chi(\cdot/R)$ be a smooth, even cutoff function with support in $[-R,R]$. The coefficient $((1+s^{2})^{-3/2} - (1+f'^2)^{-3/2})$ belongs to $H^1$ and thus
\begin{align*}
&\| (1 - \chi_R) ((1+s^{2})^{-3/2} - (1+f'^2)^{-3/2}) (\tilde{h}_g - x(\tilde{h}_g)') \|_{H^{1}} \\
&\qquad \lesssim \| (1 - \chi_R) ((1+s^{2})^{-3/2} - (1+f'^2)^{-3/2}) \|_{H^1} \| \tilde{h}_g - x(\tilde{h}_g)' \|_{H^1} \longrightarrow 0 \text{ as } R \to \infty.
\end{align*}
By Rellich compactness, $ \tilde{h}_g \mapsto ((1+s^{2})^{-3/2} - (1+f'^2)^{-3/2}) (\tilde{h}_g - x(\tilde{h}_g)')$ can be approximated by compact operators, and thus is compact from $X_h$ to $H^1_e$. Therefore $g \mapsto ((1+s^{2})^{-3/2} - (1+f'^2)^{-3/2}) (\tilde{h}_g - x(\tilde{h}_g)')$ is also compact.
\end{proof}

\color{black}
\subsection{Uniform interval of existence in $s$ for all $\epsilon$}

The regularized system always has the trivial solution $(s, h, \varphi) = (0,0,0)$. Now we use a continuation argument to get existence of solutions for all $|s| \leq s_0$, with $s_0$ uniform in $\epsilon$. We introduce a degree for the map $\mathcal{F}_\epsilon$ to track the existence of a solution in $s$ and prove this quantity is locally constant. Before we state and prove this result, we present here the concepts from degree theory that we will use. We follow the construction in \cite{PejsachowiczRabier1998, BenevieriFuriPera2000} and, for the standard Brouwer degree, we refer to \cite[Chapter 1]{Deimling1985}.

  Given $F: U \subset X \to Y$ a $C^1$ Fredholm map of index zero and a compact zero set $K \subset U$, we can always choose a finite-dimensional subspace $E \subset Y$ such that
  \begin{equation*}
    \operatorname{ran} D F(u) + E = Y, \qquad \text{ for $u$ in a neighborhood $W$ of $K$}.
  \end{equation*}
  The set $M_E = W \cap F^{-1}(E)$ is a $C^1$ manifold of dimension $\dim E$. We can now consider the map $F|_{M_E}$ restricted to this manifold and, since this is finite-dimensional, we can define the degree through the classical Brouwer degree of the coordinate map. If $O$ is an open set with no zeros on $\partial O$ in $M_E$ and $\phi$ is a coordinate chart whose domain contains $\overline{O}$, then we set 
\begin{equation*}
  \deg_2 (F|_{M_E}, O, 0) = \deg_B (F|_{M_E} \circ \phi^{-1}, \phi(O), 0) \mod 2.
\end{equation*}
Here $\deg_B$ is the classical Brouwer degree for a map defined in $\mathbb{R}^m$, since we can identify $E \cong \mathbb{R}^{m}$. This definition is independent of the choices of $\phi$, $E$ and $W$. If $0$ is a regular value of $F$ and $F^{-1}(0)$ is compact, we can enclose its points in disjoint open sets $O_i$ with coordinate charts $\phi_i$ and define the degree as
\begin{equation*}
  \deg_2 (F|_{M_E}, M_E, 0) = \sum_i \deg_2 (F|_{M_E}, O_i, 0) \mod 2.
\end{equation*}
Whenever $0$ is a singular value, the definition can be adapted by perturbing this value. Choose a relatively compact open neighborhood $O \subset M_E$ of $K$, with no zeros on $\partial O$, and perturb $F|_{M_E}$ to a sufficiently close regular value $e \in E$, counting only preimages in $O$. The degree is then defined as
\begin{equation*}
  \deg_2 (F|_{M_E}, O, 0) = \deg_2 (F|_{M_E}, O, e) \mod 2.
\end{equation*}
For a map $F: U \to Y$ and an open set $V \subset U$ for which the zero set $K=F^{-1}(0) \cap V$ is compact, choose the reduction neighborhood $W \subset V$ around $K$ and define the degree as
\begin{equation*}
  \deg_2 (F, V, 0) = \deg_2 (F|_{M_E}, M_E, 0).
\end{equation*}
The main properties of this degree are stated in the following proposition, which is a consequence of the construction in \cite{PejsachowiczRabier1998}.
\begin{prop}\label{prop:degree_properties}
  Let $F:U \to Y$ be a $C^1$ Fredholm map of index zero, where $U \subset X$ is open. For every $V \subset U$ open with compact $K = F^{-1}(0) \cap V$, the degree $\deg_2(F, V, 0)$ is well-defined and has the following properties: 
  \begin{enumerate}[label=\alph*.,ref=\alph*]
  \item\label{degree_properties:normalization} If $DF(u)$ is invertible at every $u\in K$, then $K$ is finite and
  \begin{equation*}
    \deg_2(F,V,0)=\#K\pmod2\,.
  \end{equation*}
  \item\label{degree_properties:existence} If $\deg_2(F,V,0) = 1$, then $K\neq\emptyset$.
  \item\label{degree_properties:excision} If $V_0\subset V$ is open and contains $K$, then $\deg_2(F,V_0,0)=\deg_2(F,V,0)$. If we have disjoint open sets $V_1,\ldots,V_j\subset V$ whose union contains $K$, then
  \begin{equation*}
    \deg_2(F,V,0) = \sum_{i=1}^j\deg_2(F,V_i,0) \mod 2.
  \end{equation*}
\item\label{degree_properties:homotopy} Let $H:[0,1]\times V\to Y$ extend to a $C^1$ map on an open neighborhood of $[0,1]\times V$. Suppose each $D_u H(t,u)$ is Fredholm of index zero and $K=\{(t,u)\in[0,1]\times V:H(t,u)=0\}$ is compact in $[0,1]\times V$. Then 
  \begin{equation*}
    \deg_2(H(0,\cdot),V,0)=\deg_2(H(1,\cdot),V,0)\,.
  \end{equation*}
\end{enumerate} 
\end{prop}

Now that we have the necessary tools from degree theory, we can prove the existence of solutions to the regularized system \eqref{eq:profile_reg} for $|s| \leq s_0$, with $s_0$ independent of $\epsilon$.

\begin{prop}\label{prop:existence_small_unif}
  There exist constants $s_0>0$ and $C<\infty$, independent of $\epsilon$ and $s$, such that for every $0<\epsilon\leq1$ and every $|s| \leq s_0$, the regularized system \eqref{eq:profile_reg} has a real even solution $(f,\varphi)$ satisfying
  \begin{equation*}
    (f-sf_{\mathrm{ref}},\varphi) \in X,
    \qquad 
    M(f,\varphi)\leq C|s|.
  \end{equation*}
\end{prop}
\begin{proof}
  With $\delta$ fixed as above, choose $s_0$ small enough so that $C s_0 < \delta/2$, where $C$ is the constant in Proposition \ref{prop:uniform_slope_bounds}. The functional $\mathcal{F}_\epsilon$ can be defined on the set
  \begin{equation*}
    V = \left\{ (s, h, \varphi) \in \mathbb{R} \times X : M(sf_{\mathrm{ref}} + h, \varphi) < \delta \right\}.
  \end{equation*}
  This set is open by continuity of $M$. By Proposition \ref{prop:compactness_solution_set}, the set of solutions 
  \begin{equation*}
    K_\epsilon = \left\{ (s, h, \varphi) \in \left[ -s_0, s_0 \right] \times X : \mathcal{F}_\epsilon(s, h, \varphi) = 0, \: M(sf_{\mathrm ref} + h, \varphi) \leq \delta \right\}
  \end{equation*}
  is compact. At a fixed $s$, define the slices of the solution set by $K_{\epsilon, s} = \left\{ (h, \varphi) : (s, h, \varphi) \in K_\epsilon \right\}$. In the same way, consider the slices of the domain $V_s = \left\{ (h, \varphi) : (s, h, \varphi) \in V \right\}$. By Proposition \ref{prop:fredholm_index_zero} the map $\mathcal{F}_\epsilon(s, \cdot, \cdot): V_s \to Y$ is Fredholm of index zero. Moreover, this map is $C^1$ and its zero set $K_{\epsilon, s}$ is compact, so we can define the mod-two degree
  \begin{equation*}
    \deg_2 \left( \mathcal{F}_\epsilon(s, \cdot, \cdot), V_s, 0 \right), \qquad \text{ for } \left| s \right| \leq s_0.
  \end{equation*}

  We now show that this degree is locally constant in $s$. We will use homotopy invariance of the degree (Proposition \ref{prop:degree_properties}\ref{degree_properties:homotopy}), and this requires that no solutions cross the boundary of $V_s$ as $s$ varies. The choice of $s_0$ such that $Cs_0 < \delta/2 $ keeps all the solutions away from the boundary of $V_s$ since Proposition \ref{prop:uniform_slope_bounds} gives $M(sf_{\mathrm{ref}}+h,\varphi)\leq C|s|<\delta/2$.

  Fix $s_* \in \left[ -s_0, s_0 \right]$. If the solution set $K_{\epsilon, s_*}$ is empty, then the degree is zero (Proposition \ref{prop:degree_properties}\ref{degree_properties:existence}). Compactness of solutions implies that for $s$ sufficiently close to $s_*$, the set $K_{\epsilon, s}$ is also empty, so the degree is zero in a neighborhood of $s_*$. We focus on the case $K_{\epsilon, s_*} \neq \emptyset$. We will apply the homotopy invariance of the degree, but this property compares maps on a fixed open set while here $V_s$ varies with $s$. To fix this, we use compactness of $K_{\epsilon, s_*}$ to take an open set $V_* \subset X$ and a relative neighborhood in $[-s_0, s_0]$ of $s_*$, denoted by $I$, such that $K_{\epsilon, s_*} \subset V_*$ and $\overline{V_*} \subset V_s$. After shrinking $I$ if necessary, compactness of $K_\epsilon$ implies that $K_{\epsilon, s} \subset V_*$ for all $s \in I$.

  For every $s \in I$, we have $K_{\epsilon, s} \subset V_* \subset V_s$. By excision (Proposition \ref{prop:degree_properties}\ref{degree_properties:excision}), the degree agrees on $V$ and $V_s$. Thus we compute
  \begin{equation*}
    \deg_2 \left( \mathcal{F}_\epsilon(s, \cdot, \cdot), V_s, 0 \right) = \deg_2 \left( \mathcal{F}_\epsilon(s, \cdot, \cdot), V_*, 0 \right).
  \end{equation*}
  For $J \subset I$ closed, the zero set on $J \times V_*$ is $K_\epsilon \cap (J \times X)$ and this is compact. Now we can use homotopy invariance (Proposition \ref{prop:degree_properties}\ref{degree_properties:homotopy}) to conclude that $\deg_2 \left( \mathcal{F}_\epsilon(s, \cdot, \cdot), V_*, 0 \right)$ is constant for $s \in I$. This shows that the degree is locally constant in $s$.

  In total, the function $s \mapsto \deg_2 \left( \mathcal{F}_\epsilon(s, \cdot, \cdot), V_s, 0 \right)$ takes values in $\mathbb{Z}_2$, is locally constant, and defined on a connected domain, so it is constant. We have that $D_{(h, \varphi)} \mathcal{F}_\epsilon (0,0,0) = \mathcal{L}_{\epsilon, 0}$ and this operator is invertible by Lemma \ref{lem:inverse_linear}. Thus, for $s = 0$, Proposition \ref{prop:degree_properties}\ref{degree_properties:normalization} yields that the degree is one since the inequality $M(f, \varphi) \leq 0$ gives that the only solution there is $(f,\varphi) = (0,0)$. We conclude that for all $s \in \left[ -s_0, s_0 \right]$ the degree is also one. In particular, for every $s \in \left[ -s_0, s_0 \right]$ there exists a solution $(h, \varphi) \in V_s$ to $\mathcal{F}_\epsilon(s, h, \varphi) = 0$.
\end{proof}

\section{Proof of the main theorem: limit \texorpdfstring{$\epsilon$ to $0$}{ε to 0} and corollaries}\label{sec:main_theorem_proof}
The goal of this section is to prove the main Theorem \ref{thm:main_theorem} and its corollaries.

\subsection{Proof of the main theorem}

\begin{proof}[ Proof of \ref{thm:main_theorem}]
  We are going to divide this proof in five steps: In Step 1 we are going to find a strong enough converging subsequence, specifically $(f_j, \varphi_j)$ converging in $C^1_{loc}\times C^0_{loc}$. In Step 2 we are going to prove the limit satisfies the kinematic equation distributionally, in Step 3 the same but for the dynamic equation. Finally in Step 4 we will prove that the solutions have $C^\infty$ local regularity, hence the equations are satisfied in the strong sense, and in Step 5 that the energy and the uniform slope and weighted velocity bounds are satisfied in the limit too.

  In all the proof $s$ is going to be a fixed number such that $|s|\leq s_0$.  
  Taking the sequence $\epsilon_j = 1/j$, since $|s|\leq s_0$ by Proposition \ref{prop:existence_small_unif} we have an even solution $(f_j, \varphi_j)$ with $(f_j-sf_{\mathrm{ref}},\varphi_j)\in X$ to the $\epsilon_j$ regularized system \eqref{eq:profile_reg}, such that
  \begin{equation}\label{eq:energy_bounds_sec7}
    \Hone{f_j - s|x|} + \Honehalf{\varphi_j} \lesssim |s|,\quad \text{ and } \quad 
    \Linf{f_j'} + \Linf{\frac{\nabla \Psi_j (x, f_j(x))}{\langle x \rangle}} \lesssim |s|,
  \end{equation}
  where $\Psi_j$ is the harmonic extension to $y<f_j(x)$ with trace $\varphi_j$. Moreover, the energy estimates of Proposition \ref{prop:energy_bounds} give
  \begin{equation}\label{eq:reg_energy_bounds_sec7}
    \epsilon_j( \Ltwo{G(f_j)\varphi_j} + \Ltwo{\kappa_j - N_j}) \lesssim \sqrt{\epsilon_j} |s| \quad \text{ and } \quad \int_{\Omega_{f_j}} |\nabla \Psi_j|^2 \lesssim s^2,
  \end{equation}
  where $\kappa_j = \kappa[f_j]$ is the curvature and $N_j = N[f_j, \varphi_j]$ is the nonlinearity in the Bernoulli equation defined in \eqref{def:curvature_N}. During the proof we are going to take multiple subsequences and we are going to abuse notation and not relabel them.

  \textit{Step 1:} \textbf{Strong local convergence.}
  We are now going to prove convergence of $f_j$ in $C^1$ and $\varphi_j$ in $C^0$ in every bounded interval $I_R = [-R, R]$, with $R\geq 1$. First we have
  \begin{equation*}
    \LinfR{f_j} \leq \Linf{f_j-s|x|} + |s|R \lesssim_R |s|, 
  \end{equation*}
 where we bounded the $L^\infty$ by the $H^1$, which  is bounded by $|s|$. Also recall that $\Linf{f_j'} \lesssim |s|$. Now we go to $\varphi_j'$, by definition
 \begin{equation*}
    \LinfR{\varphi_j'} = \LinfR{(\partial_x\Psi_j+f_j'\partial_y\Psi_j)(\cdot,f_j(\cdot))} \lesssim (1+\delta)|s|R,
 \end{equation*}
 where we used the bound on $\nabla \Psi_j$, which also gives $\LinfR{G(f_j)\varphi_j}\lesssim_R |s|$ since $G(f_j)\varphi_j=(\partial_y\Psi_j-f_j'\partial_x\Psi_j)|_{y=f_j}$. Moreover,
 \begin{equation*}
  \begin{aligned}
    \LinfR{\varphi_j} 
    &\leq \frac12\LinfR{\int_{x-1}^{x+1} \varphi_j(x')dx' } + \frac12\LinfR{\int_{x-1}^{x+1} (\varphi_j(x)-\varphi_j(x')) dx'}\\ 
    &\lesssim \Ltwo{\varphi_j} + \|\varphi_j'\|_{L^\infty(I_{R+1})} \lesssim_R |s|,
  \end{aligned}
 \end{equation*}
 where we used Cauchy--Schwarz for the average and the derivative bound on $I_{R+1}$ for the difference. Finally, to bound the curvature, we use the Bernoulli equation \begin{equation*}
  \kappa_j = -\frac{1}{3}\varphi_j + \frac{2}{3} x\varphi_j' + N_j - \epsilon_j G(f_j)\varphi_j,
 \end{equation*} 
 to bound $f_j''$
  \begin{equation*}
  \begin{aligned}
    \LinfR{f_j''}
    &\leq (1+\delta^2)^{3/2} \LinfR{\kappa_j} \\
    &\lesssim \LinfR{\varphi_j} + R\LinfR{\varphi_j'} + \LinfR{N_j} + \epsilon_j\LinfR{G(f_j)\varphi_j}\\
    &\lesssim_R |s|(1+|s|),
  \end{aligned}
  \end{equation*}
  where we used that $\LinfR{N_j} \lesssim_R \LinfR{\varphi_j'}^2 + \LinfR{G(f_j)\varphi_j}^2$. Putting all together we have proved
  \begin{equation*}
    \| f_j\|_{W^{2, \infty}(I_R)} + \| \varphi_j\|_{W^{1, \infty}(I_R)} \lesssim_{R,s} 1,
  \end{equation*}
  Therefore, by Arzelà--Ascoli and a diagonal argument we have a subsequence that converges $f_j \to f$ in $C^1(I)$ and $\varphi_j \to \varphi$ in $C^0(I)$ for every compact interval $I \subset \R$, with $f\in W^{2,\infty}_{loc}$ and $\varphi\in W^{1,\infty}_{loc}$ by the same derivative bounds.
  Moreover, since we have the bounds 
  \begin{equation*}
    \Hone{f_j-s|x|} + \Honehalf{\varphi_j} \lesssim_s 1,
  \end{equation*}
  we have that a subsequence converges weakly in $H^1\times H^{1/2}$ to the same limit 
  \begin{equation*}
    f_j - s|x| \rightharpoonup f- s|x| \text{ in } H^1, \quad \text{ and } \quad \varphi_j \rightharpoonup \varphi \text{ in } H^{1/2}.
  \end{equation*}

  \textit{Step 2:} \textbf{Limit satisfies Kinematic equation.}
    Let us now define the flattened potential $\psi_j(x,z) := \Psi_j(x, z+f_j(x))$. Then
    \begin{equation*}
      \int_{\R} \varphi_j G(f_j)\varphi_j = \int_{z<0} A(f_j') \nabla \psi_j \cdot \nabla \psi_j, \quad \text{ where } \quad A(f_j') =
      \begin{pmatrix*}
        1 & -f_j' \\
        -f_j' & 1+f_j'^2
      \end{pmatrix*}
      .
    \end{equation*}
    The matrices $A(f_j')$ are uniformly elliptic and converge locally in $C^0$ to $A(f')$.
    Using the energy bound $\| \nabla \psi_j\|_{L^2(\R^2_-)}\lesssim \|\nabla \Psi_j\|_{L^2(\Omega_{f_j})} \lesssim |s|$, we get
    \begin{equation*}
      \int_{[-R, R]\times[-L,0]} |\psi_j|^2 \leq 2 L\int_{[-R, R]}  \varphi_j^2 dx + 2L^2 \int_{[-R,R]\times [-L,0]} |\p_z \psi_j|^2 \lesssim_L |s|^2,
    \end{equation*}  
    which gives boundedness of $\psi_j$ in $H^1$ on every such rectangle uniformly in $j$. Taking a diagonal subsequence we obtain that $\psi_j$ converges weakly there to $\psi$ with trace $\varphi$, and $\nabla\psi_j\rightharpoonup\nabla\psi$ in $L^2(\R^2_-)$. Testing the flattened harmonic equation against compactly supported zero-trace functions gives $\operatorname{div}(A(f')\nabla\psi)=0$, so $\Psi(x,y)=\psi(x,y-f(x))$ is the finite-energy harmonic extension of $\varphi$ by Lemma \ref{lem:H12_trace}. Changing variables back on interior compact sets and using interior harmonic estimates gives $\nabla\Psi_j\to\nabla\Psi$ locally uniformly in $\Omega_f$.
    Moreover, we have convergence of $G(f_j)\varphi_j$ in distributions. Indeed, for every $p \in C^\infty_c(\R)$
    \begin{equation*}
      \langle G(f_j)\varphi_j,p \rangle = \int_{\R^2_-} A(f_j') \nabla \psi_j \cdot \nabla E_p \, \longrightarrow  \, \int_{\R^2_-} A(f') \nabla \psi \cdot \nabla E_p = \langle G(f) \varphi,p \rangle,
    \end{equation*}
    where $E_p \in C^\infty_c(\overline{\R^2_-})$ is any extension of $p$. Therefore, using the regularized kinematic equation \eqref{eq:full_regularized}
    \begin{equation}\label{eq:kinematic_lim_eq}
      G(f_j)\varphi_j = -\frac{2}{3}(f_j-xf_j') + \epsilon_j(\kappa_j - N_j) \to G(f)\varphi =- \frac{2}{3}(f-xf') ,
    \end{equation}
    since $\epsilon_j \Ltwo{\kappa_j - N_j} \lesssim \sqrt{\epsilon_j}|s| \to 0$ by \eqref{eq:reg_energy_bounds_sec7} the RHS converges locally strongly in $L^2$, so also the LHS. We conclude the limiting functions $(f, \varphi)$ satisfy the kinematic equation.

  \textit{Step 3:} \textbf{Limit satisfies Dynamic equation.}
  For the dynamic equation, we are going to show the limit satisfies it in the distributional sense. For the LHS of \eqref{eq:full_regularized}, since $\varphi_j$ converges to $\varphi$ in $C^0$ on compact intervals, clearly
  \begin{equation*}
    -\frac13\varphi_j + \frac23 x \varphi_j' \to -\frac13\varphi + \frac23 x \varphi' \quad \text{ in }\quad \mathcal{D}'.
  \end{equation*}
  Then, the regularization part of the RHS $\epsilon_j G(f_j)\varphi_j$ goes to zero in $L^2$ \eqref{eq:reg_energy_bounds_sec7} so also in $D'$.
  Lastly, by \eqref{eq:energy_dynamic_stress} the other part of the RHS can be written as
  \begin{equation*}
    \partial_x \left ( \frac{f_j'}{\sqrt{1+f_j'^2}} - P_j \right ), \quad  \quad  P_j(x): = \int_{-\infty}^{f_j(x)} \partial_x \Psi_{j} \partial_y \Psi_j dy.
  \end{equation*}
  Since $f_j$ converges in $C^1$ locally, the first term converges in distributions. For $P_j$, write $d=f_j(x)-y>0$. By Lemmas \ref{lem:velocity_continuity} and \ref{lem:local_bulk} at unit scale and \eqref{eq:energy_bounds_sec7}--\eqref{eq:reg_energy_bounds_sec7}, $|\nabla\Psi_j(x,f_j(x)-d)|\lesssim |s|\langle x\rangle$ for $d<1$. For $d\geq1$ the mean-value estimate on an interior ball of radius $c d$ gives $|\nabla\Psi_j|\lesssim |s|/d$. At every positive depth the gradients converge by Step 2, so dominated convergence on $I_R\times(0,\infty)$ gives $P_j\to P=\int_{-\infty}^{f(x)}\Psi_x\Psi_y\,dy$ in $L^1(I_R)$. Therefore
  \begin{equation}\label{eq:bernoulli_lim_eq}
    \begin{aligned}
      -\frac{1}{3}\varphi_j + \frac{2}{3}x\varphi_j' &= \partial_x \left ( \frac{f_j'}{\sqrt{1+f_j'^2}} - P_j\right ) + \epsilon_j G(f_j)\varphi_j \\
      &\to \partial_x \left ( \frac{f'}{\sqrt{1+f'^2}} - P\right ) = -\frac{1}{3}\varphi + \frac{2}{3}x\varphi',
    \end{aligned}
  \end{equation}
  in the sense of distributions $D'$. 

  Due to \textit{Step 2} and \textit{Step 3} we have concluded that the limiting pair $(f, \varphi)$ is a weak solution to
  \begin{equation}\label{eq:weak_equations}
  \begin{aligned}
      -\frac{2}{3}(f-xf') &= G(f)\varphi \\
      - \frac{1}{3}\varphi + \frac{2}{3}x\varphi' &= \partial_x \left (\frac{f'}{\sqrt{1+f'^2}} - P \right ),
  \end{aligned} 
  \end{equation}
  in the distributional sense.

  \textit{Step 4:} \textbf{Smoothness bootstrap.}
  We prove that if $f \in C^{m, \alpha}_{loc}$ with $m\geq1$ and $\alpha \in (0,1)$, then if $(f, \varphi)$ are weak solutions to the equations \eqref{eq:weak_equations} it implies $\varphi \in C^{m, \alpha}_{loc}$, and $f \in C^{m+1, \alpha}_{loc}$. Therefore since we have that $f \in W^{2, \infty}_{loc}$, we can start the iteration and conclude $(f, \varphi)$ are locally $C^\infty$.

  We start by using the kinematic equation to conclude that
  \begin{equation}
    G(f)\varphi = -\frac{2}{3}(f- xf')  \in C^{m-1, \alpha}_{loc},
  \end{equation}
  so by Lemma \ref{lem:neumann_schauder} with $\partial_n\Psi=G(f)\varphi/\sqrt{1+f'^2}$ we have $\Psi \in C^{m, \alpha}_{loc}(\overline{\Omega_f})$ and therefore $\varphi \in C^{m, \alpha}_{loc}$. Then by differentiating $P$ in distributions, using $\partial_x(\Psi_x\Psi_y)=\frac12\partial_y(\Psi_x^2-\Psi_y^2)$ and the decay at large depth from Step 3, we get
  \begin{equation*}
    P' (x)= \left [ \frac{1}{2}(\Psi_x^2-\Psi_y^2) + f'\Psi_x \Psi_y\right ]_{y=f(x)} = N[f,\varphi] \in C^{m-1, \alpha}_{loc}.
  \end{equation*}
  Finally, writing the Bernoulli equation
  \begin{equation*}
    f'' = (1+f'^{2})^{3/2}(-\frac{1}{3}\varphi + \frac{2}{3}x\varphi' + P') \in C^{m-1, \alpha}_{loc},
  \end{equation*}
  so $f \in C^{m+1, \alpha}_{loc}$.

  \textit{Step 5:} \textbf{Global norm bounds.}
  Since $(f_j-s|x|, \varphi_j)$ converges weakly in $H^{1}\times H^{1/2}$, by weak lower semicontinuity we obtain
  \begin{equation*}
    \Hone{f-s|x|} + \Honehalf{\varphi} \lesssim |s|.
  \end{equation*}
  Since $f_j'\to f'$ locally uniformly and the velocity bounds of Step 3 pass to the limit at every positive depth, we get
  \begin{equation*}
    \Linf{f'} + \left\|\frac{\nabla \Psi}{\langle x \rangle }\right\|_{L^\infty(\Omega_f)} \leq C|s|.
  \end{equation*}
  Interior convergence and Fatou's lemma in \eqref{eq:reg_energy_bounds_sec7} give $\int_{\Omega_f}|\nabla\Psi|^2\lesssim s^2$. Since $f,\varphi$ are even and we have chosen $s_0$ such that $Cs_0 < \delta$, we have indeed $\Linf{f'}<\delta$ so $(f, \varphi)$ is admissible and under the hypotheses of Proposition \ref{prop:energy_balance} at $\epsilon=0$ and we can conclude $\mc{H}_s[f,\varphi] = 0$.
  
  \begin{remark}
    To conclude that $\mc{H}_s[f,\varphi]=0$ we could not use that $\epsilon_j \Ltwo{\nabla_{(f,\varphi)} \mc{H}_s[f_j,\varphi_j]}^2$ goes to zero since we don't know that a priori. To avoid that we wrote Section \ref{sec:energy} with the level of regularity of the unregularized solution.
  \end{remark}
\end{proof}

\subsection{Proof of corollaries}

\begin{proof}[Proof of Corollary \ref{cor:corner_formation}]
  Let $(f_s,\varphi_s)$ be a profile from Theorem \ref{thm:main_theorem}, with harmonic potential $\Psi_s$. The rescaling \eqref{eq:self_similar_coord} gives smooth data $(f_s,\varphi_s)$ at $t=-1$, and for $-1\leq t<0$ it gives
  \[
    \|\eta-s|\cdot|\|_{H^1}+\|\phi\|_{H^{1/2}}+\|u\|_{L^2(\Omega(t))}
    \lesssim |s|\,|t|^{1/3}.
  \]
  Thus the surface converges uniformly to the wedge and the velocity tends to zero in $L^2$. By evenness and the kinematic condition, $u(P_t,t)=(0,\partial_t\eta(0,t))$ and $\kappa(0,t)=\partial_X^2\eta(0,t)$, giving \eqref{eq:kappa_u_blowup_rate}. Also \eqref{eq:energy_high_frequency} gives $sf_s(0)>0$, since $f_s'$ is continuous and $s\sign(x)$ jumps at zero.

  To show that $sf_s''(0)>0$ for all sufficiently small nonzero $s$, we proceed by contradiction. If this fails, there are nonzero $s_n\to0$ and corresponding profiles with $f_{s_n}''(0)/s_n\leq0$ for every $n$. Along this sequence, write $s=s_n$ and set $\tilde{f}_s=f_s/s$, $\tilde{\varphi}_s=\varphi_s/s$. The energy bounds and the local bounds on $f_s''/s$ and $\varphi_s'/s$ from Step~1 of the proof of Theorem \ref{thm:main_theorem} give, along a subsequence,
  \[
    (\tilde{f}_s,\tilde{\varphi}_s)\to(\tilde{f},\tilde{\varphi})\ \text{in }C^1_{\rm loc}\times C^0_{\rm loc},\qquad
    (\tilde{f}_s-|\cdot|,\tilde{\varphi}_s)\rightharpoonup(\tilde{f}-|\cdot|,\tilde{\varphi})\ \text{in }H^1\times H^{1/2}.
  \]
  As in Step~2, the flattened potentials $\tilde{\psi}_s(x,z)=\Psi_s(x,z+f_s(x))/s$ converge weakly in $H^1$ on bounded rectangles up to $z=0$, with $\nabla \tilde{\psi}_s\rightharpoonup\nabla \tilde{\psi}$ in $L^2(\R^2_-)$. Since $A(f_s')\to I$ uniformly, $\tilde{\psi}$ is the harmonic extension of $\tilde{\varphi}$ and $G(f_s)\tilde{\varphi}_s\to |D|\tilde{\varphi}$ in distributions. In \eqref{eq:selfsimilar} the quadratic terms divided by $s$ are $O(|s|)$ on compact intervals, and the normalized curvature is $\partial_x(\tilde{f}_s'/\sqrt{1+f_s'^2})$. Thus in distributions
  \[
    \tfrac23(x\tilde{f}'-\tilde{f})=|D|\tilde{\varphi},\qquad \tfrac23 x\tilde{\varphi}'-\tfrac13\tilde{\varphi}=\tilde{f}'' .
  \]
  Rescaling $(\tilde{f},\tilde{\varphi})$ gives the flat linear evolution with terminal data $(|X|,0)$ by the same dilation estimate. At $t=-1$, the Fourier formula of Section \ref{subsec:strategy} and $\partial_t\widehat\eta=|\xi|\widehat\phi$ give
  \[
    \widehat{\tilde{\varphi}}(\xi)=-2\frac{\sin(|\xi|^{3/2})}{|\xi|^{3/2}},\qquad
    \tilde{\varphi}(0)=-\frac{4}{3\pi}\int_0^\infty\frac{\sin r}{r^{4/3}}\,dr<0 ,
  \]
  by Fourier inversion ($\widehat{\tilde{\varphi}}\in L^1$) and pairing successive half-periods in the integral.

  At $x=0$, evenness and the kinematic condition give $\varphi_s'(0)=f_s'(0)=0$ and $(G(f_s)\varphi_s)(0)=-2f_s(0)/3$. The dynamic equation at this point therefore gives
  \[
    \frac{f_s''(0)}{s}=-\frac13\tilde{\varphi}_s(0)-\frac29 s \tilde{f}_s(0)^2
    \longrightarrow -\frac13\tilde{\varphi}(0)>0 .
  \]
  This contradicts $f_{s_n}''(0)/s_n\leq0$, which holds along the subsequence as well. Thus, after taking $s_0$ smaller if necessary, $sf_s''(0)>0$ for every $0<|s|\leq s_0$. Since $P_t\to(0,0)$, \eqref{eq:kappa_u_blowup_rate} gives both lower bounds in every fixed neighborhood of the origin.
\end{proof}

\begin{proof}[Proof of Corollary \ref{cor:corner_retraction}]
    For $t>0$ take the solutions from the last Corollary \ref{cor:corner_formation} and reverse time $\eta(X, -t), -\Phi(X, Y,-t)$. Time reversibiliity of the Euler equations gives this smooth solution for $t>0$, which satisfies the same properties but when taking $t \downarrow 0 $ instead of $t \uparrow 0$.  
\end{proof}

\begin{proof}[Proof of Corollary \ref{cor:separation_two_fluids}]
  Let $s_j = \tan ((\theta_j-\pi)/2)$, and $f_j = f_{s_j}$ be self-similar profiles. We want to find upper and lower bounds on them such that we can stack one on top of the other (in opposite directions) and they do not have any intersection, even when adding a small enough tilt. Thus, returning to physical variables we will obtain that the first touching point is the corner at the origin.

  We proceed as in Lemma \ref{lem:hamiltonian_estimates}, using $\mc{H}_{s}[f, \varphi] = 0$ for either profile, but instead of doing the first-order expansion of $\sqrt{1+f'^2} - \sqrt{1+s^2}$ at $s\sign(x)$, we do it at $s\sign(x-x_0)$ for arbitrary $x_0$, and since the function $\sqrt{1+t^2}$ is strictly convex, we obtain the same conclusion
  \begin{equation*}
        \frac{2s f(x_0)}{\sqrt{1+s^2}}
        \geq \frac12\int_{\Omega_f}|\nabla\Psi|^2\,dx\,dy
        +c\int_{\R}|f'-s\sigma|^2\,dx>0,
  \end{equation*}
  for all $x_0 \in \R$. Hence, since we know that $f(x)-s|x|\to0$ as $|x|\to \infty$, and using the property we just derived on compact sets, there exists for each fixed profile with $s<0$ a constant $m(s)>0$, and for all admissible $s$ a uniform constant $C>0$ such that
  \begin{equation*}
    f(x) \leq -m(s) (1+|x|) \quad (s<0), \qquad f(x) \leq C (s)_+ (1+|x|) \quad \text{(all admissible $s$)},
  \end{equation*}
  where we used continuity of $f$ and $\Linf{f-s|x|}\lesssim |s|$ from Theorem \ref{thm:main_theorem}.

  Now consider a fixed $s_1<0$ and an arbitrary admissible $s_2$, with $0<|s_2|\leq s_0$. Then 
  \begin{equation*}
    -f_2(x)-f_1(x) \geq (m(s_1)-C(s_2)_+) (1+|x|),
  \end{equation*} 
  therefore, for every fixed $s_1<0$, if $(s_2)_+ < m(s_1)/C$ then the profiles have no intersection. 
  Moreover, if we do a tilt of angle $\delta$ around the origin of the upper profile, a sufficient condition for $-f_{2, \delta}(x) - f_1(x)>0$ becomes
  \begin{equation}
    (s_2)_+ + |\delta| < m(s_1)/C,
  \end{equation} 
  where this $C$ might be different, and we used that tilting around the origin affects the height by at most $C|\delta|(1+|x|)$ for small $\delta$.
  Converting these bounds to angles $\theta_1,\theta_2$ using $(s_2)_+\lesssim(\theta_2-\pi)_+$,  concludes the proof of the corollary, by time reversal for separation.
\end{proof}

\begin{proof}[Proof of Corollary \ref{cor:fluid_wall}]
  Restric the previous solutions to $X\geq 0$, taking $\delta=0$ at the Corollary \ref{cor:separation_two_fluids}. Since both the surface and the potential are even in $X$, we have that hte velocity $u= \nabla \Phi$ satisfies $u\cdot n =0$ at $X=0$. Hence this is a solution with non-penetration boundary conditions with a wall at $X=0$. Therefore, if we take the corner retraction soluiosn from $\theta$ to immediately smooth, in this case we get that we start at $\theta/2$ angle with the wall and immediately regularize to $90^\circ$. Time reversing, we obtain the opposite dynamical picture with the same mentioned angles.
\end{proof}

\appendix
\section{Boundary operators and elliptic estimates} 

\subsection{Estimates for harmonic extensions}\label{app:harmonic}

\begin{lemma}\label{lem:H12_trace}
  Let $f:\mathbb{R} \to \mathbb{R}$ be Lipschitz. For every $\varphi \in H^{1/2}(\mathbb{R})$, there exists a unique weakly harmonic function $\Psi$ in $\Omega_f = \left\{ (x,y) \in \mathbb{R}^{2} : y < f(x) \right\}$ such that $\nabla \Psi \in L^{2}(\Omega_f)$ and $\Psi(x, f(x)) = \varphi(x)$. Moreover,
  \begin{equation*}
   C(\| f' \|_{L^{\infty}})^{-1} \| \varphi \|_{\dot{H}^{1/2}} \leq \| \nabla \Psi \|_{L^2(\Omega_f)} \leq C(\| f' \|_{L^{\infty}}) \| \varphi \|_{\dot{H}^{1/2}},
  \end{equation*}
  where $C(\| f' \|_{L^{\infty}})$ denotes a constant depending on $\| f' \|_{L^{\infty}}$.
\end{lemma}

\begin{proof}
  This result is standard, and can for instance be found in \cite{NguyenPausader2020}. See \cite[Proposition 3.6]{NguyenPausader2020} for the existence and uniqueness of the harmonic extension, as well as the upper bound. The lower bound can be found in Theorem A.3 in the same paper.
\end{proof}

\begin{lem}\label{lem:local_bulk} 
  Let $f$ be Lipschitz, $p \in \partial \Omega_f$, and $\rho > 0$, and write $B_r = B_r(p)$. If $w \in C(\overline{\Omega}_f \cap B_{2\rho})$ is harmonic in $\Omega_f \cap B_{2\rho}$, then
  \begin{equation*}
    \|w\|_{L^\infty(\Omega_f \cap B_\rho)} \lesssim \|w\|_{L^\infty(\partial\Omega_f \cap B_{2\rho})} + \rho^{-1}\|w\|_{L^2(\Omega_f \cap B_{2\rho})}.
  \end{equation*}
  Moreover, if $w=0$ on $\partial\Omega_f \cap B_{2\rho}$, then
  \begin{equation*}
    \|w\|_{L^\infty(\Omega_f \cap B_\rho)} \lesssim \rho^{-2}\|w\|_{L^1(\Omega_f \cap B_{2\rho})}.
  \end{equation*}
\end{lem}

\begin{proof}
    Let $M = \|w\|_{L^{\infty}(\partial \Omega_f \cap B_{2\rho})}$. Consider the function $(|w| - M)_{+}$ which is extended by zero outside of $\Omega_f$. This function is continuous and subharmonic in $B_{2\rho}$. The mean value property for subharmonic functions gives
    \[(|w| - M)_+(z) \lesssim \rho^{-2} \int_{\Omega_f \cap B_{2\rho}}(|w| - M)_+ \, dx\,dy\]
    for any $z \in B_{\rho}$. The first assertion follows from Cauchy-Schwarz. The second assertion follows from the same argument, using that $M = 0$.
\end{proof}

\begin{lem}\label{lem:gradient_decay}
  Assume $f$ is Lipschitz with $\| f' \|_{L^{\infty}} \leq 1/2$. Fix $x_0 \in \mathbb{R}$ and $\rho>0$, and write $z_0 = x_0 + i f(x_0)$ and $B_r = B_r(z_0)$. If $w$ is harmonic in $\Omega_f \cap B_{4\rho}(z_0)$, continuous on $\overline{\Omega}_f \cap B_{4\rho}(z_0)$, and zero on $\partial\Omega_f \cap B_{4\rho}(z_0)$, then, for $0<d<\rho/4$,
\begin{equation*}
  \left| \nabla w(z_0 - id) \right| \lesssim d^{-1/3} \rho^{-8/3} \| w \|_{L^1(\Omega_f \cap B_{4\rho}(z_0))}. 
\end{equation*}
\end{lem}

\begin{proof}
  Without loss of generality, we may assume $z_0 = 0$. Write $M = \|w\|_{L^\infty(\Omega_f \cap B_{2\rho})}$. Consider the function $b(z)$ expressed in the polar coordinates (rotated so that $\theta = 0$ points downwards) as 
  \[b(z) = r^{2/3}\cos\left(\frac{2\theta}{3}\right).\]
  This function is positive and harmonic in the sector $y < |x|$, where $|\theta| < 3\pi/4$. Since $f(0)=0$ and $\|f'\|_{L^\infty} \leq 1/2$, we have $|f(x)|\leq |x|/2$, and hence
  \[\Omega_f \subset \{(x,y) : y < |x|/2\}, \qquad |\theta|\leq \frac{\pi}{2}+\arctan(1/2)<\frac{3\pi}{4} \quad \text{on } \Omega_f.\]
  It follows that there is some $c > 0$ such that, for every $z \in \Omega_f$,
  \[c|z|^{2/3} \leq b(z) \leq |z|^{2/3}.\]
  We apply the maximum principle to the subharmonic function
  \[|w| - \frac{M}{c(2\rho)^{2/3}}b\]
  on the region $\Omega_f \cap B_{2\rho}$. This function is continuous on the closure and nonpositive on the boundary: on the graph we have $w=0$, while on $\Omega_f\cap\partial B_{2\rho}$ we have $|w|\leq M$ and $b\geq c(2\rho)^{2/3}$. Therefore,
  \begin{equation}\label{eq:app_gradient_boundary_decay}
    |w(z)|\lesssim M\left(\frac{|z|}{\rho}\right)^{2/3}, \qquad z\in\Omega_f\cap B_{2\rho}.
  \end{equation}
  For $0<d<\rho/4$, the ball $B_{d/2}(-id)$ is contained in $\Omega_f\cap B_{3d/2}$. Indeed, in this ball $|x|<d/2$ and $y<-d/2$, whereas $f(x)\geq-|x|/2\geq-d/4$. The interior gradient estimate for harmonic functions and \eqref{eq:app_gradient_boundary_decay} give
  \begin{equation*}
    |\nabla w(-id)|\lesssim d^{-1}\|w\|_{L^\infty(B_{d/2}(-id))}
    \lesssim d^{-1}M\left(\frac{d}{\rho}\right)^{2/3}
    =d^{-1/3}\rho^{-2/3}M.
  \end{equation*}
  Finally, Lemma~\ref{lem:local_bulk} at scale $2\rho$ gives
  \[M\lesssim \rho^{-2}\|w\|_{L^1(\Omega_f\cap B_{4\rho})}.\]
  Combining these estimates proves the claim.
\end{proof}

\begin{lem}\label{lem:local_oblique}
  Assume $f$ is Lipschitz with $\|f'\|_{L^\infty}\leq 1/2$. Fix $x_0 \in \mathbb{R}$ and $R>0$, and write $z(x)=x+if(x)$, $B_r=B_r(z(x_0))$, and $I_{2R}=(x_0-2R,x_0+2R)$. Let $\theta,\gamma\in W^{1,\infty}(I_{2R})$. Assume $W$ is bounded and holomorphic in $\Omega_f\cap B_{2R}$, continuous on $\overline{\Omega_f}\cap B_{2R}$, and
  \begin{equation*}
    \text{Re}\big(e^{i\theta(x)}W(z(x))\big) = \gamma(x), \qquad z(x) \in B_{2R}.
  \end{equation*}
Then there is a constant $C_R$ depending only on $R$ such that
\begin{equation*}
  \| W \|_{C^{1/2}(\overline{\Omega_f} \cap \overline{B}_R)} \leq C_R\Big(\|\gamma'\|_{L^\infty(I_{2R})} + (1+\|\theta'\|_{L^\infty(I_{2R})})\|W\|_{L^\infty(\Omega_f\cap B_{2R})}\Big).
\end{equation*}
\end{lem} 

\begin{proof}
  Let $M = \|W\|_{L^\infty(\Omega_f \cap B_{2R})}$ and 
  \[K = \|\gamma'\|_{L^\infty(I_{2R})} + M\|\theta'\|_{L^\infty(I_{2R})}.\]
  Fix a point $P=z(a)\in\partial\Omega_f\cap B_{3R/2}$ and set
  \[u_P(z)=\operatorname{Re}\big(e^{i\theta(a)}W(z)\big)-\gamma(a).\]
  This function is harmonic, satisfies $|\nabla u_P|=|W'|$, and is bounded by $2M$, since $|\gamma(a)|\leq M$. Note that we can express
  \[
    u_P(z(x)) = \gamma(x)-\gamma(a) + \operatorname{Re}\Big(\big(e^{i\theta(a)}-e^{i\theta(x)}\big)W(z(x))\Big).  
  \]
  Therefore, by definition of $K$,
  \[|u_P(z(x))| \leq K|x-a|\leq K|z(x)-P|\]
  Let $b_P$ be the harmonic function of \Cref{lem:gradient_decay}, centered at $P$. It is positive and harmonic in $\Omega_f$, with
  \[c|z-P|^{2/3}\leq b_P(z)\leq |z-P|^{2/3}\]
  for some fixed $c>0$. Since $B_{R/4}(P)\subset B_{2R}$, we can apply the comparison principle to the subharmonic function
  \[|u_P| - C\big(MR^{-2/3}+KR^{1/3}\big)b_P\]
  on $\Omega_f\cap B_{R/4}(P)$ for $C$ large enough. Indeed, on $\partial \Omega_f\cap B_{R/4}(P)$, we have that $|u_P(z)| \leq K|z - P|$ and on the rest of the boundary of $ \Omega_f\cap B_{R/4}(P)$, we have $|u_P|\leq 2M$. The maximum principle therefore gives
  \begin{equation}\label{eq:app_oblique_boundary_decay}
    |u_P(z)|\lesssim \big(MR^{-2/3}+KR^{1/3}\big)|z-P|^{2/3},
    \qquad z\in\Omega_f\cap B_{R/4}(P).
  \end{equation}
  Now let $z\in\Omega_f\cap B_{5R/4}$ and put $\delta=\operatorname{dist}(z,\partial\Omega_f)$. If $\delta<R/16$, the nearest graph point $P$ on $\partial \Omega_f$ is in $B_{21R/16}\subset B_{3R/2}$, and $B_{\delta/2}(z)\subset\Omega_f\cap B_{R/4}(P)$. Therefore we can apply the interior gradient estimate to $u_P$ on $B_{\delta/2}(z)$ and use \eqref{eq:app_oblique_boundary_decay} to obtain
  \begin{equation}\label{eq:app_oblique_gradient_decay}
    |W'(z)|=|\nabla u_P(z)|
    \lesssim \big(MR^{-2/3}+KR^{1/3}\big)\delta^{-1/3}.
  \end{equation}
  If $\delta\geq R/16$, the interior gradient estimate on $B_{R/32}(z)$ gives $|W'(z)|\lesssim M/R$. Since $\delta\leq5R/4$, this also implies \eqref{eq:app_oblique_gradient_decay}. \par
  Finally, we will obtain the desired H\"older estimate by integrating \eqref{eq:app_oblique_gradient_decay}. Let $z_1, z_2 \in \overline{\Omega_f}\cap\overline B_R$. First suppose that $|z_1 - z_2| < R / 16$. Then we can connect $z_1$ and $z_2$ by a path consisting of three line segments: a vertical segment from $z_1$ downwards by $2|z_1 - z_2|$, a vertical segment from $z_2$ downwards by $2|z_1 - z_2|$, and a segment connecting the two endpoints of these segments. This path remains in $\overline{\Omega_f} \cap B_{9R/8}$ because $\|f'\|_{L^\infty} \leq 1/2$. Integrating \eqref{eq:app_oblique_gradient_decay} along these paths gives
  \begin{equation*}
    \begin{aligned}
      |W(z_1)-W(z_2)|
      &\lesssim \big(MR^{-2/3}+KR^{1/3}\big)
      \left(2\int_0^{2|z_1-z_2|}t^{-1/3}\,dt+|z_1 - z_2|\,|z_1 - z_2|^{-1/3}\right)\\
      &\lesssim \big(MR^{-2/3}+KR^{1/3}\big)|z_1 - z_2|^{2/3}.
    \end{aligned}
  \end{equation*}
  If $|z_1 - z_2|\geq R/16$, the bound $|W(z_1)-W(z_2)|\leq2M$ gives the same estimate. Therefore,
  \begin{equation}\label{eq:app_oblique_holder}
    [W]_{C^{2/3}(\overline{\Omega_f}\cap\overline B_R)}
    \lesssim MR^{-2/3}+KR^{1/3}.
  \end{equation}
  Since the diameter of this set is at most $2R$, \eqref{eq:app_oblique_holder} implies
  \[\|W\|_{C^{1/2}(\overline{\Omega_f}\cap\overline B_R)}
    \lesssim M+R^{-1/2}M+R^{1/2}K\leq C_R(M+K),\]
  as wanted.
\end{proof}

\begin{lem}\label{lem:neumann_schauder}
Let $m\geq1$ be an integer and $0<\alpha<1$. Suppose $\Psi\in H^1_{\mathrm{loc}}(\overline\Omega_f)$ is weakly harmonic with Neumann data $g$:
\[\begin{cases}
  \Delta\Psi=0, &\text{in }\Omega_f,\\
  \partial_n\Psi=g, &\text{on }\partial\Omega_f,
\end{cases}.\]
If $f\in C^{m,\alpha}_{\mathrm{loc}}(\R)$ and $g\in C^{m-1,\alpha}_{\mathrm{loc}}(\R)$, then $\Psi \in C^{m,\alpha}_{\mathrm{loc}}(\ov{\Omega_f})$. 
\end{lem}

\begin{proof}
  Let $\Theta$ denote the harmonic conjugate of $\Psi$ in $\Omega_f$ and write $\theta(x) = \Theta(x, f(x))$. Then
  \[\theta' = \Theta_x + f' \Theta_y = f' \Psi_x - \Psi_y = -\sqrt{1 + f'^2} \, \partial_n \Psi \in C_{\mathrm{loc}}^{m-1,\alpha}(\R).\] It follows that $\Theta$ has $C_{\mathrm{loc}}^{m,\alpha}$ Dirichlet data. By \cite[Corollary 8.36]{GT}, we have that $\Theta \in C^{1, \alpha}_{\mathrm{loc}}(\ov{\Omega_f})$. Then by \cite[Corollary 6.19]{GT} for $m \geq 2$, we can upgrade the regularity to $\Theta \in C^{m,\alpha}_{\mathrm{loc}}(\ov{\Omega_f})$. The Cauchy-Riemann equations then give the same regularity for $\Psi$ as wanted.
\end{proof}

\subsection{Dirichlet-to-Neumann estimates}\label{app:DtN}

This section contains some standard estimates for the Dirichlet-to-Neumann operator adapted to our setting of Lipschitz graphs.

\begin{lemma}\label{lem:DtN_H12_bound} 
  Let $f: \mathbb{R} \to \mathbb{R}$ be Lipschitz. For $\varphi \in H^{1/2}$, the Dirichlet-to-Neumann operator satisfies $G(f)\varphi\in\dot H^{-1/2}(\R)$ and there exists a constant $C(\| f' \|_{L^{\infty}})$ depending on $\| f' \|_{L^{\infty}}$ such that
  \begin{equation*}
    \|G(f)\varphi\|_{\dot H^{-1/2}} \leq C(\|f'\|_{L^\infty})\|\varphi\|_{\dot H^{1/2}},
  \end{equation*}
\end{lemma}
\begin{proof}
  We use the estimate in \cite[Proposition 3.6]{NguyenPausader2020}, whose space for boundary data  is $\dot H^{1/2}$ by Remark 3.1 of the same paper. This result applies to arbitrary Lipschitz graphs. Let $\Psi$ and $U$ be the finite-energy harmonic extensions of $\varphi$ and $\eta\in C_c^\infty(\R)$, respectively. Using Green's identity and the cited estimate in the last step,
  \begin{equation*}
    |\langle G(f)\varphi,\eta\rangle| = \left|\int_{\Omega_f}\nabla\Psi\cdot\nabla U\,dx\,dy\right| \leq \|\nabla\Psi\|_{L^2(\Omega_f)}\|\nabla U\|_{L^2(\Omega_f)} \leq C(\|f'\|_{L^\infty}) \|\varphi\|_{\dot H^{1/2}}\|\eta\|_{\dot H^{1/2}}.
  \end{equation*}
  Taking the supremum over $\eta$ proves the claim by duality.
\end{proof}

\begin{lemma}\label{lem:Gf_Lp}
  For every $1<p<\infty$ there exists $\delta_p > 0$ such that, when $f$ is Lipschitz with $\| f' \|_{L^{\infty}} \leq \delta_p$ and $\varphi \in H^{1/2}$, we have $\varphi' \in L^{p}(\mathbb{R})$ if and only if $G(f)\varphi \in L^{p}(\mathbb{R})$. In that case,
  \begin{equation*}
    \| \varphi' \|_{L^{p}} \lesssim_{p} \| G(f)\varphi\|_{L^p} \lesssim_{p} \|\varphi'\|_{L^p}.
  \end{equation*}
\end{lemma}
\begin{proof}
  We apply~\cite[Theorem 2.2]{DindosPipher2024} to the Laplacian in $\Omega_f$, using its unbounded graph case $r_0=\infty$. The coefficient Carleson measure vanishes, so sufficiently small $\|f'\|_{L^\infty}$ gives both the $L^p$ regularity and Neumann estimates. Since the nontangential maximal function controls both derivatives at the boundary, the finite-energy harmonic extension $\Psi$ of $\varphi$ satisfies
  \begin{equation*}
    \|\partial_n\Psi\|_{L^p(d\sigma)} \lesssim_p \|\partial_\tau\Psi\|_{L^p(d\sigma)} \lesssim_p \|\partial_n\Psi\|_{L^p(d\sigma)}.
  \end{equation*}
  Here $n$ is the unit outward normal and $\tau$ the unit tangent in the direction increasing $x$. Since in our case $d\sigma=\sqrt{1+f'^2}dx$, $\partial_\tau\Psi=\varphi'/\sqrt{1+f'^2}$ and $\partial_n\Psi=G(f)\varphi/\sqrt{1+f'^2}$, the bound $1\leq\sqrt{1+f'^2}\leq\sqrt{1+\delta_p^2}$ gives the desired estimates.
\end{proof}

The following argument is a standard description of the Dirichlet-to-Neumann operator through the Hilbert transform; a precise reference is, for instance, \cite{Wu:almost-global-wellposedness-2d}. If $\varphi$ is the trace of the velocity potential $\Psi$ on the boundary of the fluid domain $\Omega_f$, then the Dirichlet-to-Neumann operator is
\begin{equation*}
  G(f)\varphi = \bigl(\Psi_y - f' \Psi_x\bigr)(x, f(x)).
\end{equation*}
The function $\Psi$ has a harmonic conjugate $\Theta$ in $\Omega_f$. Define 
\begin{equation*}
  z = x + iy, \qquad F(z) = \Psi(x,y) + i \Theta(x,y), \qquad \Theta_b(x) = \Theta(x, f(x)), \qquad w(x) = \varphi(x) + i \Theta_b(x).
\end{equation*}
For $f$ Lipschitz and $\varphi \in H^{1/2}$, the harmonic conjugate can be normalized so that $\Theta_b \in L^{2}$, and
\begin{equation}\label{eq:W_trace_varphi}
 \| w \|_{L^{2}} \lesssim_{\| f' \|_{L^{\infty}}} \| \varphi \|_{L^{2}}.    
\end{equation}
The Cauchy-Riemann equations give
\begin{equation}\label{eq:DtN_harmonic_conjugate}
  G(f)\varphi = -(\Theta_x + f' \Theta_y)(x,f(x)) = -\Theta_b'(x).
\end{equation}
Now, using the Cauchy integral formula, we have
\begin{equation}\label{eq:DtN_Cauchy_interior}
  F(z) = \frac{1}{2\pi i} \int_{\mathbb{R}} \frac{1 + if'(t)}{(x+iy) - (t + i f(t))} w(t) \; d t,
\end{equation}
and taking the limit in $\Omega_f$ as $z \to x + i f(x)$, we obtain, as an identity in $L^{2}$, 
\begin{equation}\label{eq:DtN_Cauchy_boundary}
  w(x) = \frac{1}{2} w(x) + \frac{1}{2\pi i} \text{p.v.} \int_{\mathbb{R}} \frac{1 + if'(t)}{(x-t) + i(f(x) - f(t))} w(t) d t.
\end{equation}
Recalling that $w(x) = \varphi(x) + i \Theta_b(x)$ and comparing the real and imaginary parts in this formula, we can get an expression for $\Theta_b(x)$. Using \eqref{eq:DtN_harmonic_conjugate} gives the following integral formula for the Dirichlet-to-Neumann operator
\begin{equation}\label{eq:DtN_integral_formula}
  G(f)\varphi(x) = \Re \partial_x \left[ \frac{1}{\pi} \operatorname{p.v.} \int_{\mathbb{R}} \frac{1 + if'(t)}{(x - t) + i(f(x) - f(t))} (\varphi(t) + i\Theta_b(t)) \; d t  \right].
\end{equation}
It is important to note that in the case of a flat interface, $f \equiv 0$, the previous formula reduces to the Hilbert transform through
\begin{equation}\label{eq:half_laplacian_integral_formula}
  G(0)\varphi(x) = \partial_x \left[ \frac{1}{\pi} \operatorname{p.v.} \int_{\mathbb{R}} \frac{\varphi(y)}{x - y} \; d y  \right] = H(\varphi')(x).
\end{equation}
In Fourier variables, this is the multiplier $|\xi|$, which will be denoted by $\Lambda = H \partial_x$.

\begin{lemma}\label{lem:DtN_analytic}
  There exists $\delta>0$ such that the map $f'\mapsto G(f)$ is real analytic from the set $\{f'\in W^{1,\infty}(\R):\|f'\|_{L^\infty}<\delta\}$ into $\mathcal{L}(H^2,H^1)$. Moreover, in this domain,
  \begin{equation}\label{eq:DtN_f_dependance_estimate}
    \| (G(f_1) - G(f_2)) \varphi \|_{L^{2}} \lesssim_{\delta} \| f_1' - f_2' \|_{L^{\infty}} \| \varphi' \|_{L^{2}}, \qquad \text{ for } \varphi \in H^{1}.
  \end{equation}
\end{lemma}

\begin{proof}
  This proof adapts the argument in \cite{CoifmanMeyer1985} to our setting. Differentiating \eqref{eq:DtN_Cauchy_boundary} and integrating by parts in $t$ gives $\mathcal C_f w'=iw'$, where
  \begin{equation*}
    \mathcal C_f u(x) = \frac{1+if'(x)}{\pi}\operatorname{p.v.} \int_{\R}\frac{u(t)}{x-t + i(f(x)-f(t))}\,dt = A_f u(x) + iB_f u(x),
  \end{equation*}
  Since $w'=\varphi'-iG(f)\varphi$ by \eqref{eq:DtN_harmonic_conjugate}, taking real parts gives
  \begin{equation*}
    (I-B_f)G(f)\varphi=A_f\varphi'.
  \end{equation*}
  Now \cite{CoifmanMeyer1985} gives analytic dependence of the Cauchy integral operator on $f'$ in $\mathcal L(L^2,L^2)$. The same result is true for $\mathcal C_f$ since this operator just differs by multiplication by $1+if'$ and its inverse. Since $B_0=0$, a Neumann series argument makes $(I-B_f)^{-1}A_f$ analytic for small $\|f'\|_{L^\infty}$. Its Taylor expansion therefore yields
  \begin{equation*}
    G(f)\varphi=\sum_{n\geq0}T_n(f',\ldots,f')\varphi'\,,
  \end{equation*}
  whose symmetric multilinear coefficients, with $T_0=H$, satisfy (cf. \cite[p.~75, (1.4)]{CoifmanMeyer1985})
  \begin{equation*}
    \|T_n(f_1',\ldots,f_n')u\|_{L^2} \leq C^{n+1}\prod_{j=1}^n\|f_j'\|_{L^\infty}\|u\|_{L^2}\,.
  \end{equation*}
  Thus, for $\|f_i'\|_{L^\infty}\leq\delta$, we have
\begin{equation*}
  \|(G(f_1)-G(f_2))\varphi\|_{L^2} \leq \sum_{n\geq1} nC^{n+1}\delta^{n-1} \|f_1'-f_2'\|_{L^\infty}\|\varphi'\|_{L^2} \lesssim_\delta \|f_1'-f_2'\|_{L^\infty}\|\varphi'\|_{L^2},
\end{equation*}
  provided $C\delta<1$, and this proves the estimate \eqref{eq:DtN_f_dependance_estimate}. Moreover,
  \begin{equation*}
    \partial_x\bigl(T_n(f',\ldots,f')u\bigr) = T_n(f',\ldots,f')u' + \sum_{j=1}^n T_n(f',\ldots,\underset{j}{f''},\ldots,f')u,
  \end{equation*}
  where in the $j$-th summand the derivative falls on the $j$-th entry. Applying the $L^2$ bound to each term gives
  \begin{align*}
    \|T_n(f',\ldots,f')u\|_{H^1} \leq C^{n+1}\Big(\|f'\|_{L^\infty}^n\|u\|_{H^1} + n\|f'\|_{L^\infty}^{n-1} \|f''\|_{L^\infty}\|u\|_{L^2}\Big).
  \end{align*}
  The bound is summable locally uniformly in $W^{1,\infty}$ whenever $C\|f'\|_{L^\infty}<1$. Thus the series converges locally uniformly in $\mathcal L(H^1,H^1)$ and is analytic. Composing with $\partial_x:H^2\to H^1$ proves the claim for $\delta<C^{-1}$.
\end{proof}

\begin{lemma}\label{lem:DtN_commutator}
  There exists $\delta>0$ such that, for every $M>0$ and $f$ Lipschitz with $f' \in W^{1, \infty}$, $\| f' \|_{L^{\infty}} \leq \delta$, $\| f'' \|_{L^{\infty}} \leq M$, $h \in W^{2, \infty}$ and $\varphi \in H^{1}$, we have     
  \begin{equation*}
    \| \left[ G(f), h \right] \varphi \|_{H^{1}} \lesssim_{\delta, M} (\| h' \|_{L^{\infty}} + \| h'' \|_{L^{\infty}}) \| \varphi \|_{H^{1}}. 
  \end{equation*} 
\end{lemma}
\begin{proof}
We map the fluid domain $\Omega_f$ to the upper half-plane. In these variables, the harmonic extension $\Psi$ of $\varphi$ is given by $U(x,t) = \Psi(x, f(x) - t)$. This function satisfies the elliptic equation
\begin{equation*}
\left\{
  \begin{alignedat}{2}
  - \operatorname{div}(A_f(x) \nabla U) &= 0, &\qquad &\text{ in } \mathbb{R}^{2}_+ \\
  U(x,0) &= \varphi(x), &\qquad &\text{ on } \mathbb{R},
\end{alignedat}
\right.
\qquad \text{ with } \quad
  A_f(x) = \begin{pmatrix}
    1 & f'(x) \\
    f'(x) & 1 + f'(x)^2
  \end{pmatrix}.
\end{equation*}
This matrix is real, symmetric and independent of the vertical variable. Since $\| f' \|_{L^{\infty}} \leq \delta$, this matrix is uniformly elliptic and bounded. Thus, we can apply \cite[Theorem 1.1]{HofmannZhang2021} and Lemma~\ref{lem:Gf_Lp} to get both
\begin{equation}\label{eq:DtN_commutator_L2}
  \|G(f)\varphi\|_{L^2} \lesssim \|\varphi'\|_{L^2}, \qquad \|[G(f),h]\varphi\|_{L^2} \lesssim \|h'\|_{L^\infty}\|\varphi\|_{L^2}\,,
\end{equation}
uniformly for sufficiently small slopes.

Now, we want to extend the estimate to $H^1$. In the flat case $G(0) = \Lambda$, and the commutator estimate is classical. We are going to use that $f'$ is small to extend the argument using an absorption argument. Fix $f$ and consider the horizontal and vertical velocity traces of the harmonic extension $\Psi$ of $\varphi$
\begin{equation*}
  V(\varphi)=\Psi_x(x,f(x)), \qquad B(\varphi)=\Psi_y(x,f(x)).
\end{equation*}
Since we have the identities
\begin{equation}\label{eq:app_boundary_identities}
  \varphi'=V(\varphi)+f'B(\varphi), \qquad G(f)\varphi=B(\varphi)-f'V(\varphi),
\end{equation}
the fact that $\Psi$ is harmonic gives 
\begin{equation}\label{eq:DtN_horizontal_velocity}
  B(\varphi)'=\Psi_{xy}+f'\Psi_{yy}
  =\Psi_{xy}-f'\Psi_{xx}=G(f)V(\varphi).
\end{equation}
Taking commutators in \eqref{eq:app_boundary_identities} yields
\begin{equation}\label{eq:horizontal_velocity_commutator}
\begin{aligned}
  [V,h]\varphi+f'[B,h]\varphi &= h'\varphi,\\
  [B,h]\varphi-f'[V,h]\varphi &= [G(f),h]\varphi.
\end{aligned}
\end{equation}
This system has determinant $1+f'^2$, so the $L^2$ commutator bound \eqref{eq:DtN_commutator_L2} implies
\begin{equation*}
  \Ltwo{[V,h]\varphi}+\Ltwo{[B,h]\varphi} \lesssim \Ltwo{h'\varphi}+\Ltwo{[G(f),h]\varphi} \lesssim \Linf{h'}\Ltwo{\varphi}.
\end{equation*}

Next, we differentiate the second identity in \eqref{eq:horizontal_velocity_commutator}. Consider for now all functions to be smooth. Then we get
\begin{equation}\label{eq:DtN_commutator_derivative}
  \partial_x([G(f),h]\varphi) =\partial_x([B,h]\varphi)-f'\partial_x([V,h]\varphi)-f''[V,h]\varphi.
\end{equation}
We estimate each of these terms now. Differentiating $[B,h]\varphi$ and using \eqref{eq:DtN_horizontal_velocity} gives
\begin{equation*}
  \partial_x([B,h]\varphi)
  =G(f)([V,h]\varphi)+[G(f),h]V(\varphi)-h'B(\varphi).
\end{equation*}
Since $\Ltwo{V(\varphi)}+\Ltwo{B(\varphi)}\lesssim\Hone{\varphi}$, the two estimates in~\eqref{eq:DtN_commutator_L2} yield
\begin{equation*}
  \Ltwo{\partial_x([B,h]\varphi)}
  \leq C_*\Ltwo{\partial_x([V,h]\varphi)}
  +C_*\Linf{h'}\Hone{\varphi},
\end{equation*}
where $C_*$ is uniform for sufficiently small slopes and independent of $M$. Differentiating the first identity in \eqref{eq:horizontal_velocity_commutator} gives
\begin{equation*}
  \partial_x([V,h]\varphi) = (h'\varphi)'-f'\partial_x([B,h]\varphi)-f''[B,h]\varphi,
\end{equation*}
and hence
\begin{equation*}
  \Ltwo{\partial_x([V,h]\varphi)} \leq \delta\Ltwo{\partial_x([B,h]\varphi)} + C_{\delta,M}\|h'\|_{W^{1,\infty}}\Hone{\varphi}.
\end{equation*}
Choose $\delta$ so that $C_*\delta\leq 1/2$. Substitution into the preceding estimate and absorption then give
\begin{equation}\label{eq:velocity_commutator_derivative_bound}
  \Ltwo{\partial_x([V,h]\varphi)}+\Ltwo{\partial_x([B,h]\varphi)} \lesssim_{\delta,M} \|h'\|_{W^{1,\infty}}\Hone{\varphi}.
\end{equation}
Finally substituting \eqref{eq:velocity_commutator_derivative_bound} into \eqref{eq:DtN_commutator_derivative} and using the $L^2$ bound for $[V,h]\varphi$ gives the desired estimate for smooth functions.

The absorption step above requires derivatives of commutators to be in $L^2$, which was not justified. To make this rigorous, mollify $f$ and $h$ and approximate $\varphi$ in $H^1$ by smooth compactly supported functions. The estimate \eqref{eq:DtN_f_dependance_estimate}, together with Lemma \ref{lem:Gf_Lp} for $p=2$, gives convergence of these commutators to $[G(f),h]\varphi$ in $L^2$. Weak lower semicontinuity in $H^1$ gives the desired estimate.
\end{proof}

\bibliographystyle{amsalpha}
\bibliography{ref}

\end{document}